%% file: article.tex
\documentclass[graybox]{svmult}

\usepackage{mathptmx}       
\usepackage{helvet}         
\usepackage{courier}        
\usepackage{type1cm}        
\usepackage{makeidx}         
\usepackage{graphicx}        
\usepackage{multicol}        
\usepackage[bottom]{footmisc}
\usepackage{amsmath,amsfonts,latexsym}
\usepackage{enumerate}
\usepackage{cite}
\usepackage{comment}

\makeindex             

\def\bi{{\mathbf i}}
\def\BbbR{{\mathbb R}}

\def\BbbN{{\mathbb N}}

\def\EI{{\mathcal E}_I}
\def\EB{{\mathcal E}_B}
\newcommand{\bsigma}{\boldsymbol \sigma}
\newcommand{\bSigma}{\boldsymbol \Sigma}
\newcommand{\btau}{\boldsymbol \tau}
\newcommand{\jump}[1]{[\![ #1 ]\!]}
\newcommand{\average}[1]{\{\!\!\!\{ #1 \}\!\!\!\}}
\newenvironment{numberedproof}[1]{\emph{Proof of #1.} }{\ \qed}
\DeclareMathSymbol \lesssim {\mathord}{AMSa}{"2E}
\newcommand{\eremk}{\hbox{}\hfill\rule{0.8ex}{0.8ex}}
\spnewtheorem{Theorem}{Theorem}[section]{\bfseries}{\itshape}
\spnewtheorem{Lemma}[Theorem]{Lemma}{\bfseries}{\itshape}
\spnewtheorem{Corollary}[Theorem]{Corollary}{\bfseries}{\itshape}
\spnewtheorem{Example}[Theorem]{Example}{\bfseries}{\itshape}
\spnewtheorem{Remark}[Theorem]{Remark}{\bfseries}{\itshape}
\spnewtheorem{Proposition}[Theorem]{Proposition}{\bfseries}{\itshape}
\newif\iftechreport
\techreporttrue    
\iftechreport
\else
\fi

\begin{document}

\iftechreport
\title*{On stability of discretizations of the Helmholtz equation (extended version)}
\fi 
\author{S.~Esterhazy and J.M. Melenk}
\institute{Vienna University of Technology, Institute for Analysis and Scientific Computing. Wiedner Hauptstrasse 8-10, A-1040 Vienna. \email{s.esterhazy@tuwien.ac.at, melenk@tuwien.ac.at}}
%
%
\maketitle

\abstract{ 
We review the stability properties of several discretizations of the Helmholtz equation at 
large wavenumbers. For a model problem in a polygon, a complete $k$-explicit stability 
(including $k$-explicit stability of the continuous problem) and
convergence theory for high order finite element methods is developed. 
In particular, 
quasi-optimality is shown for a fixed number of degrees of freedom per wavelength 
if the mesh size $h$ and the approximation order $p$ are selected such that 
$kh/p$ is sufficiently small and $p = O(\log k)$, and, additionally, 
appropriate mesh refinement is used near the vertices. We also review the stability 
properties of two classes of numerical schemes that use piecewise solutions of the 
homogeneous Helmholtz equation, namely, Least Squares methods and 
Discontinuous Galerkin (DG) methods. The latter includes the Ultra Weak Variational Formulation. 
}
\section{Introduction}
%
A fundamental equation describing acoustic or electromagnetic 
phenomena is the time-dependent wave equation 
$$ 
\frac{\partial^2 w}{\partial t^2} -  c^2 \Delta w = g, 
$$
given here for homogeneous, isotropic media whose propagation speed
of waves is $c$. It arises in many applications, for example,  
radar/sonar detection, noise filtering, optical fiber design, 
medical imaging and seismic analysis. A commonly encountered 
setting is the time-harmonic case, in which 
the solution $w$ (and correspondingly the right-hand side $g$) is assumed 
to be of the form 
$\operatorname*{Re} \left(e^{- \bi \omega t} u(x)\right)$ for a frequency
$\omega$. Upon introducing the {\em wavenumber} $k = \omega/c$
and the {\em wave length} $\lambda:= 2\pi/k$, the resulting equation for 
the function $u$, which depends solely on the 
spatial variable $x$, is then the {\em Helmholtz equation}
\begin{equation}
\label{eq:helmholtz}
-\Delta u -  k^2 u = f.  
\end{equation}
\iftechreport 
In many high frequency situations of large $k$ 
the solution $u$ is highly oscillatory but has some 
multiscale character that can be captured, for example,
by means of asymptotic analysis; a classical reference 
in this direction is \cite{babich-buldyrev91}. 
\fi

\begin{figure}[tb]
\includegraphics[width=0.5\textwidth]{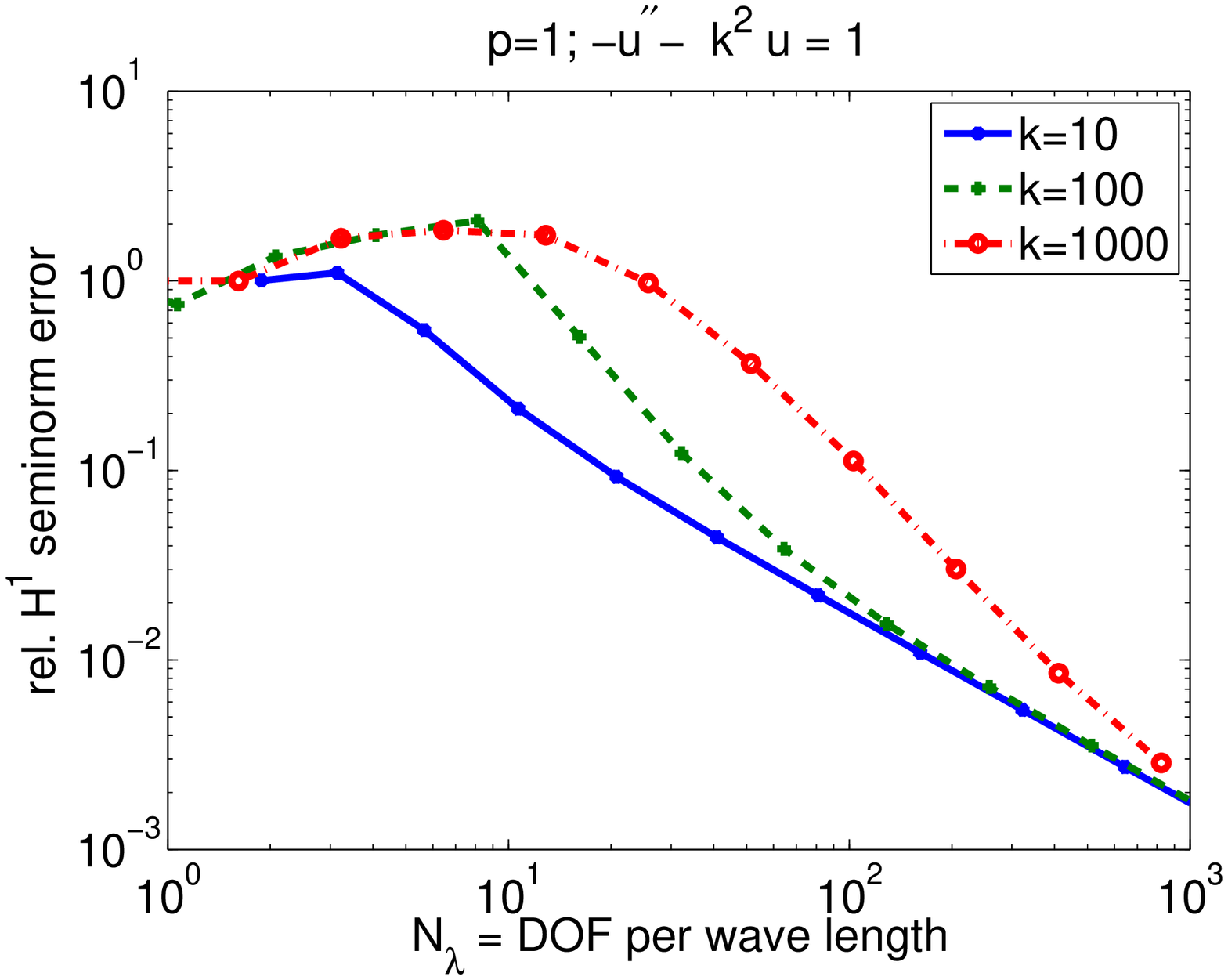}
\includegraphics[width=0.5\textwidth]{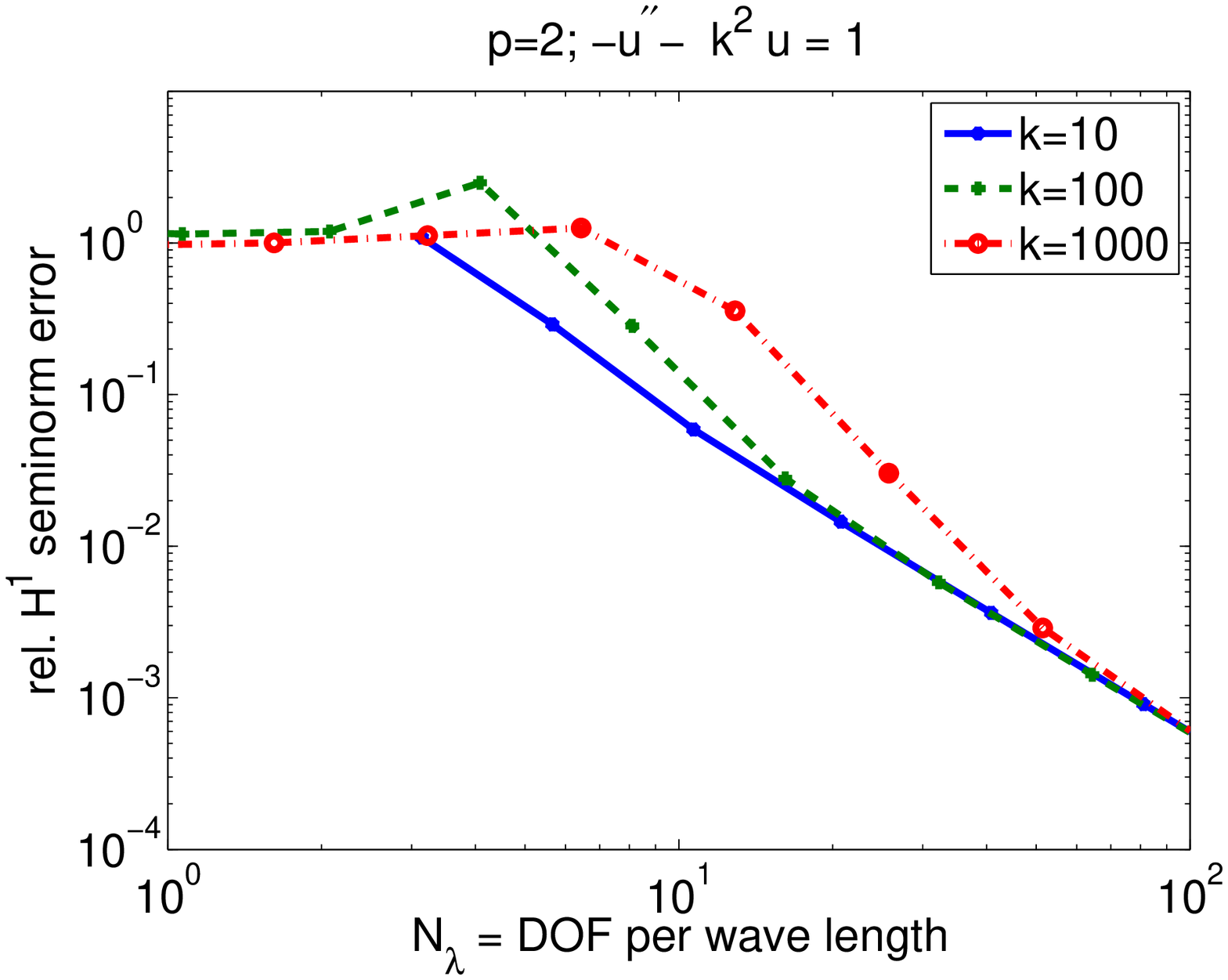}
\includegraphics[width=0.5\textwidth]{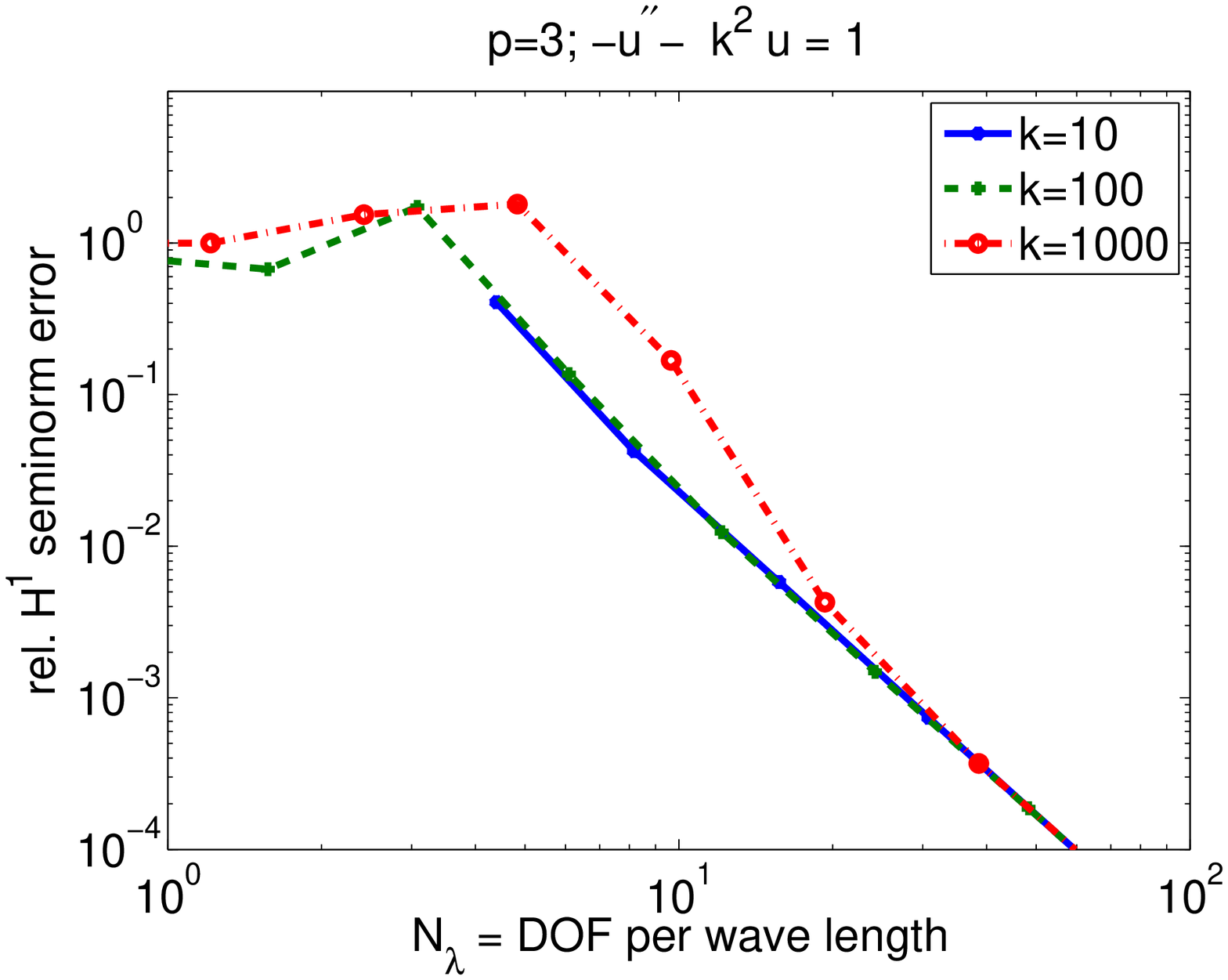}
\includegraphics[width=0.5\textwidth]{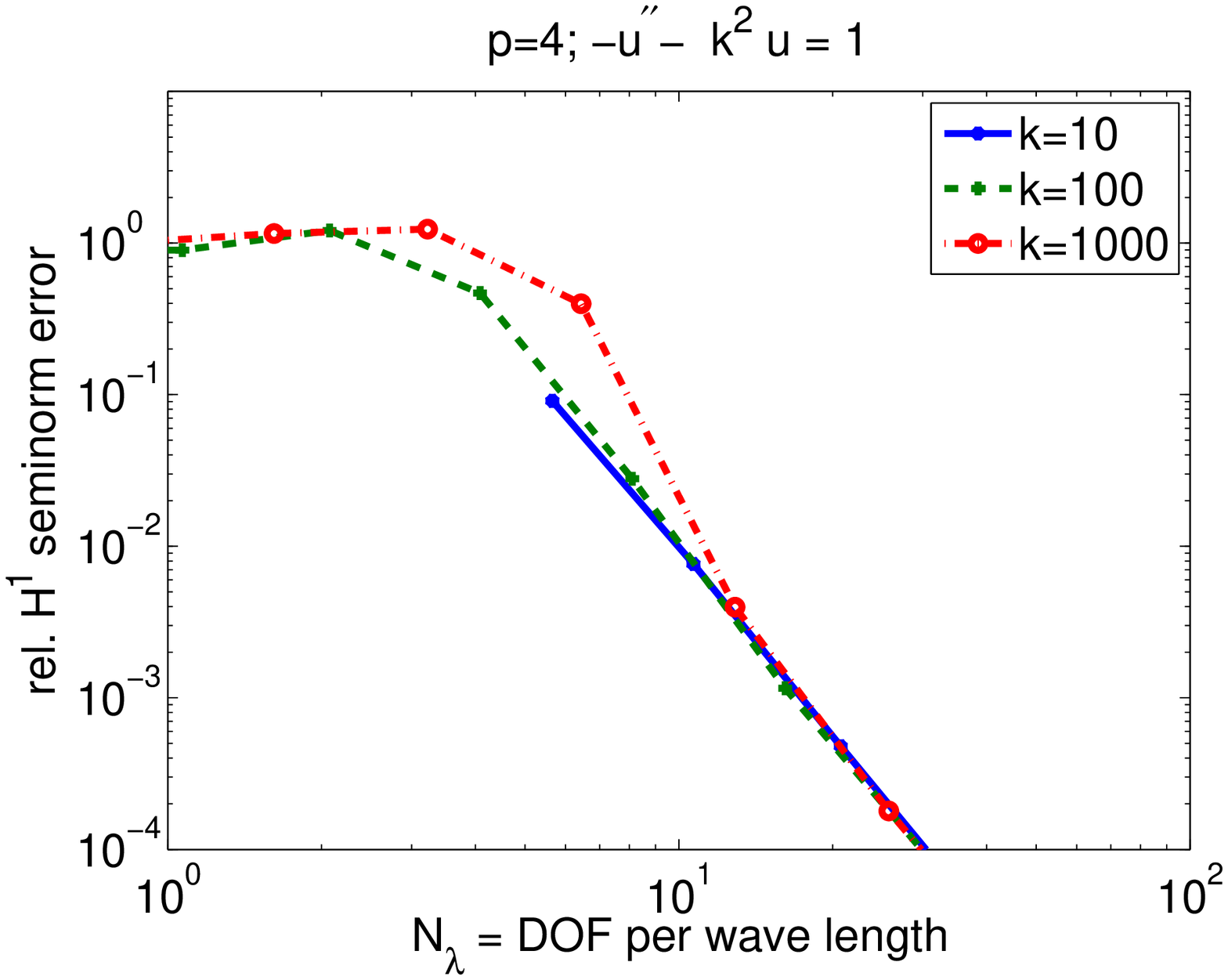}
\caption{\label{fig:1D-example} Performance of $h$-FEM for (\ref{eq:1D-model}).  
Top: $p=1$, $p=2$. Bottom: $p=3$, $p=4$ 
(see Example~\ref{ex:1d-pollution}).}
\end{figure}

In this article, we concentrate on numerical schemes for 
the Helmholtz equation at large wavenumbers $k$. Standard 
discretizations face several challenges, notably:
\begin{enumerate}[(I)]
\item 
\label{item:I}
For large wavenumber $k$, the solution $u$ is highly oscillatory. 
Its resolution, therefore, requires fine meshes, namely, 
at least $N = k^d$ degrees of freedom, where $d$ is the 
spatial dimension. 
\item 
\label{item:II}
The standard $H^1$-conforming variational formulation is indefinite, and 
stability on the discrete level is therefore an issue. A manifestation of 
this problem is the so-called ``pollution'', which expresses the observation
that much more stringent conditions on the discretization have to be met 
than the minimal $N = O(k^d)$ to achieve a given accuracy.
\end{enumerate} 
The second point, which will be the focus of the article, is best seen in the following, 
one-dimension example: 
\begin{Example}
\label{ex:1d-pollution}
{\rm 
For the boundary value problem
\begin{equation}
\label{eq:1D-model}
-u^{\prime\prime} - k^2 u = 1 \quad \mbox{ in $(0,1)$, } 
\quad u(0) = 0,\qquad u^\prime(1) - \bi k u(1) = 0,
\end{equation}
we consider the $h$-version finite element method (FEM) on uniform meshes with 
mesh size $h$ for different
approximation orders $p \in \{1,2,3,4\}$ and wavenumbers 
$k \in \{1,10,100\}$. Fig.~\ref{fig:1D-example} shows the 
relative error in the $H^1(\Omega)$-semi norm 
(i.e., $|u - u_N|_{H^1(\Omega)}/|u|_{H^1(\Omega)}$, where $u_N$
is the FEM approximation) versus the number of degrees of freedom
per wavelength $N_\lambda:= N/\lambda = 2 \pi N/k$ 
with $N$ being the dimension of the finite element space employed. 
We observe several effects in Fig.~\ref{fig:1D-example}: Firstly, 
since the solution $u$ of (\ref{eq:1D-model}) is smooth, higher order methods lead 
to higher accuracy for a given number of degrees of freedom per wavelength than 
lower order methods. Secondly, {\em asymptotically}, the FEM is quasioptimal with 
the finite element error $|u - u_N|_{H^1(\Omega)}$ satisfying 
\begin{equation}
\label{eq:1D-example:asymptotic-regime}
|u - u_N|_{H^1(\Omega)}\approx C_p N_\lambda^{-p} |u|_{H^1(\Omega)} 
\end{equation} 
for a constant $C_p$ independent of $k$. 
Thirdly, the performance of the FEM as measured 
in ``error vs. number of degrees of freedom per wavelength'' does depend 
on $k$: As $k$ increases, the preasymptotic range with reduced FEM performance becomes
larger. Fourthly, higher order methods are less sensitive to $k$
than lower order ones, i.e., for given $k$, high order methods enter 
the asymptotic regime in which 
(\ref{eq:1D-example:asymptotic-regime}) holds for smaller values of 
$N_\lambda$ than lower order methods. 
\eremk
}
\end{Example}
The behavior of the FEM in Example~\ref{ex:1d-pollution} has been analyzed 
in \cite{ihlenburg-babuska97,ihlenburg98}, where error bounds of the 
form (see \cite[Thm.~{4.27}]{ihlenburg98}) 
\begin{equation}
\label{eq:ihlenburg-estimate}
|u  -u_N|_{H^1(\Omega)} \leq C_p \left( 1 + k^{p+1} h^p\right) h^p |u|_{H^{p+1}(\Omega)}
\end{equation}
are established for a constant $C_p$ depending only on the approximation
order $p$. In this particular example, it is also easy to see that 
$|u|_{H^{p+1}(\Omega)}/|u|_{H^1(\Omega)} \sim k^p$, so that (\ref{eq:ihlenburg-estimate})
can be recast in the form 
\begin{equation}
\label{eq:ihlenburg-estimate-2}
|u  -u_N|_{H^1(\Omega)}
\leq C_p \left( 1 + k^{p+1} h^p\right) (kh)^p |u|_{H^1(\Omega)} 
\sim \left(1+ k N_\lambda^{-p}\right) N_\lambda^{-p} |u|_{H^1(\Omega)}. 
\end{equation}
This estimate goes a long way to explain the above observations. 
The presence of the factor $1 + k N_\lambda^{-p}$ explains 
the ``pollution effect'', i.e., the observation
that for fixed $N_\lambda$, the (relative) error of the FEM as compared
with the best approximation (which is essentially proportional to $N_\lambda^{-p}$
in this example) increases with $k$. The estimate (\ref{eq:ihlenburg-estimate-2})
also indicates that the asymptotic convergence behavior 
(\ref{eq:1D-example:asymptotic-regime}) is reached for $N_\lambda = O(k^{1/p})$. 
This confirms the observation made above that higher order methods are less prone
to pollution than lower order methods. 
Although Example~\ref{ex:1d-pollution} is restricted to 1D, 
similar observations have been made in the literature also for 
multi-d situations as early as \cite{bayliss85}. 
We emphasize that for uniform meshes (as in Example~\ref{ex:1d-pollution}) or, 
more generally, translation invariant meshes, complete and detailed dispersion analyses 
are available in an $h$-version setting, 
\cite{ihlenburg-babuska97,ihlenburg98,deraemaeker-babuska-bouillard99,ainsworth-monk-muniz06}, 
and in a $p$/$hp$-setting, \cite{ainsworth-monk-muniz06,ainsworth04,ainsworth-wajid09}, 
that give strong mathematical evidence for the superior ability of 
high order methods to cope with the pollution effect. 

The present paper, which 
discusses and generalizes the work \cite{melenk-sauter10,melenk-sauter11}, proves that also
on unstructured meshes, high order methods are less prone to pollution. 
More precisely, for a large class of Helmholtz problems, stability and quasi-optimality is given 
under the scale resolution condition 
\begin{equation}
\label{eq:scale-resolution}
\frac{kh}{p} \leq c_1 \qquad \mbox{ and } \qquad 
p \ge c_2 \log k,  
\end{equation}
where $c_1$ is sufficiently small and $c_2$ sufficiently large. For piecewise smooth
geometries (e.g., polygons), additionally appropriate mesh refinement near the singularities 
is required. 

We close our discussion of Example~\ref{ex:1d-pollution} by remarking
that its analysis and, in fact, the analysis of significant parts 
of this article rests on $H^1$-like norms. Largely,
this choice is motivated by the numerical scheme, namely, an 
$H^1$-conforming FEM. 
\subsection{Non-standard FEM} 
The limitations of the classical FEM mentioned above in (\ref{item:I}) and (\ref{item:II}) 
have sparked a significant amount of research in the past decades to overcome or at 
least mitigate them. This research focuses on two techniques that are often considered
in tandem: firstly, the underlying approximation by classical piecewise
polynomials is replaced with special, problem-adapted functions such as systems of 
plane waves; secondly, the numerical scheme is based on a variational
formulation different from the classical $H^1$-conforming Galerkin approach. 
Before discussing these ideas in more detail, we point the 
reader to the interesting work \cite{babuska-sauter00}, which shows for a model
situation on regular, infinite grids in 2D that no 
9-point stencil (i.e., a numerical method based on connecting the value at a node with 
those of the 8 nearest neighbors)
generates a completely
pollution-free method; the 1D situation is special and discussed 
briefly in 
\iftechreport 
Section~\ref{sec:1D}. 
\fi

Work that is based on a new or modified variational formulation but rests 
on the approximation properties of piecewise polynomials includes 
the Galerkin Least Squares Method \cite{harari-hughes92,harari06}, the methods of 
\cite{babuska-paik-sauter94}, and Discontinuous
Galerkin Methods (\cite{feng-wu09,feng-wu11,feng-xing10} and references there). 
Several methods have been proposed that are based on the approximation properties
of special, problem-adapted systems of functions such as systems of plane waves.
In an $H^1$-conforming Galerkin setting, this idea has been pursued in the 
Partition of Unity Method/Generalized FEM by several authors, e.g., \cite%
{melenk95,babuska-melenk96b, laghrouche-bettes00,
laghrouche-bettess-astley02, perrey-debain-laghrouche-bettes-trevelyan04,
astley-gamallo05,ortiz04, strouboulis-babuska-hidajat06, huttunen-gamallo-astley09}.  
A variety of other methods that are based on problem-adapted ansatz functions 
leave the $H^1$-conforming Galerkin setting and enforce the jump across element
boundaries in a weak sense. This can be done by least squares techniques (%
\cite{stojek98, monk-wang99, li07, betcke-barnett10, desmet07} and
references there), by Lagrange multiplier techniques as in the
Discontinuous Enrichment Method \cite{farhat-harari-hetmaniuk03,
farhat-tezaur-weidemann04, farhat-tezaur06} or by Discontinuous Galerkin
(DG) type methods, 
\cite{cessenat-despres98,cessenat-despres03,
huttunen-monk07,luostari-huttunen-monk09,
buffa-monk07,gittelson-hiptmair-perugia09,
hiptmair-moiola-perugia09a, hiptmair-moiola-perugia11, moiola09,
monk-schoeberl-sinwel10}; in these last references, we have included the work on
the Ultra Weak Variational Formulation (UWVF) since it can be understood
as a special DG method as discussed in  \cite{buffa-monk07,gittelson-hiptmair-perugia09}. 
%
\subsection{Scope of the article}
The present article focuses on the stability properties of 
numerical methods for Helmholtz problems and exemplarily 
discusses three different approaches in more detail for their
differences in techniques. The first
approach, studied in Section~\ref{sec:discrete-stability-H1-conforming}, 
is that of the classical $H^1$-based Galerkin method for 
Helmholtz problems. The setting is that of a G{\aa}rding inequality so that 
stability of a numerical method can be inferred from the stability of the 
continuous problem by perturbation arguments. This motivates us to study
for problem (\ref{eq:problem}), which will serve as our 
model Helmholtz problem in this article, the stability properties of the 
continuous problem in Section~\ref{sec:stability-continuous-problem}. 
In order to
make the perturbation argument explicit in the wavenumber $k$, a detailed, 
$k$-explicit regularity analysis for Helmholtz problems is necessary. This
is worked out in Section~\ref{sec:regularity} for our model problem (\ref{eq:problem})
posed on 
polygonal domains. These results generalize a similar regularity theory 
for convex polygons or domains with analytic boundary of 
\cite{melenk-sauter10,melenk-sauter11}.  Structurally similar results have been 
obtained in connection with boundary integral formulations 
in \cite{melenk10a,melenk-loehndorf10}.  

We discuss in Sections~\ref{sec:least-squares} and \ref{sec:DG} somewhat briefly
a second and a third approach to stability of numerical schemes. In contrast
to the setting discussed above, where stability is only ensured asymptotically
for sufficiently fine discretizations, these methods 
are stable by construction and can even feature quasioptimality in appropriate
residual norms. We point out, however, that relating this residual norm to a 
more standard norm such as the $L^2$-norm for the error is a non-trivial task. 
Our presentation for these methods will follow the works 
\cite{monk-wang99,buffa-monk07,gittelson-hiptmair-perugia09,hiptmair-moiola-perugia11}.

\iftechreport 
Many aspects of discretizations for Helmholtz problems are not addressed
in this article. For recent developments in boundary element techniques 
for this 
problem class, we refer to the survey article \cite{chandler-wilde-graham09}. 
The model problem (\ref{eq:problem}) discussed here involves the rather 
simple boundary condition (\ref{eq:problem-b}), which can be understood
as an approximation to a Dirichlet-to-Neumann operator that provides a 
coupling to a homogeneous Helmholtz equation in an exterior domain together
with appropriate radiation conditions at infinity. A variety of techniques 
for such problems are discussed in \cite{givoli92}. Further methods 
include FEM-BEM coupling, the PML due to B\'erenger 
(see \cite{collino-monk98,bramble-pasciak08} and references therein), 
infinite elements \cite{demkowicz-gerdes98},
and methods based on the pole condition, \cite{hohage-nannen09}. 
Another topic not addressed here is the solution of the arising linear system; 
we refer the reader to \cite{erlangga08,engquist-ying10} for a discussion 
of the state of the art. 
Further works with survey character includes 
\cite{engquist-runborg03,ihlenburg98,ihlenburg09,thompson05}.
\fi

\subsection{Some notation}
\iftechreport 
We employ standard notation for Sobolev spaces, 
\cite{adams75a,necas67,brenner-scott94,schwab98}.  
\fi 
For domains $\omega$ and $k \ne 0$ we denote  
\begin{equation}
\label{eq:norm}
\|u\|_{1,k,\omega}^2:=k^2 \|u\|^2_{L^2(\omega)} + \|\nabla u\|^2_{L^2(\omega)}. 
\end{equation}
This norm is equivalent to the standard $H^1$-norm. The presence of the weight $k$ 
in the $L^2$-part leads to a 
balance between the $H^1$-seminorm and the $L^2$-norm for functions with the expected 
oscillatory behavior such as plane waves $e^{\bi k {\bf d} \cdot x}$ 
(with ${\bf d}$ being a unit vector). 
Additionally, the bilinear form $B$ considered below is bounded uniformly 
in $k$ with respect to this ($k$-dependent) norm. 

Throughout this work, a standing assumption will be 
\begin{equation}
\label{eq:k0}
|k| \ge k_0 > 0; 
\end{equation}
our frequently used phrase ``independent of $k$'' will still implicitly assume (\ref{eq:k0}). 
We denote by $C$ a generic constant. If not stated otherwise, $C$ will be independent
of the wavenumber $k$ but may depend on $k_0$. 
For smooth functions 
$u$ defined on a $d$-dimensional manifold, we employ the notation
$
\displaystyle 
|\nabla^n u(x)|^2:= \!\!\!\!\!\sum_{\alpha \in \BbbN_0^d\colon |\alpha| = n} \frac{|\alpha|!}{\alpha!} |D^\alpha u(x)|^2. 
$

\subsection{A model problem}
\label{sec:model-problem}
In order to fix ideas, we will use the following, 
specific model problem: For a bounded Lipschitz domain 
$\Omega \subset \BbbR^d$, $d \in \{2,3\}$, we study 
for $k \in \BbbR$, $|k| \ge k_0$, 
\iftechreport 
the boundary value problem 
\fi
\begin{subequations}
\label{eq:problem}
\begin{eqnarray}
\label{eq:problem-a}
-\Delta u - k^2 u &=& f \mbox{ in $\Omega$}, \\
\label{eq:problem-b}
\partial_n u + \bi k u &=& g \mbox{ on $\partial\Omega$}.
\end{eqnarray}
\end{subequations}
Henceforth, to simplify the
notation, we assume $k \ge k_0 >0$ but point out that the choice 
of the sign of $k$ is not essential. The weak formulation for 
(\ref{eq:problem}) is:
\begin{equation}
\label{eq:problem-weak}
\mbox{Find $u \in H^1(\Omega)$ s.t.} 
\qquad 
B(u,v) = l(v) \qquad \forall v \in H^1(\Omega),
\end{equation}
where,
for $f \in L^2(\Omega)$ and $g \in L^2(\partial\Omega)$, 
$B$ and $l$ are given by 
\begin{equation}
\label{eq:B}
B(u,v):=\int_\Omega (\nabla u \cdot \nabla \overline{v} - k^2 u \overline{v}) + 
\bi k \int_{\partial\Omega} u \overline{v}, 
\quad  l(v) := (f,v)_{L^2(\Omega)} + (g,v)_{L^2(\partial\Omega)}. 
\end{equation}
As usual, 
if $f \in (H^1(\Omega))^\prime$ and $g \in H^{-1/2}(\partial\Omega)$, then the 
$L^2$-inner products $(\cdot,\cdot)_{L^2(\Omega)}$ and $(\cdot,\cdot)_{L^2(\partial\Omega)}$
are understood as duality pairings. The multiplicative trace inequality proves 
continuity of $B$; in fact, there exists $C_B  > 0$ independent of $k$ such that 
(see, e.g., \cite[Cor.~{3.2}]{melenk-sauter10} for details)
\begin{equation}
\label{eq:CB}
|B(u,v)| \leq C_B \|u\|_{1,k,\Omega} \|v\|_{1,k,\Omega} 
\qquad \forall u,v \in H^1(\Omega).
\end{equation}
\section{Stability of the continuous problem}
\label{sec:stability-continuous-problem}
Helmholtz problems can often be cast in the form 
``coercive + compact perturbation'' where the compact perturbation
is $k$-dependent. In other words, a G{\aa}rding inequality is satisfied.
For example, the sesquilinear form $B$ of 
(\ref{eq:B}) is of this form since 
\begin{equation}
\label{eq:garding-inequality}
\operatorname*{Re} B(u,u) + 2 k^2 (u,u)_{L^2(\Omega)} = \|u\|^2_{1,k,\Omega} 
\end{equation}
and the embedding $H^1(\Omega)\subset L^2(\Omega)$ is compact by Rellich's theorem. 
Classical Fredholm theory (the ``Fredholm alternative'') then yields 
unique solvability of (\ref{eq:problem-weak}) for all 
$f \in (H^1(\Omega))^\prime$ and $g \in H^{-1/2}(\partial\Omega)$, if
one can show uniqueness. Uniqueness in turn is often obtained by 
exploiting analyticity of the solutions of homogeneous Helmholtz equation, 
or, more generally, the unique continuation principle for elliptic problems, 
(see, e.g., \cite[Chap.~{4.3}]{leis86}): 
\begin{Example}[Uniqueness for (\ref{eq:problem})]
\label{example:uniqueness}
{\rm 
Let $f = 0$ and $g = 0$ in (\ref{eq:problem}). 
Then, any solution $u\in H^1(\Omega)$ 
of (\ref{eq:problem}) satisfies $u|_{\partial\Omega} = 0$ since 
$0 = \operatorname*{Im} B(u,u) = k \|u\|^2_{L^2(\partial\Omega)}$
(see Lemma~\ref{lemma:a-priori}).
Hence, the trivial extension $\widetilde u$ 
to $\BbbR^2$ satisfies $\widetilde u \in H^1(\BbbR^2)$. The observations
$B(u,v) = 0$ for all $v \in H^1(\Omega)$ and $u|_{\partial\Omega} = 0$ 
show
$$
\int_{\BbbR^2} \nabla \widetilde u \cdot \nabla \overline v-  k^2 \widetilde u \overline{v}
 = 0 \qquad \forall v \in C^\infty_0(\BbbR^2). 
$$
Hence, $\widetilde u$ is a solution of the homogeneous Helmholtz equation and 
$\widetilde u$ vanishes on $\BbbR^2 \setminus \overline{\Omega}$. Analyticity
of $\widetilde u$ (or, more generally, the unique continuation principle 
presented in \cite[Chap.~{4.3}]{leis86})
then implies that $\widetilde u \equiv 0$, which in turn yields $u \equiv 0$. 
\eremk
}
\end{Example}
The arguments based on the Fredholm alternative 
do not give any indication of how the solution operator depends 
on the wavenumber $k$. Yet, it is clearly of interest to know how $k$ 
enters bounds for the 
solution operator. It turns out that both the geometry and the type of 
boundary conditions strongly affect these bounds. For example, 
for an exterior Dirichlet problem, 
\cite{betcke-chandler-wilde-graham-langdon-lindner10}
exhibits a geometry and a sequence of wavenumber $(k_n)_{n \in \BbbN}$ tending
to infinity such that the norm of the solution operator for these wavenumbers
is bounded from below by an exponentially growing term $C e^{b k_n}$ 
for some $C$, $b > 0$. These geometries 
feature so-called ``trapping'' or near-trapping and are not convex. 
For convex or at least star-shaped geometries, the $k$-dependence is much 
better behaved. An important ingredient of the analysis on such geometries are 
special test functions in the variational formulation. For example, 
assuming in the the model problem (\ref{eq:problem-weak}) that $\Omega$ is 
star-shaped with respect to the origin (and has a smooth boundary), one may take 
as the test function $v(x) = x \cdot \nabla u(x)$, where $u$ is the exact solution.  
An integration by parts (more generally, the so-called ``Rellich identities''
\cite[p.~{261}]{necas67} or an identity due to Poho\v{z}aev, \cite{pohozaev65}) 
then leads to the following estimate for the model
problem (\ref{eq:problem-weak}): 
\begin{equation}
\label{eq:a-priori-star-shaped}
\|u\|_{1,k,\Omega} \leq C \left[ \|f\|_{L^2(\Omega)} + \|g\|_{L^2(\partial\Omega)}\right]; 
\end{equation}
this was shown in \cite[Prop.~{8.1.4}]{melenk95} (for $d=2$) and subsequently 
by \cite{cummings-feng06} for $d = 3$. Uniform in $k$ bounds were established 
in \cite{hetmaniuk07} for star-shaped domains and certain boundary conditions of 
mixed type by related techniques. The same test function was also crucial 
for a boundary integral setting in \cite{monk-chandler-wilde08}. A refined version
of this test function that goes back to Morawetz and Ludwig, 
\cite{morawetz-ludwig68} was used recently 
in a boundary integral equations context
(still for star-shaped domains), 
\cite{spence-chandler-wilde-graham-smyshlyaev10}. 

While (\ref{eq:a-priori-star-shaped}) does not make minimal assumptions on the 
regularity of $f$ and $g$, the estimate (\ref{eq:a-priori-star-shaped}) can be used
to show that (for star-shaped domains) the sesquilinear form $B$ of (\ref{eq:problem-weak}) 
satisfies an inf-sup condition with inf-sup constant $\gamma = O(k^{-1})$---this can be 
shown using the arguments presented in the proof Theorem~\ref{thm:inf-sup}.

\bigskip
An important ingredient of the regularity and stability theory will be the concept of 
{\em polynomial well-posedness} by which we mean polynomial-in-$k$-bounds for the 
norm of the solution operator. 
The model problem (\ref{eq:problem}) on star-shaped domains with 
the {\sl a priori} bound (\ref{eq:a-priori-star-shaped}) is an example.  
The following Section~\ref{sec:polynomial-well-posedness} shows polynomial 
well-posedness for the model problem (\ref{eq:problem}) on general Lipschitz domains 
(Thm.~\ref{thm:a-priori}). 
It is thus not the geometry but the type of boundary conditions in our 
model problem (\ref{eq:problem}), namely, Robin boundary conditions that makes 
it polynomially well-posed. In contrast, the Dirichlet boundary conditions in conjunction
with the lack of star-shapedness in the examples given in 
\cite{betcke-chandler-wilde-graham-langdon-lindner10} 
make these problem not polynomially well-posed.  
\subsection{Polynomial well-posedness for the model problem (\ref{eq:problem})}
\label{sec:polynomial-well-posedness}
\begin{Lemma}
\label{lemma:a-priori}
Let $\Omega \subset \BbbR^d$ be a bounded Lipschitz domain. Let $u \in H^1(\Omega)$ be a 
weak solution of (\ref{eq:problem}) with $f = 0$ and $g \in L^2(\partial\Omega)$. Then 
$
\|u\|_{L^2(\partial\Omega)} \leq  k^{-1} \|g\|_{L^2(\partial\Omega)}. 
$
\end{Lemma}
\begin{proof}
Selecting $v = u$ in the weak formulation (\ref{eq:problem-weak}) and considering
the imaginary part yields 
$$
k \|u\|^2_{L^2(\partial\Omega)} = 
\operatorname*{Im} \int_{\partial\Omega} g \overline{u} \leq \|g\|_{L^2(\partial\Omega)} \|u\|_{L^2(\partial\Omega)}. 
$$
This concludes the argument. 
\qed
\end{proof}

Next we use results on layer potentials for the Helmholtz equation from \cite{melenk10a} to prove the following lemma:
\begin{Lemma}
\label{lemma:use-of-potentials}
Let $\Omega\subset\BbbR^d$ be a bounded Lipschitz domain, 
$u \in H^1(\Omega)$ solve (\ref{eq:problem}) with $f = 0$. 
Assume $u|_{\partial\Omega} \in L^2(\partial\Omega)$ and 
$\partial_n u \in L^2(\partial\Omega)$. Then there exists $C > 0$ independent 
of $k$ and $u$ such that
\begin{eqnarray*}
\|u\|_{L^2(\Omega)} &\leq& C 
k\left( \|u\|_{L^2(\partial\Omega)}+\|\partial_n u\|_{H^{-1}(\partial\Omega)}\right),\\
\|u\|_{1,k,\Omega}  & \leq& C \left[ 
k^2 \|u\|_{L^2(\partial\Omega)} + k^2 \|\partial_n u\|_{H^{-1}(\partial\Omega)} + 
k^{-2} \|\partial_n u\|_{L^2(\partial\Omega)}\right].
\end{eqnarray*}
\end{Lemma}
\begin{proof}
With the single layer and double layer potentials $\widetilde V_k$ and $\widetilde K_k$ we have 
the representation formula 
$
u = \widetilde V_k \partial_n u - \widetilde K_k u 
$. From \cite[Lemmata~{2.1}, {2.2}, Theorems~{4.1}, {4.2}]{melenk10a} we obtain 
\begin{eqnarray*}
\|\widetilde V_k \partial_n u\|_{L^2(\Omega)} \leq C k \|\partial_n u\|_{H^{-1}(\partial\Omega)}, 
\qquad 
\|\widetilde K_k u\|_{L^2(\Omega)} \leq C k \|u\|_{L^{2}(\partial\Omega)}. 
\end{eqnarray*}
\iftechreport
Thus, 
$$
\|u\|_{L^2(\Omega)} \leq C  k \left( \|u\|_{L^2(\partial\Omega)} 
+ \|\partial_n u\|_{H^{-1}(\partial\Omega)}\right). 
$$
\fi
Next, using $v = u$ in the weak formulation (\ref{eq:problem-weak})
yields 
$$
\|\nabla u\|^2_{L^2(\Omega)} \leq C \left[ k^2 \|u\|^2_{L^2(\Omega)} + 
\|\partial_n u\|_{L^2(\partial\Omega)} \|u\|_{L^2(\partial\Omega)}
\right]
$$
and therefore 
\begin{eqnarray*}
\|\nabla u\|^2_{L^2(\Omega)} + k^2 \|u\|^2_{L^2(\Omega)} &\leq& C \left[ 
k^4 \|u\|^2_{L^2(\partial\Omega)} + k^4 \|\partial_n u\|^2_{H^{-1}(\partial\Omega)} + 
k^{-4} \|\partial_n u\|^2_{L^2(\partial\Omega)}\right],
\end{eqnarray*}
which concludes the proof. 
\qed
\end{proof}
\begin{Theorem}
\label{thm:a-priori}
Let $\Omega\subset\BbbR^d$, $d \in \{2,3\}$ be a bounded Lipschitz domain. Then there exists $C > 0$ 
(independent of $k$) such that for $f \in L^2(\Omega)$ and $g \in L^2(\partial\Omega)$
the solution $u\in H^1(\Omega)$ of (\ref{eq:problem}) satisfies 
$$
\|u\|_{1,k,\Omega}
\leq C \left[ 
k^{2} \|g\|_{L^2(\partial\Omega)} + k^{5/2} \|f\|_{L^2(\Omega)} \right].
$$
\end{Theorem}
\begin{proof}
We first transform the problem to one with homogeneous right-hand side $f$
in the standard way. A particular solution of the equation (\ref{eq:problem-a}) is 
given by the Newton potential 
$
u_0:= G_k \star f
$; here, $G_k$ is a Green's function for the Helmholtz equation and 
we tacitly extend $f$ by zero outside $\Omega$.  
Then $u_0 \in H^2_{loc}(\BbbR^d)$ and by the analysis of the Newton potential
given in \cite[Lemma~{3.5}]{melenk-sauter10} we have 
\begin{equation}
\label{eq:thm:a-priori-10}
k^{-1}\|u_0\|_{H^2(\Omega)} + \|u_0\|_{H^1(\Omega)}  + k\|u_0\|_{L^2(\Omega)} 
\leq C \|f\|_{L^2(\Omega)}. 
\end{equation}
The difference $\widetilde u:= u - u_0$ then satisfies 
\iftechreport 
\begin{subequations}
\label{eq:thm:a-priori-100}
\begin{eqnarray}
-\Delta \widetilde u  - k^2 \widetilde u &=& 0 \quad \mbox{ in $\Omega$}, \\
\partial_n \widetilde u  + \bi k \widetilde u &=& g - 
\left( \partial_n u_0 + \bi k u_0
\right) =:\widetilde g. 
\end{eqnarray}
\end{subequations}
\fi 
We have with the multiplicative trace inequality  
\begin{eqnarray}
\nonumber 
\|\widetilde g\|_{L^2(\partial\Omega)} &\leq& C \left[ \|g\|_{L^2(\partial\Omega)} + 
\|u_0\|_{H^2(\Omega)}^{1/2} \|u_0\|_{H^1(\Omega)}^{1/2} + 
k \|u_0\|_{H^1(\Omega)}^{1/2} \|u_0\|_{L^2(\Omega)}^{1/2}  
\right]
\\
&\leq & C \left[ \|g\|_{L^2(\partial\Omega)} + k^{1/2} \|f\|_{L^2(\Omega)}\right]. 
\label{eq:thm:a-priori-200}
\end{eqnarray}
To get bounds on $\widetilde u$, we employ Lemma~\ref{lemma:a-priori} and 
(\ref{eq:thm:a-priori-200})
to conclude 
\begin{eqnarray}
\|\widetilde u\|_{L^2(\partial\Omega)} &\leq& C  k^{-1} \|\widetilde g\|_{L^2(\partial\Omega)} 
\leq C \left[ k^{-1} \|g\|_{L^2(\partial\Omega)} + k^{-1/2} \|f\|_{L^2(\Omega)}\right], \\
\|\partial_n \widetilde u\|_{L^2(\partial\Omega)} &\leq& C \left[ 
\|\widetilde g\|_{L^2(\partial\Omega)} + k \|\widetilde u\|_{L^2(\partial\Omega)} \right]
\leq C \left[ \|g\|_{L^2(\partial\Omega)} + k^{1/2}\|f\|_{L^2(\Omega)}\right].  
\qquad 
\end{eqnarray}
Lemma~\ref{lemma:use-of-potentials} and 
the generous estimate $\|\partial_n \widetilde u\|_{H^{-1}(\partial\Omega)} 
\leq C \|\partial_n \widetilde u\|_{L^2(\partial\Omega)}$ produce
\begin{eqnarray}
\label{eq:thm:a-priori-2000}
\|\widetilde u\|_{H^1(\Omega)} + k \|\widetilde u\|_{L^2(\Omega)} 
& \leq & C \left[ k^{2} \|g\|_{L^2(\partial\Omega)} + k^{5/2} \|f\|_{L^2(\Omega)} \right].
\end{eqnarray}
Combining (\ref{eq:thm:a-priori-10}), (\ref{eq:thm:a-priori-2000}) finishes the argument.
\qed
\end{proof}
The {\sl a priori} estimate of 
Theorem~\ref{thm:a-priori} does not make minimal assumptions on the regularity of 
$f$ and $g$. However, it can be used to obtain estimates on the inf-sup 
and hence {\sl a priori} bounds for $f \in (H^1(\Omega))^\prime$ and 
$g \in H^{-1/2}(\partial\Omega)$ as we now show: 
\begin{Theorem}
\label{thm:inf-sup}
Let $\Omega\subset\BbbR^d$, $d \in \{2,3\}$ be a bounded Lipschitz domain. 
Then there exists $C > 0$ (independent of $k$) such that the sesquilinear 
form $B$ of (\ref{eq:B}) satisfies 
\begin{equation}
\label{eq:thm:inf-sup-1}
\inf_{0 \ne u \in H^1(\Omega)} \sup_{0 \ne v \in H^1(\Omega)} 
\frac{\operatorname*{Re} B(u,v)}{\|u\|_{1,k,\Omega} \|v\|_{1,k,\Omega}}
\ge C k^{-7/2}. 
\end{equation}
Furthermore, for every $f \in (H^1(\Omega))^\prime$ and $g \in H^{-1/2}(\partial\Omega)$
the problem 
(\ref{eq:problem-weak}) is uniquely solvable, and its solution $u\in H^1(\Omega)$ 
satisfies the {\sl a priori} bound 
\begin{equation}
\label{eq:thm:inf-sup-2}
\|u\|_{1,k,\Omega} \leq C k^{7/2}\left[ \|f\|_{(H^1(\Omega))^\prime} + \|g\|_{H^{-1/2}(\partial\Omega)}\right].
\end{equation}
If $\Omega$ is convex or if $\Omega$ is star-shaped and has a smooth boundary, then 
the following, sharper estimate holds:
\begin{equation}
\label{eq:thm:inf-sup-100}
\inf_{0 \ne u \in H^1(\Omega)} \sup_{0 \ne v \in H^1(\Omega)} 
\frac{\operatorname*{Re} B(u,v)}{\|u\|_{1,k,\Omega} \|v\|_{1,k,\Omega}}
\ge C k^{-1}. 
\end{equation}
\end{Theorem}
\begin{proof}
The proof relies on standard arguments for sesquilinear forms satisfying a G{\aa}rding inequality. 
For simplicity of notation, we write $\|\cdot\|_{1,k}$ for $\|\cdot\|_{1,k,\Omega}$. 

Given $u \in H^1(\Omega)$ we define 
$z \in H^1(\Omega)$ as the solution of 
$$
2 k^2 (\cdot,u)_{L^2(\Omega)} = B(\cdot,z). 
$$
Theorem~\ref{thm:a-priori} implies 
$
\|z\|_{1,k} \leq C k^{9/2} \|u\|_{L^2(\Omega)}, 
$
and $v = u + z$ satisfies 
$$
\operatorname*{Re} B(u,v) = \operatorname*{Re} B(u,u) + 
\operatorname*{Re} B(u,z) = \|u\|^2_{1,k} - 2 k^2 \|u\|^2_{L^2(\Omega)} + 
\operatorname*{Re} B(u,z) = \|u\|^2_{1,k}. 
$$
Thus,
\begin{eqnarray*}
\operatorname*{Re} B(u,v) &=& \|u\|^2_{1,k}, \\
\|v\|_{1,k} &=& \|u + z\|_{1,k} \leq \|u\|_{1,k} + \|z\|_{1,k} 
\leq \|u\|_{1,k} + C k^{9/2} \|u\|_{L^2(\Omega)} \leq C k^{7/2} \|u\|_{1,k}. 
\end{eqnarray*}
Therefore, 
$$
\operatorname*{Re} B(u,v) = \|u\|^2_{1,k} 
\ge \|u\|_{1,k} C k^{-7/2} \|v\|_{1,k}, 
$$
which concludes the proof of (\ref{eq:thm:inf-sup-1}). Example~\ref{example:uniqueness}
provides unique solvability for (\ref{eq:problem}) so that 
(\ref{eq:thm:inf-sup-1}) gives the {\sl a priori} estimate (\ref{eq:thm:inf-sup-2}). 
Finally, (\ref{eq:thm:inf-sup-100}) is shown by the same arguments
using (\ref{eq:a-priori-star-shaped}).
\qed
\end{proof}
\section{$k$-explicit regularity theory} 
\label{sec:regularity}
\subsection{Regularity by decomposition} 
Since the Sobolev regularity of elliptic problems is determined by the 
leading order terms of the differential equation and the boundary conditions, 
the Sobolev regularity properties of our model problem (\ref{eq:problem}) 
are the same as those for the Laplacian. However, regularity results that 
are explicit in the wavenumber $k$ are clearly of interest; for example, 
we will use them in Section~\ref{sec:stability-of-hpFEM} below to 
quantify how fine the discretization has to be (relative to $k$) so 
that the FEM is stable and quasi-optimal.

The $k$-explicit regularity theory developed in \cite{melenk-sauter10,melenk-sauter11} 
(and, similarly, for integral equations in \cite{melenk-loehndorf10,melenk10a})
takes the form of an additive splitting of the solution into a part with 
finite regularity but $k$-independent bounds and a part that is analytic and for which $k$-explicit
bounds for all derivatives are available. Below, we will present a similar regularity theory
for the model problem (\ref{eq:problem}) for polygonal $\Omega\subset\BbbR^2$, thereby extending
the results of \cite{melenk-sauter11}, which restricted its analysis of polygons to 
the convex case. In order to motivate the ensuing developments, we quote 
from \cite{melenk-sauter10} a result that shows in a simple setting the key features of our 
$k$-explicit ``regularity by decomposition'': %
\begin{Lemma}[\protect{\cite[Lemma~{3.5}]{melenk-sauter10}}]
\label{lemma:decomposition-lemma}
Let $B_R(0) \subset \BbbR^d$, $d \in \{1,2,3\}$ be the ball of radius $R$ centered at the origin. 
Then, there exist $C$, $\gamma > 0$ such that for all $k$ (with $k \ge k_0$) the following 
is true: For all $f \in L^2(\BbbR^d)$ with $\operatorname*{supp} f \subset B_R(0)$ 
the solution $u$ of 
$$
-\Delta u - k^2 u = f \quad \mbox{ in $\BbbR^d$ }, 
$$
subject to the 
Sommerfeld radiation condition 
$$
\lim_{|x| \to \infty} |x|^{\frac{d-1}{2}} \Big( \frac{\partial u}{\partial |x|} - iku \Big) = 0 
\quad \text{for} \quad  |x| \rightarrow \infty,
$$
has the following regularity properties: 
\begin{enumerate}[(i)]
\item 
\label{item:decomposition-i}
$u|_{B_{2R}(0)} \in H^2(B_{2R}(0))$ and $\|u\|_{H^2(B_{2R}(0))} \leq C k \|f\|_{L^2(B_R(0))}$. 
\item 
\label{item:decomposition-ii}
$u|_{B_{2R}(0)}$ can be decomposed as $u = u_{H^2} + u_{\mathcal A}$ for a  $u_{H^2} \in H^2(B_{2R})$
and an analytic $u_{\mathcal A}$ together with the bounds 
\begin{eqnarray*}
k\|u_{H^2}\|_{1,k,B_{2R}(0)}+\|u_{H^2}\|_{H^2(B_{2R}(0))} &\leq& C \|f\|_{L^2(B_R(0))},\\
\|\nabla^n u_{\mathcal A}\|_{L^2(B_{2R}(0))} &\leq& C \gamma^n \max\{n,k\}^{n-1} \|f\|_{L^2(B_R(0))} 
\quad \forall n \in \BbbN_0. 
\end{eqnarray*}
\end{enumerate}
\end{Lemma}
A few comments concerning Lemma~\ref{lemma:decomposition-lemma} are in order. 
For general $f \in L^2(B_R(0))$, one cannot expect 
better regularity than $H^2$-regularity for the solution $u$ and, indeed, both regularity results 
(\ref{item:decomposition-i}) and (\ref{item:decomposition-ii}) assert this. The estimate
(\ref{item:decomposition-i}) is sharp in its dependence on $k$ as the following simple example shows: 
For the fundamental solution $G_k$ (with singularity at the origin) and a cut-off function 
$\chi\in C^\infty_0(\BbbR^d)$ with $\operatorname*{supp} \chi \subset B_{2R}(0)$ 
and $\chi \equiv 1$  on $B_R(0)$, 
the functions $u:= (1-\chi) G_k$ and $f:= - \Delta u - k^2 u$ satisfy 
$\|u\|_{H^2(B_{2R}(0))} = O(k^2)$ 
and $\|f\|_{L^2(B_R(0))} = O(k)$. Compared to (\ref{item:decomposition-i}), the 
regularity assertion (\ref{item:decomposition-ii}) is finer in that its $H^2$-part $u_{H^2}$ 
has a better $k$-dependence. The $k$-dependence of the analytic part $u_{\mathcal A}$ is not improved
(indeed, $\|u_{\mathcal A}\|_{H^2(B_{2R}(0))} \leq C k \|f\|_{L^2(B_R(0))}$), but the analyticity 
of $u_{\mathcal A}$ is a feature that higher order methods can exploit.

The decomposition in (\ref{item:decomposition-ii}) of Lemma~\ref{lemma:decomposition-lemma}
is obtained by a decomposition of the datum $f$ using low pass and high pass filters, 
i.e., $f = L_{\eta k} f + H_{\eta k} f$, where the low pass filter $L_{\eta k}$ cuts off 
frequencies beyond $\eta k$ (here, $\eta > 1$) and $H_{\eta k}$ eliminates the frequencies
small than $\eta k$. Similar frequency filters will be important tools in our analysis
below as well (see Sec.~\ref{sec:filters}). The regularity properties stated 
in (\ref{item:decomposition-ii}) then follow from this decomposition and the 
explicit solution formula $u = G_k \star f$ 
(see \cite[Lemma~{3.5}]{melenk-sauter10} for details). 

Lemma~\ref{lemma:decomposition-lemma} serves as a prototype for 
``regularity theory by decomposition''. Similar decompositions have been developed recently
for several Helmholtz problems in \cite{melenk-sauter11} and 
\cite{melenk10a,melenk-loehndorf10}. Although they vary in their details, these 
decomposition are structurally similar in that they have the form of an additive splitting 
into a part with finite regularity with $k$-independent bounds and an analytic part 
with $k$-dependent bounds. The basic ingredients of these decomposition results are 
(a) (piecewise) analyticity of the geometry (or, more generally, the data of the problem) and 
(b) {\sl a priori} bounds for solution operator. The latter appear only in the estimate 
for the analytic part of the decomposition, and the most interesting case is that of polynomially
well-posed problems.  We illustrate the construction of the decomposition 
for the model problem (\ref{eq:problem}) in polygonal domains $\Omega\subset\BbbR^2$. This result 
is an extension to general polygons of the results \cite{melenk-sauter11}, which restricted its 
attention to the case of convex polygons.  We emphasize that the choice of the 
boundary conditions (\ref{eq:problem-b}) is not essential for the form of the decomposition and 
other boundary conditions could be treated using similar techniques. 
\subsection{Setting and main result}
Let $\Omega \subset \BbbR^2$ be a bounded, polygonal Lipschitz domain with vertices
$A_j$, $j=1,\ldots,J$, and interior angles $\omega_j$, $j=1,\ldots,J$. 
We will require the countably normed spaces introduced 
in \cite{babuska-guo86a,melenk02}. These space are designed to capture important 
features of solutions of elliptic partial differential equations posed on 
polygons, namely, analyticity of the solution and the singular behavior 
at the vertices. Their characterization in terms of these countably normed spaces
also permits proving exponential convergence of piecewise polynomial approximation 
on appropriately graded meshes. 

These countably normed spaces are defined with the aid of 
weight functions $\Phi_{p,\overrightarrow{\beta},k}$ that we now define. 
For $\beta \in [0,1)$, $n \in {\mathbb{N}}_0$, $k >0$, 
and $\overrightarrow{\beta} \in [0,1)^J$, we set  
\begin{eqnarray}
\nonumber 
\Phi_{n,\beta,k}(x) &=& 
\min \left\{ 1,\frac{\left\vert x\right\vert }{\min
\left\{ 1,\frac{\left\vert n\right\vert +1}{k+1}\right\} }\right\} ^{n+\beta
}, \\
\label{eq:weight-fct-polygon}
\Phi_{n,\overrightarrow{\beta},k}(x) &=& \prod_{j=1}^J \Phi_{n,\beta_j,k}(x
- A_j).
\end{eqnarray}
Finally, we denote by $H^{1/2}_{pw}(\partial\Omega)$ the space of functions 
whose restrictions of the edges of $\partial\Omega$ are in $H^{1/2}$. 

We furthermore introduce the constant $C_{sol}(k)$ as a suitable norm of the 
solution operator for (\ref{eq:problem}). That is, $C_{sol}(k)$ is such that 
for all $f \in L^2(\Omega)$, $g \in L^2(\partial\Omega)$ and corresponding 
solution $u$ of (\ref{eq:problem}) there holds 
\begin{equation}
\label{eq:def-Csol}
\|u\|_{1,k,\Omega} \leq C_{sol}(k) \left[ \|f\|_{L^2(\Omega)} + \|g\|_{L^2(\partial\Omega)}\right]. 
\end{equation}
We recall that Theorem~\ref{thm:a-priori} gives $C_{sol}(k) = O(k^{5/2})$ for general
polygons and $C_{sol}(k) = O(1)$ by \cite[Prop.~{8.1.4}]{melenk95} for convex polygons. 
Our motivation for using the notation $C_{sol}(k)$ is emphasize in the following theorem
how {\sl a priori} estimates for Helmholtz problems affect the decomposition result: 
\begin{Theorem} 
\label{thm:decomposition-polygon}
Let $\Omega\subset \BbbR^2$ be a polygon with vertices $A_j$, $j=1,\ldots,J$. 
Then there exist constants $C$, $\gamma > 0$, $\overrightarrow{\beta} \in [0,1)^J$ independent of $k \ge k_0$ 
such that for every $f \in L^2(\Omega)$ 
and $g \in H^{1/2}_{pw}(\partial\Omega)$ the solution $u$ of 
(\ref{eq:problem})
can be written as $u = u_{H^2} + u_{\mathcal A}$ with 
\begin{eqnarray*}
k\|u_{H^2}\|_{1,k,\Omega} + \|u_{H^2}\|_{H^2(\Omega)} &\leq&  C C_{f,g} \\
\|u_{\mathcal A}\|_{H^1(\Omega)} &\leq& \left( C_{sol}(k) +1\right) C_{f,g} \\
k \|u_{\mathcal A}\|_{L^2(\Omega)} &\leq& 
\left( C_{sol}(k) +k\right) C_{f,g} \\
\|\Phi_{n,\overrightarrow{\beta},k}\nabla^{n+2} u_{\mathcal A}\|_{L^2(\Omega)} &\leq&  
C (C_{sol}(k) +1) k^{-1} \max\{n,k\}^{n+2} \gamma^n C_{f,g} 
\quad \forall n \in \BbbN_0 
\end{eqnarray*}
with $C_{f,g}:= \|f\|_{L^2(\Omega)} + \|g\|_{H^{1/2}_{pw}(\partial\Omega)}$ and 
$C_{sol}(k)$ introduced in (\ref{eq:def-Csol}). 
\end{Theorem} 
\begin{proof}
The proof is relegated to Section~\ref{sec:proof-thm-decomposition-polygon}.
We mention that the $k$-dependence of our bounds on $\|u_{\mathcal A}\|_{L^2(\Omega)}$
is likely to be suboptimal due to our method of proof. 
\qed
\end{proof}
Theorem~\ref{thm:decomposition-polygon} may be viewed as the analog of 
Lemma~\ref{lemma:decomposition-lemma}, (\ref{item:decomposition-ii}); 
we conclude this section with the analog of Lemma~\ref{lemma:decomposition-lemma}, 
(\ref{item:decomposition-i}): 
\begin{Corollary}
\label{cor:H2-regularity}
Assume the hypotheses of Theorem~\ref{thm:decomposition-polygon}. Then there exist 
constants $C > 0$, $\overrightarrow{\beta} \in [0,1)^J$ independent of $k$ such that for all 
$f \in L^2(\Omega)$, $g \in H^{1/2}_{pw}(\partial\Omega)$
the solution $u$ of (\ref{eq:problem}) satisfies 
$\|u\|_{1,k,\Omega} \leq 
C C_{sol}(k) \left[ \|f\|_{L^2(\Omega)} + \|g\|_{L^2(\partial\Omega)}\right]$ as well as 
$$
\|\Phi_{0,\overrightarrow{\beta},k} \nabla^2 u\|_{L^2(\Omega)} 
\leq C k (C_{sol}(k) + 1) \left[\|f\|_{L^2(\Omega)} + \|g\|_{H^{1/2}_{pw}(\partial\Omega)}\right]. 
$$
\end{Corollary}
\begin{proof}
The estimate for $\|u\|_{1,k,\Omega}$ expresses (\ref{eq:def-Csol}). The estimate for 
the second derivatives of $u$ follows from Theorem~\ref{thm:decomposition-polygon}
since $u = u_{H^2} + u_{\mathcal A}$. 
\qed
\end{proof}

\subsection{Auxiliary results}
Just as in the proof of Lemma~\ref{lemma:decomposition-lemma}, an important
ingredient of the proof of Theorem~\ref{thm:decomposition-polygon} 
are high and low pass filters. The underlying reason is that the Helmholtz
operator $-\Delta - k^2$ acts very differently on low frequency and 
high frequency functions. Here, the dividing line between 
high and low frequencies
is at $O(k)$. For this reason, appropriate high and low pass filters are 
defined and analyzed in Section~\ref{sec:filters}. Furthermore, when applied 
to high frequency functions the Helmholtz operator behaves similarly to 
the Laplacian $-\Delta$ or the modified Helmholtz operator 
$-\Delta + k^2 $. This latter operator, being positive definite, is easier
to analyze and yet provides insight into the behavior of the 
Helmholtz operator restricted to high frequency functions. The modified Helmholtz 
operator will therefore be a tool for the proof of 
Theorem~\ref{thm:decomposition-polygon} and is thus analyzed in 
Section~\ref{sec:modified-helmholtz}. 
\subsubsection{High and low pass filters, auxiliary results}
\label{sec:filters}
For the polygonal domain $\Omega\subset\BbbR^2$ we introduce for $\eta > 1$ the following 
two low and high pass filters in terms of the Fourier transform 
${\mathcal F}$: 
\begin{enumerate}
\item
The low and high pass filters $L_{\Omega, \eta} f : L^2(\Omega) \rightarrow L^2(\Omega)$ 
and $H_{\Omega,\eta} : L^2(\Omega)\rightarrow L^2(\Omega)$ are defined by 
$$
L_{\Omega,\eta} f = ({\mathcal F}^{-1} \chi_{B_{\eta k}(0)} {\mathcal F}(E_\Omega f))|_{\Omega}, 
\qquad 
H_{\Omega,\eta} f = 
({\mathcal F}^{-1} \chi_{\BbbR^2 \setminus B_{\eta k}(0)} {\mathcal F}(E_\Omega f))|_{\Omega};
$$
here, $B_{\eta k}(0)$ is the ball of radius $\eta k$ with center $0$,  the characteristic function 
of a set $A$ is $\chi_{A}$, 
and $E_{\Omega}$ denotes the Stein extension operator of \cite[Chap.~{VI}]{stein70}.
\item 
Analogously, we define 
$L_{\partial \Omega, \eta} f : L^2(\partial \Omega) \rightarrow L^2(\partial \Omega)$ 
and $H_{\partial \Omega,\eta} : L^2(\partial \Omega)\rightarrow L^2(\partial \Omega)$ in
an {\em edgewise} fashion. Specifically, identifying an edge $e$ of $\Omega$ with an interval
and letting $E_e$ be the Stein extension operator for the interval $e\subset\BbbR$ to the 
real line $\BbbR$, we can define with the univariate Fourier transformation 
${\mathcal F}$ the operators 
$L_{e,\eta}$ and $H_{e,\eta}$ by 
$$
L_{e,\eta} g = ({\mathcal F}^{-1} \chi_{B_{\eta k}(0)} {\mathcal F}(E_e g))|_{e}, 
\qquad 
H_{e,\eta} g = 
({\mathcal F}^{-1} \chi_{\BbbR \setminus B_{\eta k}(0)} {\mathcal F}(E_e f))|_{e};
$$
the operators $L_{\partial\Omega,\eta}$ and $H_{\partial\Omega,\eta}$ are then defined 
edgewise by 
$(L_{\partial\Omega,\eta} g)|_e = L_{e,\eta} g$ and 
$(H_{\partial\Omega,\eta} g)|_e = H_{e,\eta} g$ for all edges $e\subset\partial\Omega$. 
\end{enumerate}
These operators provide stable decompositions of $L^2(\Omega)$ and $L^2(\partial\Omega)$. 
For example, one has $L_{\Omega,\eta} + H_{\Omega,\eta} = \operatorname*{Id}$ on $L^2(\Omega)$
and the bounds 
$$
\|L_{\Omega,\eta} f \|_{L^2(\Omega)} + 
\|H_{\Omega,\eta} f \|_{L^2(\Omega)}  
\leq C \|f\|_{L^2(\Omega)} \qquad \forall f \in L^2(\Omega), 
$$
where $C > 0$ depends solely on $\Omega$ (via the Stein extension operator $E_\Omega$). 
The operators $H_{\Omega,\eta}$ and $H_{\partial\Omega,\eta}$ have furthermore approximation
properties if the function they are applied to has some Sobolev regularity. 
We illustrate
this for $H_{\partial\Omega,\eta}$: 
\begin{Lemma}
\label{lemma:filters}
Let $\Omega\subset\BbbR^2 $ be a polygon. Then there exists $C > 0$ independent of 
$k$ and $\eta > 1$ such that for all $g \in H^{1/2}_{pw}(\partial\Omega)$ 
\begin{eqnarray*}
k^{1/2} ( 1 + \eta^{1/2}) \|H_{\partial\Omega,\eta} g\|_{L^2(\partial\Omega)} + 
\|H_{\partial\Omega,\eta} g\|_{H^{1/2}_{pw}(\partial\Omega)} &\leq& C \|g\|_{H^{1/2}_{pw}(\partial\Omega)}. 
\end{eqnarray*}
\end{Lemma}
\begin{proof}
We only show the estimate for $\|H_{\partial\Omega,\eta} g\|_{L^2(\partial\Omega)}$. 
We consider first the case of an interval $I\subset\BbbR$. We define  $H_{I,\eta} g$ by 
$H_{I,\eta} g = {\mathcal F}^{-1} \chi_{\BbbR\setminus B_{\eta k}(0)} {\mathcal F} E_I g$, where  
$\chi_{\BbbR\setminus B_{\eta k}(0)}$ is the characteristic function for 
$\BbbR\setminus (-\eta k,\eta k)$ 
and $E_I$ is the Stein extension operator for the interval $I$. Since, by Parseval,
${\mathcal F}$ is an isometry on $L^2(\BbbR)$ we have 
\begin{eqnarray*}
\lefteqn{
\|H_{I,\eta} g\|^2_{L^2(I)} \leq \|H_{I,\eta} g\|^2_{L^2(\BbbR)} 
= \int_{\BbbR\setminus B_{\eta k}(0)} |{\mathcal F} E_I g|^2\,d\xi 
}\\
&=& \int_{\BbbR\setminus B_{\eta k}(0)} \frac{(1 + |\xi|^2)^{1/2}}{(1 + |\xi|^2)^{1/2}} |{\mathcal F} E_I g|^2\,d\xi 
\leq \frac{1}{(1 + (\eta k)^2)^{1/2}} 
\int_{\BbbR} (1 + |\xi|^2)^{1/2}|{\mathcal F} E_I g|^2\,d\xi. 
\end{eqnarray*}
The last integral can be bounded by $C \|E_I g\|^2_{H^{1/2}(\BbbR)}$. The stability
properties of the extension operator $E_I$ then imply furthermore 
$\|E_I g\|_{H^{1/2}(\BbbR)} \leq C \|g\|_{H^{1/2}(I)}$. In total, 
$$
\|H_{I,\eta} g\|_{L^2(I)} \leq C \frac{1}{(1 + (\eta k)^2)^{1/4}} \|g\|_{H^{1/2}(I)} 
\leq C k^{-1/2} ( 1 + \eta)^{-1/2} \|g\|_{H^{1/2}(I)}, 
$$
where, in the last estimate, the constant $C$ depends additionally 
on $k_0$. From this
estimate, we obtain the desired bound for 
$\|H_{\partial \Omega,\eta} g\|_{L^2(\partial\Omega)}$ by identifying each edge of $\Omega$
with an interval. 
\qed
\end{proof}
\subsubsection{Corner singularities}
We recall the following result harking back to the work
by Kondratiev and Grisvard: 
\begin{Lemma} 
\label{lemma:corner-singularities-laplacian}
Let $\Omega\subset\BbbR^d$ be a polygon with vertices $A_j$, $j=1,\ldots,J$, 
and interior angles $\omega_j$, $j=1,\ldots,J$. Define for each vertex $A_j$ the singularity
function $S_j$ by 
\begin{equation}
\label{eq:Si}
S_j(r_j,\varphi_j) = r_j^{\pi/\omega_j} \cos \left(\frac{\pi}{\omega_j} \varphi_j\right),  
\end{equation}
where $(r_j, \varphi_j)$ are polar coordinates centered at the vertex $A_j$ such that 
the edges of $\Omega$ meeting at $A_j$ correspond to $\varphi_j = 0$ and $\varphi_j = \omega_j$. 
Then every solution $u$ of 
$$
-\Delta u =f \quad \mbox{ in $\Omega$}, 
\qquad \partial_n u  = g \quad \mbox{ on $\partial\Omega$}, 
$$
can be written as $u = u_0 + \sum_{j=1}^J a_j^{\Delta}(f,g) S_j$ with the {\sl a priori} bounds 
\begin{equation}
\label{eq:lemma:corner-singularities-laplacian-10}
\|u_0\|_{H^2(\Omega)} + \sum_{j=1}^J |a_j^\Delta(f,g)| \leq C \left [ \|f\|_{L^2(\Omega)} + 
\|g\|_{H^{1/2}_{pw}(\partial\Omega)} + \|u\|_{H^1(\Omega)}\right]. 
\end{equation}
The $a_j^\Delta$ are linear functionals, and $a_j^\Delta = 0$ for convex corners $A_j$
(i.e., if $\omega_j < \pi$). 
\end{Lemma}
\begin{proof}
This classical result is comprehensively treated in \cite{grisvard85a}. 
\qed
\end{proof}
\subsubsection{The modified Helmholtz equation} 
\label{sec:modified-helmholtz}
We consider the modified Helmholtz equation in both a bounded domain with Robin boundary
conditions and in the full space $\BbbR^2$. The corresponding solution operators will
be denoted $S^+_{\Omega}$ and $S^{+}_{\BbbR^2}$: 
\begin{enumerate}
\item 
The operator $S^+_\Omega: L^2(\Omega) \times H^{1/2}_{pw}(\partial\Omega)\rightarrow H^1(\Omega)$ is
the solution operator for  
\begin{equation}
\label{eq:yukawa}
-\Delta u + k^2 u  = f \quad \mbox{ in $\Omega$}, 
\qquad 
\partial_n u + \bi k u = g \quad \mbox{ on $\partial\Omega$.}
\end{equation}
\item
The operator $S^+_{\BbbR^2}:L^2(\BbbR^2) \rightarrow H^1(\BbbR^2)$ is 
the solution operator for 
\begin{equation}
\label{eq:yukawa-fullspace}
-\Delta u + k^2 u  = f \quad \mbox{ in $\BbbR^2$}. 
\end{equation}
\end{enumerate}
\begin{Lemma}[properties of $S^+_\Omega$]
\label{lemma:S+Omega}
Let $\Omega\subset\BbbR^2$ be a polygon and $f \in L^2(\Omega)$, $g \in H^{1/2}_{pw}(\partial\Omega)$. 
Then  the solution $u:= S^+_\Omega(f,g)$ satisfies 
\begin{eqnarray}
\label{eq:lemma:S+Omega-1}
\|u\|_{1,k,\Omega} & \leq & k^{-1/2} \|g\|_{L^2(\partial\Omega)} + k^{-1} \|f\|_{L^2(\Omega)}. 
\end{eqnarray}
Furthermore, there exists $C > 0$ independent of $k$ and the data $f$, $g$, 
and there exists 
a decomposition $u = u_{H^2} + \sum_{i=1}^J a_i^+(f,g) S_i$ 
for some linear functionals $a_i^+$ with 
\begin{equation}
\label{eq:lemma:S+Omega-2}
\|u_{H^2}\|_{H^2(\Omega)} + \sum_{i=1}^J |a_i^+(f,g)| \leq 
C \left[ \|f\|_{L^2(\Omega)} + \|g\|_{H^{1/2}_{pw}(\partial\Omega)} + k^{1/2} \|g\|_{L^2(\partial\Omega)}\right].
\end{equation}
\end{Lemma}
\begin{proof}
The estimate (\ref{eq:lemma:S+Omega-1}) for $\|u\|_{1,k,\Omega}$ follows by Lax-Milgram -- see 
\cite[Lemma~{4.6}]{melenk-sauter11} for details. 
Since $u$ satisfies 
$$
-\Delta u =  f - k^2 u \quad \mbox{ in $\Omega$}, 
\qquad \partial_n u = g - \bi k u \quad \mbox{ on $\partial\Omega$}, 
$$
the standard regularity theory for the Laplacian 
(see Lemma~\ref{lemma:corner-singularities-laplacian})
permits us to decompose 
$u = u_{H^2} + \sum_{i=1}^J a_i^\Delta(f - k^2 u, g - \bi k u)  S_i$. 
The continuity of the linear functionals $a_i^\Delta$ reads 
$$
\sum_{i=1}^J
|a_i^\Delta(f - k^2 u, g - \bi k u)| \leq C \left[ \|f - k^2 u\|_{L^2(\Omega)} 
                                     + \|g - \bi k u\|_{H^{1/2}_{pw}(\partial\Omega)}\right]. 
$$
Since $(f,g) \mapsto S^+_\Omega(f,g)$ is linear, the map $(f,g) \mapsto a_i^+(f,g):= a_i^\Delta(f - k^2 u, g - \bi k u)$
is linear, and  (\ref{eq:lemma:S+Omega-1}), 
(\ref{eq:lemma:corner-singularities-laplacian-10}) give the desired estimates 
for $u_{H^2}$ and $a_i^+(f,g)$. 
\qed
\end{proof}
\begin{Lemma}[properties of $S^+_{\BbbR^2}$]
\label{lemma:S+R2}
There exists $C > 0$ such that for every $\eta > 1$ and 
every $f \in L^2(\BbbR^2)$ whose Fourier transform ${\mathcal F} f$ 
satisfies $\operatorname*{supp} {\mathcal F} f \subset \BbbR^2 \setminus B_{\eta k}(0)$, 
the solution $u = S^+_{\BbbR^2} f$ of (\ref{eq:yukawa-fullspace}) satisfies 
\begin{eqnarray*}
\|u\|_{1,k,\BbbR^2}  \leq  k^{-1} \frac{1}{\sqrt{1+\eta^2}} \|f\|_{L^2(\BbbR^2)}, 
\qquad \qquad 
\|u\|_{H^2(\BbbR^2)} \leq C \|f\|_{L^2(\BbbR^2)}. 
\end{eqnarray*}
\end{Lemma}
\begin{proof}
The result follows from Parseval's theorem and the weak formulation for $u$ as follows
(we abbreviate the Fourier transforms by $\widehat f = {\mathcal F} f$ and $\widehat u = {\mathcal F} u$):
\begin{align*}
& \|u\|^2_{1,k,\BbbR^2} = (f,u)_{L^2(\BbbR^2)} = (\widehat f,\widehat u)_{L^2(\BbbR^2)} \\
&\leq \sqrt{\int_{\BbbR^2} (|\xi|^2 + k^2)^{-1} |\widehat f|^2\,d\xi} 
     \sqrt{\int_{\BbbR^2} (|\xi|^2 + k^2) |\widehat u|^2\,d\xi} \\
&= 
\sqrt{\int_{\BbbR^2\setminus B_{\eta k}(0)} (|\xi|^2 + k^2)^{-1} |\widehat f|^2\,d\xi} \|u\|_{1,k,\BbbR^2}
\leq \frac{1}{k \sqrt{1+\eta^2}}\|\widehat f\|_{L^2(\BbbR^2)} \|u\|_{1,k,\BbbR^2}, 
\end{align*}
where, in the penultimate step, we used the support properties of $\widehat f$. Appealing
again to Parseval, we get the desired claim for $\|u\|_{1,k,\BbbR^2}$. The estimate
for $\|u\|_{H^2(\BbbR^2)}$ now follows from $f \in L^2(\BbbR^2)$ and the standard 
interior regularity for the Laplacian. 
\qed
\end{proof}
\subsection{Proof of Theorem~\ref{thm:decomposition-polygon}}
\label{sec:proof-thm-decomposition-polygon}
We denote by $S:(f,g) \mapsto S(f,g)$ the solution operator 
for (\ref{eq:problem}). Concerning some of its properties, we  have the 
following result taken essentially from \cite[Lemma~{4.13}]{melenk-sauter11}: 
\begin{Lemma}[analytic regularity of $S(f,g)$]
\label{lemma:analytic-regularity}
Let $\Omega$ be a polygon. Let $f$ be analytic on $\Omega$ and $g \in L^2(\partial\Omega)$
be piecewise analytic and satisfy for some  constants 
$\widetilde C_f$, $\widetilde C_g$, $\gamma_f$, $\gamma_g > 0$ 
\begin{subequations}
\label{eq:analytic-data}
\begin{eqnarray}
\|\nabla^n f\|_{L^2(\Omega)} &\leq& \widetilde C_f \gamma_f^n \max\{n,k\}^n \qquad \forall n \in \BbbN_0\\
\|\nabla_T^n g\|_{L^2(e)} &\leq& \widetilde C_g \gamma_g^n \max\{n,k\}^n \qquad \forall n \in \BbbN_0
\quad \forall e \in {\mathcal E},
\end{eqnarray}
\end{subequations}
where ${\mathcal E}$ denotes the set of edges of $\Omega$ and $\nabla_T$ tangential differentiation. 
Then there exist $\overrightarrow \beta \in [0,1)^J$ 
(depending only on $\Omega$) and constants 
$C$, $\gamma > 0$ (depending only on $\Omega$, $\gamma_f$, $\gamma_g$, $k_0$) 
such that the following is true with the constant $C_{sol}(k)$ 
of (\ref{eq:def-Csol}): 
\begin{eqnarray}
\label{eq:lemma:analytic-regularity-1}
\|u\|_{1,k,\Omega} &\leq& C_{sol}(k) (\widetilde C_f + \widetilde C_g)\\
\label{eq:lemma:analytic-regularity-2}
\|\phi_{n,\overrightarrow\beta,k} \nabla^{n +2} u\|_{L^2(\Omega)} &\leq& 
C C_{sol}(k) k^{-1} (\widetilde C_f +\widetilde C_g) \gamma^n \max\{n,k\}^{n+2} \quad \forall n \in \BbbN_0.
\end{eqnarray}
\end{Lemma}
\begin{proof}
The estimate (\ref{eq:lemma:analytic-regularity-1}) is simply a restatement of 
the definition of $C_{sol}(k)$. The estimate (\ref{eq:lemma:analytic-regularity-2})
will follow from \cite[Prop.~{5.4.5}]{melenk02}. To simplify the presentation, 
we assume by linearity that $g$ vanishes on all edges of $\Omega$ with the exception
of one edge $\Gamma$. Furthermore, we restrict our attention to the vicinity of 
one vertex, which we take to be the origin; we assume $\Gamma \subset (0,\infty) \times\{0\}$,  
and that near the origin, $\Omega$ is above $(0,\infty) \times \{0\}$, i.e., 
$\{(r \cos \varphi,r \sin\varphi)\colon 0 < r < \rho, 0 < \varphi < \omega\} \subset\Omega$
for some $\rho$, $\omega  > 0$.  

Upon setting $\varepsilon :=1/k$, we note that $u$ solves
\begin{equation*}
-\varepsilon ^{2}\Delta u-u=\varepsilon ^{2} f  \quad \text{on }\Omega , 
\qquad 
\varepsilon ^{2}\partial _{n}u=\varepsilon (\varepsilon g-\bi 
u)  \quad \text{on }\partial \Omega .%
\end{equation*}%
On the edge $\Gamma$, the function $g$ is the restriction of 
$G_{1,0}(x,y):= g(x) e^{-y/\varepsilon}$ to $\Gamma$. The assumptions on 
$f$ and $g$ then imply that \cite[Prop.~{5.4.5}]{melenk02} is applicable with
the following choice of constants appearing in \cite[Prop.~{5.4.5}]{melenk02}:
\begin{equation*}
\begin{array}{lllll}
C_{f}=\varepsilon ^{2} \widetilde C_f, & C_{G_{1}}=%
\varepsilon \varepsilon^{1/2} \widetilde C_g, & C_{G_{2}}=%
\varepsilon , & C_{b}=0, & C_{c}=1, \\ 
\gamma _{f}=O(1), & \gamma _{G_{1}}=O(1), & \gamma _{G_{2}}=O(1), & \gamma
_{b}=0, & \gamma _{c}=0,%
\end{array}%
\end{equation*}%
resulting in the existence of constants $C$, $K > 0$ and $\overrightarrow{%
\beta} \in [0,1)^J$ with
\begin{equation*}
\Vert \Phi _{n,\overrightarrow{\beta },k}\nabla ^{n+2} u\Vert
_{L^{2}(\Omega )}\leq C K^{n+2}\max \{n+2,k\}^{n+2}\left( k^{-2} \widetilde C_f 
+k^{-1}\| u\|_{1,k,\Omega} + k^{-3/2} \widetilde C_g \right)
\end{equation*}%
for all $n\in \mathbb{N}_{0}$. 
We conclude the argument by 
inserting (\ref{eq:lemma:analytic-regularity-1}) and estimating generously 
$k^{-1} \widetilde C_f + k^{-1/2} \widetilde C_g \leq 
C\left(\widetilde C_f + \widetilde C_g\right)$. 

We remark that this last generous estimate comes from the precise form of our 
stability assumption (\ref{eq:def-Csol}). Its form (\ref{eq:def-Csol}) is 
motivated by the estimates {\em available} for the star-shaped case, but could
clearly be replaced with other assumptions. 
\qed
\end{proof}
\begin{Corollary}[analytic regularity of $S(L_{\Omega,\eta} f,L_{\partial\Omega,\eta}g)$]
\label{cor:analytic-regularity}
Let $\Omega$ be a polygon and $\eta > 1$. Then there exist 
$\overrightarrow{\beta} \in [0,1)^J$ (depending only on $\Omega$) and 
$C$, $\gamma > 0$ (depending only on $\Omega$, $k_0$, 
and $\eta > 1$) such that 
for every $f \in L^2(\Omega)$ and $g \in L^2(\partial\Omega)$, the function 
$u = S(L_{\Omega,\eta} f,L_{\partial\Omega,\eta} g)$ satisfies 
with $C_{f,g}:= \|f\|_{L^2(\Omega)} + \|g\|_{L^2(\partial\Omega)}$
\begin{eqnarray}
\|u\|_{1,k,\Omega}& \leq & C C_{sol}(k) C_{f,g} \\
\|\Phi_{n,\overrightarrow{\beta},k} \nabla^{n+2}u\|_{L^2(\Omega)}& \leq & C C_{sol}(k) k^{-1} 
\gamma^n \max\{n,k\}^{n+2} C_{f,g} 
\qquad \forall n \in \BbbN_0.
\end{eqnarray}
\end{Corollary}
\begin{proof}
The definitions of $L_{\Omega,\eta}f $ and $L_{\partial\Omega,\eta}$ imply 
with Parseval 
\begin{eqnarray*}
\|\nabla^n L_{\Omega,\eta} f\|_{L^2(\Omega)} &\leq& C (\eta k)^n \|f\|_{L^2(\Omega)} 
\qquad \forall n \in \BbbN_0, \\
\|\nabla_T^n L_{\partial\Omega,\eta} g\|_{L^2(\partial\Omega)} &\leq& C (\eta k)^n 
\|g\|_{L^2(\partial\Omega)} 
\qquad \forall n \in \BbbN_0, 
\end{eqnarray*}
where again $\nabla_T$ is the (edgewise) tangential gradient. The desired estimates 
now follow from Lemma~\ref{lemma:analytic-regularity}. 
\qed
\end{proof}

Key to the proof of Theorem~\ref{thm:decomposition-polygon} is the following
contraction result: 
\begin{Lemma}[contraction lemma]
\label{lemma:contraction-lemma}
Let $\Omega\subset\BbbR^2$ be a polygon. Fix $q \in (0,1)$. Then one can find 
$\overrightarrow\beta \in [0,1)^J$ (depending solely on $\Omega$) and constants 
$C$, $\gamma > 0$ independent of $k$ such that for every $f \in L^2(\Omega)$ and 
every $g \in H^{1/2}_{pw}(\partial\Omega)$, 
the solution $u$ of (\ref{eq:problem}) can be decomposed as 
$
u = u_{H^2} + \sum_{i=1}^J a_i(f,g) S_i + u_{\mathcal A} + r
$, 
where $u_{H^2} \in H^2(\Omega)$, the $a_i$ are linear functionals, 
and $u_{\mathcal A} \in C^\infty(\Omega)$. 
These functions satisfy 
\begin{align*}
& k\|u_{H^2}\|_{1,k,\Omega} + \|u_{H^2}\|_{H^2(\Omega)} + \sum_{i=1}^J |a_i(f,g)| 
 \leq  C \left[ \|f\|_{L^2(\Omega)} + \|g\|_{H^{1/2}_{pw}(\partial\Omega)}\right], \\
& \|u_{\mathcal A}\|_{1,k,\Omega} 
\leq C C_{sol}(k) \left[ \|f\|_{L^2(\Omega)} + \|g\|_{L^2(\partial\Omega)}\right],\\
& \|\Phi_{n,\overrightarrow{\beta},k} \nabla^{n+2} u_{\mathcal A}\|_{L^2(\Omega)} 
\leq C C_{sol}(k) k^{-1} \gamma^n \max\{n,k\}^{n+2} 
\left[ \|f\|_{L^2(\Omega)} + \|g\|_{L^2(\partial\Omega)}\right] 
\end{align*}
for all $n \in \BbbN_0$.  Finally, the remainder $r$ satisfies 
\begin{eqnarray*}
-\Delta r - k^2 r = \widetilde f \quad \mbox{ on $\Omega$}, 
\qquad \partial_n r + \bi k r = \widetilde g
\end{eqnarray*}
for some $\widetilde f \in L^2(\Omega)$ and $\widetilde g \in H^{1/2}_{pw}(\partial\Omega)$
with 
$$
\|\widetilde f\|_{L^2(\Omega)} + \|\widetilde g\|_{H^{1/2}_{pw}(\partial\Omega)} \leq q 
\left(\|f\|_{L^2(\Omega)} + \|g\|_{H^{1/2}_{pw}(\partial\Omega)}\right).
$$
\end{Lemma}
\begin{proof}
We start by decomposing 
$(f,g) = (L_{\Omega,\eta} f, L_{\partial\Omega,\eta} g) 
+ (H_{\Omega,\eta} f, H_{\partial\Omega,\eta} g)$ with a parameter $\eta > 1$ that will be 
selected below. We set 
$$
u_{\mathcal A}:= S(L_{\Omega,\eta}f,L_{\partial\Omega,\eta} g), 
\qquad 
u_{1}:= S^+_{\BbbR^2}(H_{\Omega,\eta} f), 
$$
where we tacitly extended $H_{\Omega,\eta} f$ (which is only defined on $\Omega$) by zero outside
$\Omega$. 
Then $u_{\mathcal A}$ satisfies the desired estimates
by Corollary~\ref{cor:analytic-regularity}. 
For  $u_1$ we have 
by Lemma~\ref{lemma:S+R2} and the stability $\|H_{\Omega,\eta} f\|_{L^2(\Omega)} \leq C \|f\|_{L^2(\Omega)}$
(we note that $C > 0$ is independent of $k$ and $\eta$) the {\sl a priori} estimates
\begin{eqnarray*}
\|u_1\|_{1,k,\BbbR^2} &\leq&  C k^{-1} (1+\eta^2)^{-1/2} \|H_{\Omega,\eta} f\|_{L^2(\Omega)}
                     \leq C k^{-1} (1+\eta)^{-1} \|f\|_{L^2(\Omega)},\\
\|u_1\|_{H^2(\BbbR^2)} &\leq&  C \|H_{\Omega,\eta} f\|_{L^2(\Omega)}
                     \leq C \|f\|_{L^2(\Omega)}. 
\end{eqnarray*}
The trace and the multiplicative trace inequalities imply for 
$g_1:= \partial_n u_1 + \bi k u_1$: 
\begin{eqnarray*}
k^{1/2} (1 + \eta)^{1/2} \|g_1\|_{L^{2}(\partial\Omega)} + 
\|g_1\|_{H^{1/2}_{pw}(\partial\Omega)} &\leq& C \|f\|_{L^2(\Omega)}.  
\end{eqnarray*}
For $g_2:= H_{\partial\Omega,\eta} g - g_1$ we then get from 
Lemma~\ref{lemma:filters} and the triangle inequality
\begin{eqnarray*}
k^{1/2} (1+\eta)^{1/2} \|g_2\|_{L^2(\partial\Omega)} + \|g_2\|_{H^{1/2}_{pw}(\partial\Omega)} 
&\leq & C \left[ \|g\|_{H^{1/2}_{pw}(\partial\Omega)} + \|f\|_{L^2(\Omega)} \right]. 
\end{eqnarray*}
Lemma~\ref{lemma:S+Omega} yields
for $u_2:= S^+_\Omega(0, g_2)$, 
$$
\|u_2\|_{1,k,\Omega} \leq C k^{-1/2} \|g_2\|_{L^2(\partial\Omega)} 
\leq C k^{-1} (1 + \eta)^{-1/2} 
\left[ \|f\|_{L^2(\Omega)}  + \|g\|_{H^{1/2}_{pw}(\partial\Omega)}\right], 
$$
and furthermore we can write 
$
u_2 = u_{H^2} + \sum_{i=1}^J a_i^+(0,g_2) S_i, 
$
with 
$$
\|u_{H^2}\|_{H^2(\Omega)} + \sum_{i=1}^J |a_i^+(0,g_2)| \leq 
C \left[\|f\|_{L^2(\Omega)} + \|g\|_{H^{1/2}_{pw}(\partial\Omega)}\right]. 
$$
We then define $a_i(f,g):= a_i^+(0,g_2)$ and note that $(f,g) \mapsto a_i(f,g)$ 
is linear by linearity of the maps $a_i^+$ and $(f,g) \mapsto g_2$. 
The above shows that $u_{H^2}$ and the $a_i$ satisfy the required estimates. 
Finally, the function $\widetilde u:= u - (u_{\mathcal A} + u_1 + u_2)$ satisfies 
$$
-\Delta \widetilde u - k^2 \widetilde u = 2 k^2 (u_1+u_2) =:\widetilde f, 
\qquad 
\partial_n \widetilde u + \bi k \widetilde u = 0 =:\widetilde g, 
$$
together with 
$$
\|\widetilde f\|_{L^2(\Omega)} \leq 
2 k^2 \left( \|u_1\|_{L^2(\Omega)} + \|u_2\|_{L^2(\Omega)}\right) 
\leq C (1 + \eta)^{-1/2} \left[ \|f\|_{L^2(\Omega)} + \|g\|_{H^{1/2}_{pw}(\partial\Omega)}\right].
$$
Hence, selecting $\eta > 1$ sufficiently large so that for the chosen $q \in (0,1)$ we have 
$C (1+ \eta)^{-1/2} \leq q$ allows us to conclude the proof. 
\qed
\end{proof}
\begin{numberedproof}{Theorem~\ref{thm:decomposition-polygon}}
The contraction property of Lemma~\ref{lemma:contraction-lemma}
can be restated as 
$S(f,g) = u_{H^2} + \sum_{i=1}^J a_i(f,g) S_i + u_{\mathcal A} + S(\widetilde f,\widetilde g)$,
where, for a chosen $q \in (0,1)$, we have 
$\|\widetilde f\|_{L^2(\Omega)} + \|\widetilde g\|_{H^{1/2}_{pw}(\partial\Omega)} \leq 
q \left[ \|f\|_{L^2(\Omega)} + \|g\|_{H^{1/2}_{pw}(\partial\Omega)}\right]$ together with appropriate
estimates for $u_{H^2}$, $a_i(f,g)$, and $u_{\mathcal A}$. This consideration can be repeated 
for $S(\widetilde f,\widetilde g)$. We conclude that a geometric
series argument can be employed to write 
$u = S(f,g) = u_{H^2} + \sum_{i=1}^J \widetilde a_{i}(f,g) S_i + \widetilde u_{\mathcal A}$, 
where $u_{H^2} \in H^2(\Omega)$, $\widetilde u_{\mathcal A} \in C^\infty(\Omega)$, 
and the coefficients $\widetilde a_i$ are in fact linear functionals of the data $(f,g)$. 
Furthermore, we have with the abbreviation 
$C_{f,g} := \|f\|_{L^2(\Omega)} + \|g\|_{H^{1/2}_{pw}(\partial\Omega)}$ 
\begin{align*}
& \|\widetilde u_{\mathcal A}\|_{1,k,\Omega} \leq C C_{f,g} \\
& \|\Phi_{n,\overrightarrow\beta,k} \nabla^{n+2} \widetilde u_{\mathcal A}\|_{L^2(\Omega)} 
\leq C C_{sol}(k)k^{-1} C_{f,g} \gamma^n \max\{n,k\}^{n+2} \qquad \forall n \in \BbbN_0,\\
& k \|u_{H^2}\|_{1,k,\Omega}+\|u_{H^2}\|_{H^2(\Omega)} 
+ \sum_{i=1}^{J}|\widetilde a_{i}(f,g)| \leq C C_{f,g}. 
\end{align*}
Finally, Lemma~\ref{lemma:bbeta-regularity-singularity-fcts} below allows us to 
absorb the contribution $\sum_{i=1}^J \widetilde a_{i}(f,g) S_i$ in the analytic part
by setting $u_{\mathcal A}:= \widetilde u_{\mathcal A} 
+ \sum_{i=1}^J \widetilde a_{i}(f,g) S_i$. 
In view of $\beta_i < 1$, we have $2 - \beta_i \ge 1$ and arrive at 
\begin{eqnarray*}
\|u_{\mathcal A}\|_{H^1(\Omega)} &\leq& C (C_{sol}(k) + 1)C_{f,g}, 
\qquad \
k \|u_{\mathcal A}\|_{L^2(\Omega)} \leq C C_{f,g} (C_{sol}(k) + k), 
\\
\|\Phi_{n,\overrightarrow\beta,k} \nabla^{n+2} u_{\mathcal A}\|_{L^2(\Omega)} 
& \leq & C C_{f,g} \left[C_{sol}(k) k^{-1} + k^{-1} \right] 
\max\{n,k\}^{n+2} \quad \forall n \in \BbbN_0,
\end{eqnarray*}
which concludes the argument. 
\end{numberedproof}
\begin{Lemma}
\label{lemma:bbeta-regularity-singularity-fcts}
Let $\beta_i \in [0,1)$ satisfy $\beta_i > 1-\frac{\pi}{\omega_i}$. Then, 
for some $C$, $\gamma > 0$ independent of $k$, 
the singularity functions $S_i$ of (\ref{eq:Si}) satisfy
$\|S_i\|_{H^1(\Omega)} \leq C$ and 
\begin{eqnarray*}
\|\Phi_{n,\overrightarrow\beta,k} \nabla^{n+2} S_i\|_{L^2(\Omega)} &\leq &C k^{-(2-\beta_i)} \gamma^n \max\{n,k\}^{n+2} 
\quad \forall \in \BbbN_0 
\end{eqnarray*}
\end{Lemma}
\begin{proof}
Follows by a direct calculation. 
\iftechreport 
See Lemma~\ref{lemma:singular-fcts} for details.
\fi
\qed
\end{proof}
\section{Stability of Galerkin discretizations}
\label{sec:discrete-stability-H1-conforming}
\subsection{Abstract results}
We consider the model problem (\ref{eq:problem}) and a sequence 
$(V_N)_{N \in \BbbN} \subset H^1(\Omega)$ of finite-dimensional spaces. 
Furthermore, we assume  that $(V_N)_{N \in \BbbN}$ is such that 
for every $v \in H^1(\Omega)$ we have 
$\lim_{N \rightarrow \infty} \inf_{v_N \in V_N} \|v - v_N\|_{H^1(\Omega)} = 0$. 
The conforming approximations $u_N$ to the solution $u$ of (\ref{eq:problem}) are
then defined by: 
\begin{equation}
\label{eq:problem-weak-discrete}
\mbox{ Find $u_N \in V_N$ s.t.} \quad B(u_N,v) = l(v) \qquad \forall v \in V_N. 
\end{equation}
Since the sesquilinear form $B$ satisfies a G{\aa}rding inequality, general functional
analytic argument show that {\em asymptotically}, the discrete problem 
(\ref{eq:problem-weak-discrete}) has a unique solution $u_N$ and are quasi-optimal
(see, e.g., \cite[Thm.~{4.2.9}]{sauter-schwab10}, \cite{schatz74}).  
More precisely, there exist $N_0 > 0$ and $C > 0$ such that for all $N \ge N_0$ 
\begin{equation}
\label{eq:quasi-optimality-soft}
\|u - u_N\|_{1,k,\Omega} \leq  C \inf_{v \in V_N} \|u - v\|_{1,k,\Omega}. 
\end{equation}
This general functional analytic argument does not give any indication
of how $C$ and $N_0$ depend on discretization parameters and the 
wavenumber $k$. Inspection of the arguments reveals that it is 
the approximation properties of the spaces $V_N$ for the 
approximation of the solution of certain adjoint problems that leads
to the quasi-optimality result (\ref{eq:quasi-optimality-soft}). 
For the reader's convenience, we repeat the 
argument, which has been used previously in, e.g.,  
\cite{schatz74,aziz-kellogg-stephens88,melenk95,melenk-sauter10,melenk-sauter11,banjai-sauter07,sauter05}
and is often attributed to Schatz, \cite{schatz74}:
\begin{Lemma}[\protect{\cite[Thm.~{3.2}]{melenk-sauter11}}]
\label{lemma:abstract-quasi-optimality}
Let $\Omega\subset\BbbR^d$ be a bounded Lipschitz domain and $B$ be defined in (\ref{eq:B}). 
Denote by $S^\star:L^2(\Omega) \rightarrow H^1(\Omega)$ the solution operator for the problem
\begin{equation}
\label{eq:adjoint-problem}
\mbox{ Find $u^\star \in H^1(\Omega)$ s.t.} \qquad 
B(v,u^\star) = (v,f)_{L^2(\Omega)} \qquad \forall v \in H^1(\Omega). 
\end{equation}
Define the {\em adjoint approximation property} $\eta(V_N)$ by 
$$
\eta(V_N):= \sup_{f \in L^2(\Omega)} \inf_{v \in V_N} 
\frac{\|S^\star(f) - v\|_{1,k,\Omega}}{\|f\|_{L^2(\Omega)}}.
$$
If, for the continuity constant $C_B$ of (\ref{eq:CB}),  the space $V_N$ 
satisfies 
\begin{equation}
\label{eq:cond-for-stability}
2 C_B k \eta(V_N) \leq 1, 
\end{equation}
then the solution $u_N$ of  (\ref{eq:problem-weak-discrete}) exists and satisfies 
\begin{equation}
\label{eq:lemma:abstract-quasi-optimality-10}
\|u - u_N\|_{1,k,\Omega} \leq 2 C_B \inf_{v \in V_N} \|u - v\|_{1,k,\Omega}. 
\end{equation}
\end{Lemma}
\begin{proof}
We will not show existence of $u_N$ but restrict our attention on the quasi-optimality result
(\ref{eq:lemma:abstract-quasi-optimality-10}); we refer
to \cite[Thm.~{3.9}]{melenk-loehndorf10} for the demonstration that 
(\ref{eq:lemma:abstract-quasi-optimality-10}) in fact implies existence and uniqueness of $u_N$. 
Letting $e = u - u_N$ be the error, we start with an estimate for $\|e\|_{L^2(\Omega)}$: 
Using the definition of the operator $S^\star$ and the Galerkin orthogonality satisfied 
by $e$, we have for arbitrary $v \in V_N$
$$
\|e\|^2_{L^2(\Omega)} = (e,e)_{L^2(\Omega)} = B(e,S^\star e) = B(e,S^\star e - v) 
\leq C_B \|e\|_{1,k,\Omega} \|S^\star e - v\|_{1,k,\Omega}. 
$$
Infimizing over all $v \in V_N$ yields 
with the adjoint approximation property $\eta(V_N)$
$$
\|e\|_{L^2(\Omega)} \leq C_B \eta(V_N) \|e\|_{1,k,\Omega}. 
$$
The G{\aa}rding inequality and the Galerkin orthogonality yield for arbitrary $v \in V_N$: 
\begin{eqnarray*}
\|e\|^2_{1,k,\Omega} &=& \operatorname*{Re} B(e,e) + 2 k^2 \|e\|^2_{L^2(\Omega)} 
= \operatorname{Re} B(e,u - v) + 2 k^2 \|e\|^2_{L^2(\Omega)} \\
&\leq& C_B \|e\|_{1,k,\Omega} \|u - v\|_{1,k,\Omega} + \left( C_B k \eta(V_N)\right)^2\|e\|^2_{1,k,\Omega}. 
\end{eqnarray*}
The assumption $C_B k \eta(V_N) \leq 1/2$ allows us to rearrange this bound to get 
$\|e\|_{1,k,\Omega} \leq 2 C_B \|u -v\|_{1,k,\Omega}$. Since $v \in V_N$ is arbitrary, we arrive at 
(\ref{eq:lemma:abstract-quasi-optimality-10}). 
\qed
\end{proof} 
Lemma~\ref{lemma:abstract-quasi-optimality} informs us that the convergence analysis
for the Galerkin discretization of (\ref{eq:problem}) can be reduced to the study
of the adjoint approximation property $\eta(V_N)$, which is purely a question of 
approximation theory. In the context of piecewise polynomial approximation spaces $V_N$
this requires a good regularity theory for the operator $S^\star$. The strong form 
of the equation satisfied by $u^\star := S^\star f$ is
\begin{equation}
\label{eq:problem-adjoint}
-\Delta u^\star - k^2 u^\star = f 
\quad \mbox{ in $\Omega$}, 
\qquad \partial_n u^\star - \bi k u^\star = 0 
\quad \mbox{ on $\partial\Omega$},
\end{equation}
which is again a Helmholtz problem of the type considered in 
Section~\ref{sec:regularity}. More formally, with the solution operator $S$ of 
Section~\ref{sec:regularity}, we have $S^\star f = \overline{S(\overline{f},0)}$, where 
an overbar denotes complex conjugation. Thus, the regularity theory of Section~\ref{sec:regularity}
is applicable. 

\subsection{Stability of $hp$-FEM}
\label{sec:stability-of-hpFEM}
The estimates of Theorem~\ref{thm:decomposition-polygon} suggest that the effect 
of the corner singularities is essentially restricted to an $O(1/k)$-neighborhood 
of the vertices. This motivates us to consider meshes that are refined in a small
neighborhood of the vertices. To fix ideas, we restrict our attention to meshes 
${\mathcal T}^{geo}_{h,L}$ that are obtained in the following way: 
First, a quasi-uniform triangulation ${\mathcal T}_h$ with mesh size $h$ is selected. 
Then, the elements abutting the vertices $A_j$, $j=1,\ldots,J$, are refined further 
with a mesh that is geometrically graded towards these vertices. These geometric
meshes have $L$ layers and use a grading factor $\sigma \in (0,1)$
(see \cite[Sec.~{4.4.1}]{schwab98} for a precise formal definition). 
Furthermore, for any regular, shape-regular mesh ${\mathcal T}$, we define 
\begin{equation}
\label{eq:Sp}
S^{p}({\mathcal T}):= \{u \in H^1(\Omega)\colon u|_K \in {\mathcal P}_p \qquad \forall K \in {\mathcal T}\},
\end{equation}
where ${\mathcal P}_p$ denotes the space of polynomials of degree $p$. 
We now show that on the geometric meshes ${\mathcal T}^{geo}_{h,L}$, stability 
of the FEM is ensured if the mesh size $h$ and the polynomial degree $p$ satisfy 
the scale resolution condition (\ref{eq:scale-resolution}) and, additionally,  
$L = O(p)$ layers of geometric refinement are used near the vertices: 
\begin{Theorem}[quasi-optimality of $hp$-FEM]
\label{thm:discrete-stability-polygon}
Let ${\mathcal T}_{h,L}^{geo}$ denote the geometric meshes on the polygon $\Omega\subset\BbbR^2$
as described above. 
Fix $c_3 > 0$. Then there are constants $c_1$, $c_2 > 0$ depending solely on $\Omega$
and the shape-regularity of the mesh ${\mathcal T}_{h,L}^{geo}$ such that the following
is true: If $h$, $p$, and $L$ satisfy the conditions
\begin{equation}
\label{eq:scale-resolution-polygon}
\frac{kh}{p} \leq c_1 
\quad \mbox{ and } \quad 
p \ge c_2 \log k 
\quad \mbox{ and } \quad 
L \ge c_3 p
\end{equation}
then the $hp$-FEM based on the space $S^{p}({\mathcal T}_{h,L}^{geo})$ has a unique
solution $u_N \in S^{p}({\mathcal T}_{h,L}^{geo})$ and 
\begin{equation}
\label{eq:quasi-optimality-estimate}
\|u - u_N\|_{1,k,\Omega} 
\leq 2 C_B \inf_{v \in S^{p}({\mathcal T}_{h,L}^{geo})} \|u - v\|_{1,k,\Omega}. 
\end{equation}
\end{Theorem}
\begin{proof}
By Lemma~\ref{lemma:abstract-quasi-optimality}, we have to estimate 
$k \eta(V_N)$ with $V_N = S^{p}({\mathcal T}^{geo}_{h,L})$. Recalling the definition of 
$\eta(V_N)$ we let $f \in L^2(\Omega)$ and observe that we can decompose 
$S^\star f = u_{H^2} + u_{\mathcal A}$, where $u_{H^2}$ and $u_{\mathcal A}$ satisfy the 
bounds 
\begin{eqnarray*} 
\|u_{H^2} \|_{H^2(\Omega)} &\leq& C \|f\|_{L^2(\Omega)}, \\
\|\Phi_{n,\overrightarrow{\beta},k} \nabla^{n+2} u_{\mathcal A} \|_{L^2(\Omega)} 
&\leq& C (C_{sol}(k)+1) k^{-1} \gamma^n \max\{k,n\}^{n+2} \|f\|_{L^2(\Omega)} 
\qquad \forall n \in \BbbN_0.  
\end{eqnarray*}
Piecewise polynomial approximation on ${\mathcal T}^{geo}_{h,L}$ as 
discussed in \cite[Prop.~{5.6}]{melenk-sauter11} gives 
under the assumptions $kh/p\leq C$ and $L \ge c_3 p$:  
(inspection of the proof of \cite[Prop.~{5.6}]{melenk-sauter11} shows that only
bounds on the derivatives of order $\ge 2$ are needed): 
\begin{eqnarray*}
\inf_{v \in V_N} \|u_{H^2} - v\|_{1,k,\Omega} & \leq & C \frac{h}{p} \|f\|_{L^2(\Omega)}, \\
\inf_{v \in V_N} \|u_{\mathcal A} - v\|_{1,k,\Omega} & \leq & C 
\left[ (kh)^{1-\beta_{max}} e^{c kh - b p} + \left(\frac{kh}{\sigma_0 p}\right)^{p} 
\right]
(C_{sol}(k) +1)\|f\|_{L^2(\Omega)},
\end{eqnarray*}
where $\beta_{max} = \max_{j=1,\ldots,J} \beta_j < 1$, 
and $C$, $c$, $b > 0$ are constants independent of $h$, $p$, and $k$. From 
this, we can easily infer 
$$
k \eta(V_N) \leq C \left\{ \frac{kh}{p} + 
k (C_{sol}(k) +1)\left[ (kh)^{1-\beta_{max}} e^{c kh - b p} + \left(\frac{kh}{\sigma_0 p}\right)^{p} 
\right]\right\}.
$$
Noting that Theorem~\ref{thm:a-priori} gives $C_{sol}(k) = O(k^{5/2})$, 
and selecting $c_1$ sufficiently small as well as 
$c_2$ sufficient large allows us to make $k \eta(V_N)$ so small that 
the condition (\ref{eq:cond-for-stability}) 
in Lemma~\ref{lemma:abstract-quasi-optimality} is satisfied. 
\qed
\end{proof}
\begin{Corollary}[exponential convergence on geometric meshes]
Let $f$ be analytic on $\overline{\Omega}$ and $g$ be piecewise analytic, i.e., $f$, $g$
satisfy (\ref{eq:analytic-data}). Given $c_3 > 0$, there exist $c_1$, $c_2 > 0$ such that 
under the scale resolution conditions (\ref{eq:scale-resolution-polygon}) 
of Theorem~\ref{thm:discrete-stability-polygon}, the 
finite element approximation 
$u_N \in S^p({\mathcal T}^{geo}_{h,L})$ exists, and there are 
constants $C$, $b > 0$ independent of $k$ such that the error $u - u_N$ satisfies 
$$
\|u - u_N\|_{1,k,\Omega} \leq C e^{-b p}. 
$$ 
\end{Corollary}
\begin{proof}
In view of Theorem~\ref{thm:discrete-stability-polygon}, estimating $\|u - u_N\|_{1,k,\Omega}$ 
is purely a question of approximability for $c_1$ sufficiently small 
and $c_2$ sufficiently large. Lemma~\ref{lemma:analytic-regularity} gives that the solution 
$u = S(f,g)$ satisfies the bounds given there and, as in the proof of 
Theorem~\ref{thm:discrete-stability-polygon}, we conclude from 
\cite[Prop.~{5.6}]{melenk-sauter11} (more precisely, this follows from its proof) 
\begin{eqnarray*}
\inf_{v \in V_N} \|u_{\mathcal A} - v\|_{1,k,\Omega} & \leq & C 
\left[ (kh)^{1-\beta_{max}} e^{c kh - b p} + \left(\frac{kh}{\sigma_0 p}\right)^{p} 
\right]
(C_{sol}(k) +1)(\widetilde C_f + \widetilde C_g). 
\end{eqnarray*}
Theorem~\ref{thm:a-priori} asserts $C_{sol}(k) = O(k^{5/2})$, which implies the result
by suitably adjusting $c_1$ and $c_2$ if necessary.
\qed
\end{proof}
\begin{Remark}
{\rm 
\begin{enumerate}
\item 
The problem size $N = \operatorname*{dim} S^{p}({\mathcal T}^{geo}_{h,L})$ is 
$N = O((L+h^{-2})p^2)$. The particular choice of $L = c_3 p$ layers of geometric refinement,
approximation order  $p = c_2 \log k$, 
and mesh size $h = c_1 p/k$ in Theorem~\ref{thm:discrete-stability-polygon} ensures quasi-optimality 
of the $hp$-FEM with problem size $N = O(k^2)$, i.e., quasi-optimality can be achieved 
with a fixed number of degrees of freedom per wavelength. 
\item 
The sparsity pattern of the system matrix is that of the classical $hp$-FEM, 
i.e., each row/column has $O(p^2)$ non-zero entries. Noting that the 
scale resolution conditions 
(\ref{eq:scale-resolution}), (\ref{eq:scale-resolution-polygon})
require $p = O(\log k)$, we see that the number of non-zero entries 
entries per row/column is not independent of $k$. It is worth relating 
this observation to \cite{babuska-sauter00}. It is shown there 
for a model problem in 2D that {\em no} 9 point stencil can be found that 
leads to a pollution-free method. 
\item 
Any space $V_N$ that contains $S^{p}({\mathcal T}_{h,L}^{geo})$, where $h$, $p$, and $L$
satisfy the scale resolution condition (\ref{eq:scale-resolution-polygon})  also features
quasi-optimality. 
\item 
The factor $2$ on the right-hand side of 
(\ref{eq:quasi-optimality-estimate}) is arbitrary and can be replaced by any number greater
than $1$. 
\item
The stability analysis of Theorem~\ref{thm:discrete-stability-polygon} requires 
quite a significant mesh refinement near the vertices, namely, $L \sim p$. It is not clear 
whether this is an artifact of the proof. For a more careful numerical analysis of this 
issue, more detailed information about the stability properties of the solution operator $S$ 
is needed, e.g., estimates for $\|S(f,g)\|_{1,k,B_{1/k}(A_j)}$. 
\end{enumerate} 
}
\end{Remark}
\subsection{Numerical examples: $hp$-FEM}
\label{sec:numerics}
\begin{figure}
\includegraphics[width=0.5\textwidth]{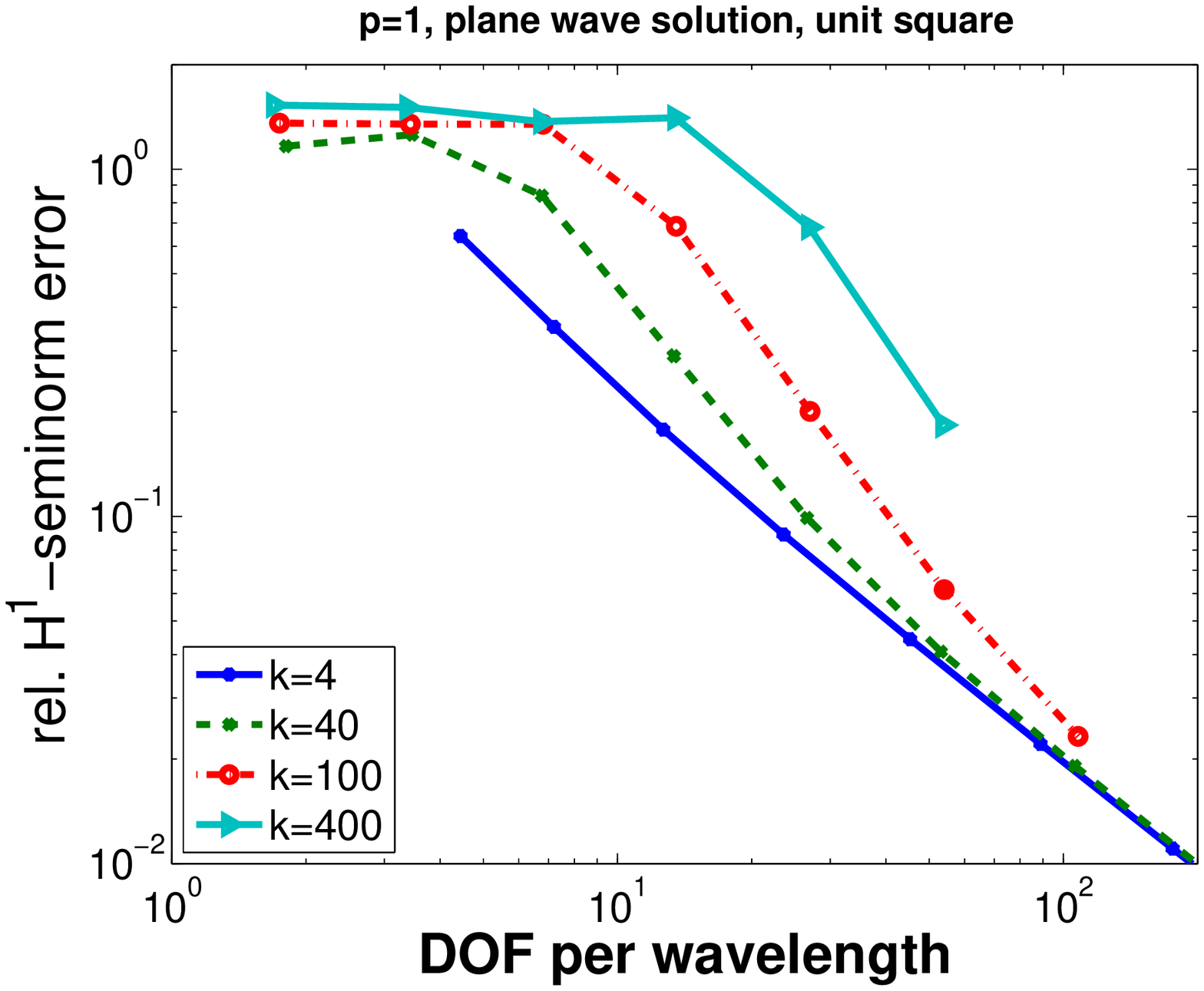}
\includegraphics[width=0.5\textwidth]{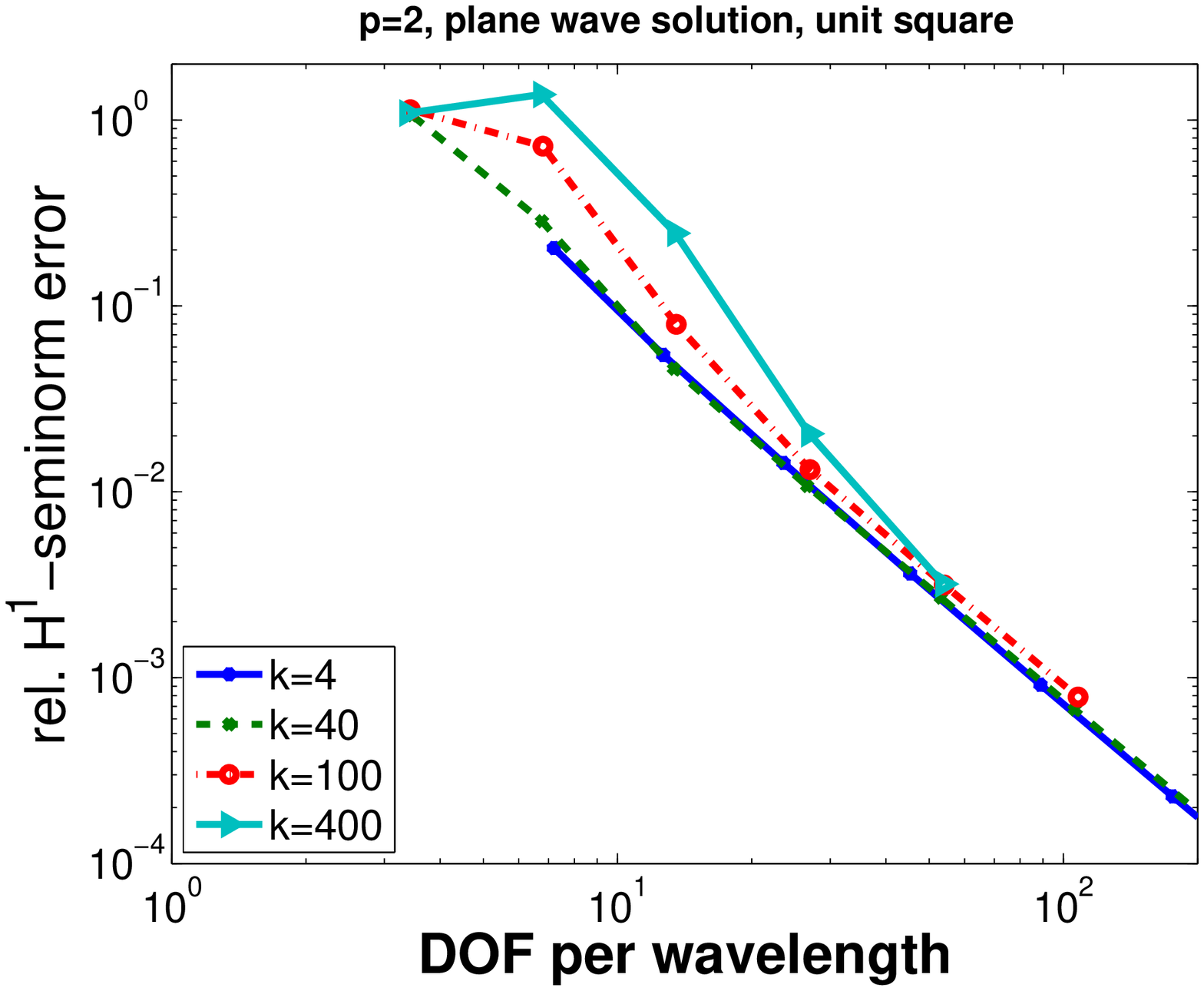}
\includegraphics[width=0.5\textwidth]{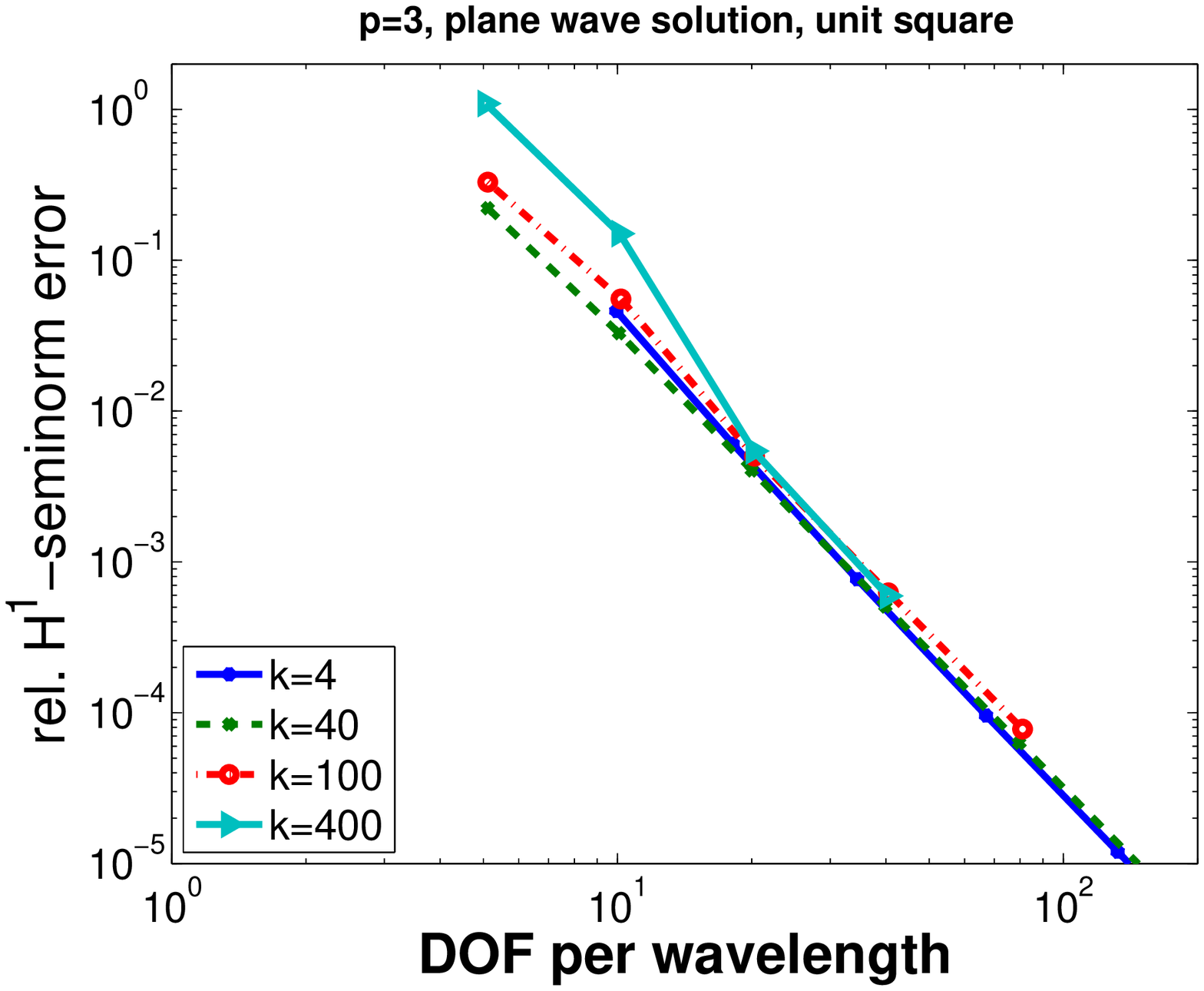}
\includegraphics[width=0.5\textwidth]{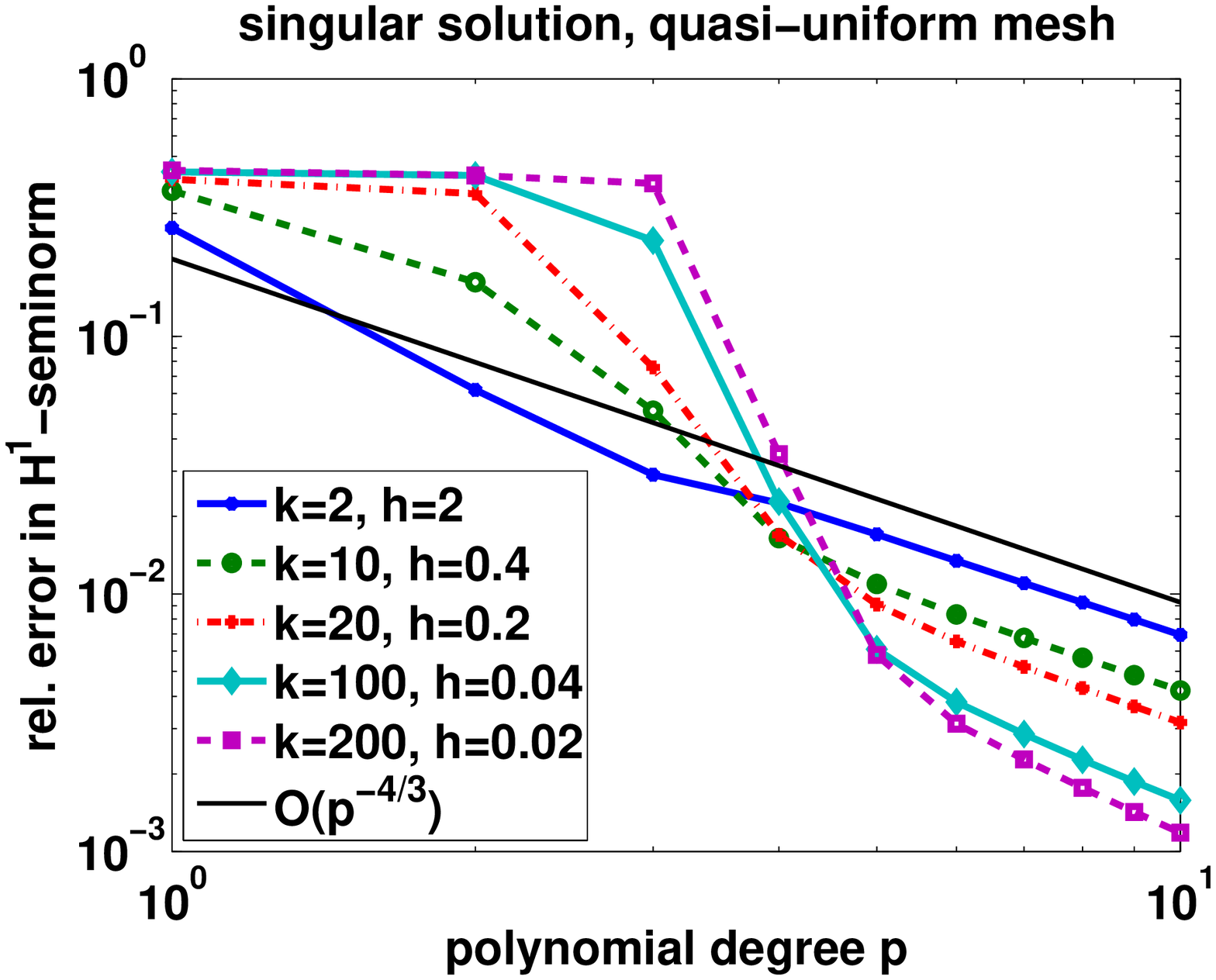}
\caption{\label{fig:square} Top: $h$-FEM with $p=1$ (left) and $p=2$ (right) as described 
in Example~\ref{ex:square}. Bottom left: $h$-FEM with $p=3$ as described in Example~\ref{ex:square}. 
Bottom right: $p$-FEM for singular solution on quasi-uniform mesh as described in Example~\ref{ex:Lshaped-singular}.}
\end{figure}

\begin{figure}
\includegraphics[width=0.5\textwidth]{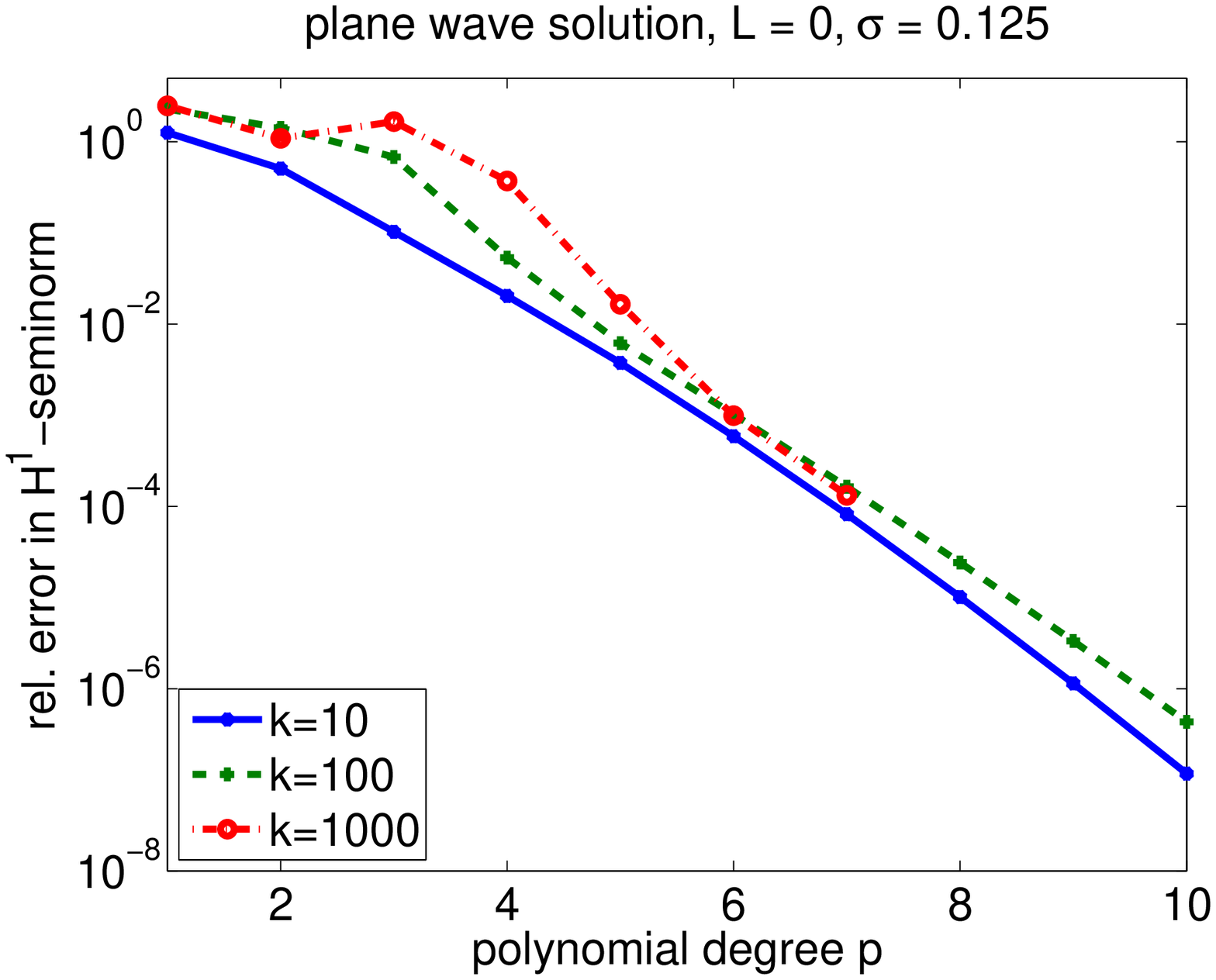}
\includegraphics[width=0.5\textwidth]{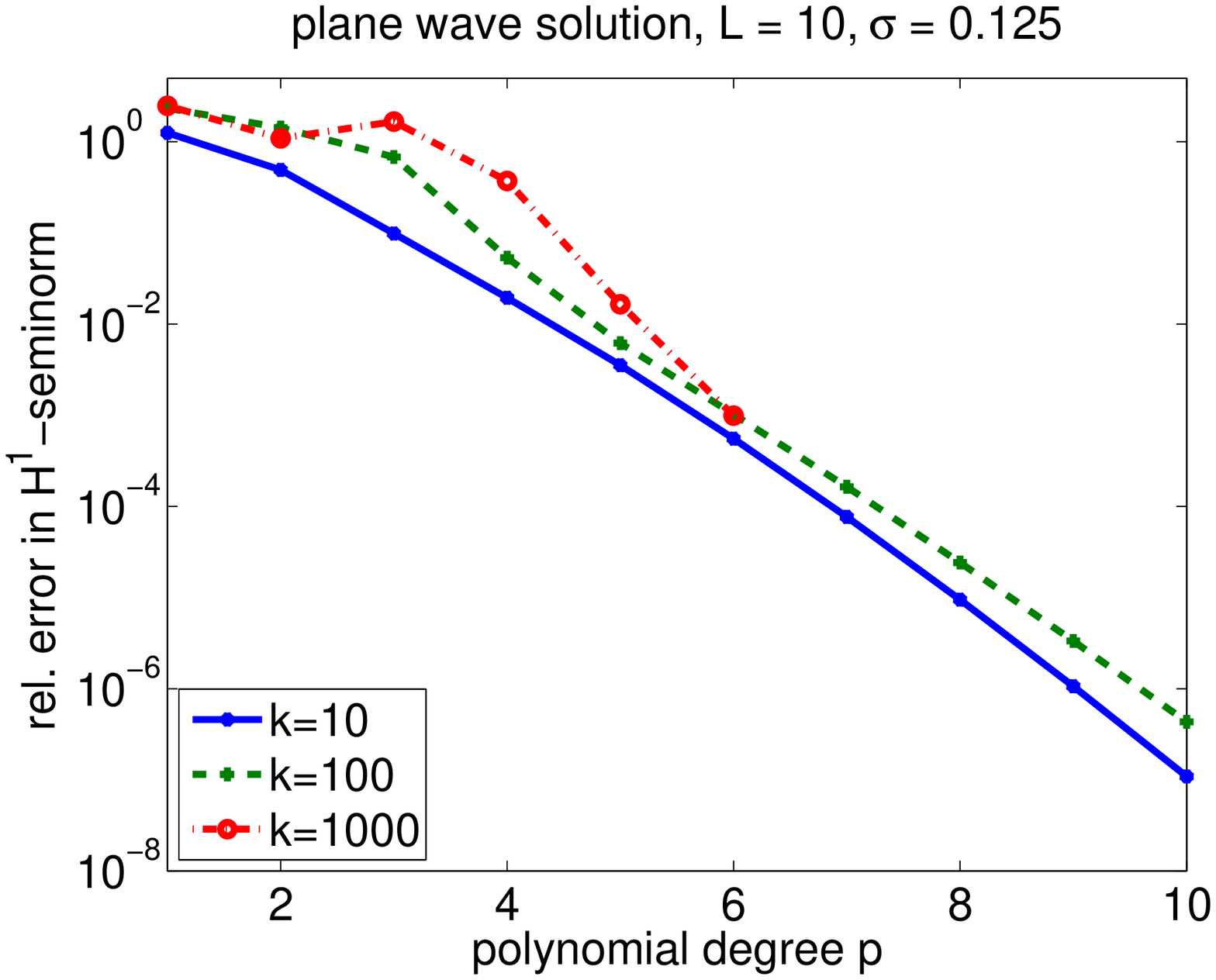}
\caption{\label{fig:numerics} $p$-FEM for plane wave solution as described in 
Example~\ref{ex:Lshaped-smooth}. Left: quasiuniform mesh ${\mathcal T}_h$ with $kh \approx 4$.  
Right: Mesh ${\mathcal T}^{geo}$ obtained from ${\mathcal T}_h$ 
by strong geometric refinement near origin.}
\end{figure}

All calculations reported in this section are performed with the $hp$-FEM
code {\sc netgen/ngsolve} by J.~Sch\"oberl, \cite{schoeberl97,schoeberl-ngsolve}. 
\begin{Example}
\label{ex:square}
{\rm 
In this 2D analog of Example~\ref{ex:1d-pollution}, we consider the 
model problem (\ref{eq:problem}) with exact solution being a plane wave 
$e^{\bi (k_1 x + k_2 y)}$, where $k_1 = - k_2 = \frac{1}{\sqrt{2}} k$ 
and $k \in \{4,40,100,400\}$. 
For fixed $p \in \{1,2,3\}$, we show in Fig.~\ref{fig:square} the 
performance of the $h$-version FEM for $p \in \{1,2,3\}$ on 
quasi-uniform meshes by displaying the relative error in the $H^1$-seminorm 
versus the number of degrees of freedom per wavelength. We observe 
that higher order methods are less prone to pollution. We note that 
the meshes are quasi-uniform, i.e., no geometric mesh refinement 
near the vertices is performed in contrast to the requirements of 
Theorem~\ref{thm:discrete-stability-polygon}. 
}
\eremk
\end{Example}
\begin{Example}{\rm 
\label{ex:Lshaped-smooth}
On the $L$-shaped domain $\Omega = (-1,1)^2\setminus (0,1)\times (-1,0)$ 
with $\Gamma$ being the union of the two edges meeting at $(0,0)$, we consider  
\begin{equation}
\label{eq:model-numerics}
-\Delta u - k^2 u =  0 \quad \mbox{ in $\Omega$}, 
\quad \partial_n u = 0 \quad \mbox{ on $\Gamma$}, 
\quad \partial_n u - \bi k u = g \quad \mbox{ on $\partial\Omega \setminus\Gamma$}, 
\end{equation}
where the Robin data $g$ are such that the exact solution is 
$
u(x,y) = e^{\bi (k_1 x + k_2 y)} $ 
with $k_1 = -k_2 = \frac{1}{\sqrt{2}} k, 
$
and 
$k \in \{10,100,1000\}$. 
We consider 
two kinds of meshes, namely, quasi-uniform meshes ${\mathcal T}_h$ with mesh size 
$h$ such that $kh \approx 4$  
and meshes ${\mathcal T}^{geo}$ that are geometrically refined near the origin. 
The meshes ${\mathcal T}^{geo}$ are derived from the quasi-uniform mesh 
${\mathcal T}_h$ by introducing a geometric grading 
on the elements abutting the origin; the grading factor is $\sigma = 0.125$ and 
the number of refinement levels is $L = 10$. 
Fig.~\ref{fig:numerics} shows 
the relative errors in the $H^1$-seminorm for the $p$-version of the FEM 
where for fixed mesh the approximation order $p$ ranges from $1$ to $10$. 
It is particularly noteworthy that the refinement 
near the origin has hardly any effect on the convergence behavior of the 
FEM; this is quite in contrast to the stability result 
Theorem~\ref{thm:discrete-stability-polygon}, which  requires 
geometric refinement near all vertices of $\Omega$. 
\eremk
}
\end{Example}
\begin{Example}{\rm 
\label{ex:Lshaped-singular}
The geometry and the boundary conditions are as described in
Example~\ref{ex:Lshaped-smooth}. The data $g$ are selected such that the
exact solution is 
$u =  J_{2/3} (kr) \cos \frac{2}{3}\varphi$, 
where $(r,\varphi)$ denote polar coordinates and $J_{\alpha}$ is a first kind 
Bessel function. $k \in \{1,10,20,100,200\}$. Our calculations are 
$p$-FEMs with $p \in \{1,\ldots,10\}$ on
the quasiuniform mesh ${\mathcal T}_h$ 
described in Example~\ref{ex:Lshaped-singular}. 
The results are displayed in the bottom right part of Fig.~\ref{fig:square}. 
The numerics illustrate that the singularity 
at the origin is rather weak: we observe that the asymptotic 
algebraic convergence behavior is 
$|u - u_N|_{H^1(\Omega)} \approx C_{k} p^{-4/3} |u|_{H^1(\Omega)}$, where 
the constant $C_k$ depends favorably on $k$. 
}\eremk
\end{Example}
\subsection{Stability of Partition of Unity Method/Generalized FEM}
\label{sec:pum}
The abstract stability result of Lemma~\ref{lemma:abstract-quasi-optimality} 
only assumes certain approximation properties of the spaces $V_N$. Particularly 
in an ``$h$-version'' setting, even non-polynomial, operator-adapted spaces may have 
sufficient approximation properties to ensure the important 
condition (\ref{eq:cond-for-stability}) for stability. We illustrate 
this effect for the PUM/gFEM, \cite{melenk95,babuska-melenk96b} with 
local approximation spaces consisting of systems of plane waves or 
generalized harmonic polynomials (see 
Section~\ref{sec:approximation-properties-Trefftz-fcts} below) and 
the classical FEM shape functions as the partition of unity. The key observation
is that for $h$ sufficiently small, the resulting space has approximation 
properties similar to the classical (low order) FEM space: 
\begin{Lemma}
\label{lemma:PUM-approximation}
Let ${\mathcal T}$ be a shape-regular triangulation of the 
polygon $\Omega\subset\BbbR^2$. Let $h$ be its mesh size; 
let $(x_i)_{i=1}^M$ be the nodes 
of the triangulation and $(\varphi_i)_{i=1}^M$ be the piecewise
linear hat functions associated with the nodes $(x_i)_{i=1}^M$. 
Assume $kh \leq C_1$ for some $C_1 > 0$. 
Let $V^{master}$ be either the space $V^p_{GHP}$ with $p \ge 0$ (see 
(\ref{eq:generalized-harmonic-polynomials}) below) 
or the space $W^p_{PW}$ with $p \ge 2$ (see 
(\ref{eq:plane-waves}) below). Define, for each $i=1,\ldots,M$,  the local approximation
$V_i$ by $V_i:= \operatorname*{span} \{ b(x - x_i)\colon b \in V^{master}\}$. 
Then the space $V_N:=\sum_{i=1}^M \varphi_i V_i$ has the following approximation
property: There exists $C > 0$ depending only on the shape regularity of ${\mathcal T}$, 
the constant $C_1$, and $V^{master}$ such that 
$$
\inf_{v \in V_N} 
\|u - v\|_{L^2(\Omega)}  + 
h \|u - v\|_{H^1(\Omega)} \leq C \left[ h^2 \|u\|_{H^2(\Omega)} 
+ (kh)^2 \|u\|_{L^2(\Omega)}\right]
\quad \forall u \in H^2(\Omega). 
$$
\end{Lemma}
\begin{proof}
\iftechreport 
We first show that each local space $V_i$ has an element $\psi_i \in V_i$ with 
\begin{equation}
\label{eq:bounds-for-psi_i}
h \|\nabla \psi_i\|_{L^\infty(\widetilde\omega_i)} + 
\|1 - \psi_i\|_{L^\infty(\widetilde \omega_i)} \leq 
C_{app} (kh)^2 
\end{equation}
for some $C_{app} > 0$
independent of $i$ and $h$; here, $\widetilde \omega_i = \operatorname*{supp} \varphi_i$ has 
diameter $O(h)$. It suffices to show (\ref{eq:bounds-for-psi_i}) for the 
set $V^{master}$. For the space of generalized harmonic polynomials, 
this follows from $J_0(kr) = 1 + O((kr)^2)$, and for the space of plane waves, 
Taylor expansion shows that for $p = 2m$ ($m \in \BbbN$) the function
$
\frac{1}{2} \left[ e^{\bi k x} + e^{-\bi k x}\right] = 1 + O((kx)^2)
$
has the desired property whereas for odd $p  =2m+1$ ($m \in \BbbN$) the observation
$$
e^{\bi k \omega_0 \cdot x} - 
\frac{1}{2c} \left[ e^{\bi k \omega_m \cdot x } + e^{\bi k \omega_{m+1} \cdot x}\right] = 
\left(1 - \frac{1}{c}\right) + O((kh)^2), 
\qquad c = \cos \frac{2 \pi m}{2m+1},
$$
can be utilized to construct $\psi_i$. 
We recall (see, e.g., \cite[Thm.~{4.8.7}]{brenner-scott94}) that for each 
$u \in H^2(\Omega)$ there is a function 
$w = \sum_{i} w(x_i) \varphi_i  \in S^{1}({\mathcal T}_h)$ with 
\begin{subequations}
\label{eq:approximation-brenner-scott}
\begin{alignat}{3}
\|w\|_{L^2(\Omega)} &\leq C \|u\|_{L^2(\Omega)},
& \quad \|w\|_{H^1(\Omega)} &\leq C \|u\|_{H^1(\Omega)}, \\
\|u - w\|_{L^2(\Omega)} &\leq C h^2 \|u\|_{H^2(\Omega)}, 
& \quad \|u - w\|_{H^1(\Omega)} &\leq C h\|u\|_{H^2(\Omega)} .
\end{alignat}
\end{subequations}
Upon setting $v:= \sum_{i} w_i(x_i) \psi_i \varphi_i \in V_N$, we get  
in view of $\sum_i \varphi_i \equiv 1$ for the error 
$u - v = u - \sum_{i} w_i(x_i) \psi_i \varphi_i = (u - w) + \sum_i \varphi_i w(x_i) (1 - \psi_i)$. 
The estimates (\ref{eq:approximation-brenner-scott})
imply $\|u - w\|_{L^2(\Omega)} + h\|u  -w\|_{H^1(\Omega)} \leq C h^2 \|u\|_{H^2(\Omega)}$. 
For the sum, we have 
\begin{eqnarray*}
\|\sum_i w_i(x_i) (1 - \psi_i) \varphi_i\|_{L^2(\Omega)} &\leq& 
C (kh)^2 \|w\|_{L^2(\Omega)}\leq 
C (kh)^2 \|u\|_{L^2(\Omega)}, \\
h \|\nabla \sum_i w_i(x_i) (1 - \psi_i) \varphi_i\|_{L^2(\Omega)} &\leq& 
C (kh)^2 \|w\|_{L^2(\Omega)}\leq 
C (kh)^2 \|u\|_{L^2(\Omega)}, 
\end{eqnarray*}
which concludes the proof. 
\fi
\qed
\end{proof}
\iftechreport 
\begin{Remark}
\label{rem:pum-approximation}
{\rm 
The approximation result of Lemma~\ref{lemma:PUM-approximation} can generalized
in various directions. For example, interpolation arguments allow one to construct, 
for $v \in H^{1+\theta}(\Omega)$ with $\theta \in (0,1)$ an approximation 
$v_{app} \in V_N$ 
such that 
$\|v - v_{app}\|_{L^2(\Omega)} + h \|v - v_{app}\|_{H^1(\Omega)} \leq 
C_1 \left( h^{1+\theta} + (kh)^2 h^\theta\right)\|v\|_{H^{1+\theta}(\Omega)} + 
C_2 (kh)^2 \|u\|_{L^2(\Omega)}$. 
(We refer the reader to the Appendix for the proof of this results.)
Furthermore, the approximation result of Lemma~\ref{lemma:PUM-approximation} can be localized, which is of interest if 
${\mathcal T}$ is not quasi-uniform. 
}\eremk
\end{Remark}
\fi 
Lemma~\ref{lemma:PUM-approximation} shows that the space $V_N$, which is 
derived from solutions of the homogeneous Helmholtz equation, nevertheless
has some approximation power for arbitrary functions with some 
Sobolev regularity. Hence, the condition (\ref{eq:cond-for-stability}) can be
met for sufficiently small mesh sizes: 
\begin{Corollary}[\protect{\cite[Prop.~{8.2.7}]{melenk95}}]
\label{cor:gFEM-stability}
Assume the hypotheses of Lemma~\ref{lemma:PUM-approximation}; in particular, 
let the space $V_N$ be constructed from systems of plane waves or generalized
harmonic polynomials.  Assume additionally that $\Omega$ is a convex polygon. 
Then there exists $C >0 $ independent of $k$ such that for $k^2 h \leq C$ the 
Galerkin method for (\ref{eq:problem}) with $f = 0$ is quasi-optimal, i.e., 
the solution $u_N \in V_N$  of (\ref{eq:problem-weak-discrete}) exists and satisfies
$$
\|u - u_N\|_{1,k,\Omega} \leq 2 C_B \inf_{v \in V_N} \|u - v\|_{1,k,\Omega}. 
$$
\end{Corollary}
\begin{proof}
In view of Lemma~\ref{lemma:abstract-quasi-optimality}, we have 
to estimate $\eta(V_N)$. To that end, we consider (\ref{eq:problem}) with $f \in L^2(\Omega)$
and $g=0$. In view of the convexity of $\Omega$, we have $C_{sol}(k) = O(1)$ and 
elliptic regularity then yields for the solution $u$ of (\ref{eq:problem}) 
$$
\|u\|_{1,k,\Omega} + k^{-1} \|u\|_{H^2(\Omega)} \leq C \|f\|_{L^2(\Omega)}. 
$$
This allows us to conclude with Lemma~\ref{lemma:PUM-approximation} that 
\begin{eqnarray*}
\inf_{v \in V_N} \|u - v\|_{1,k,\Omega} &\leq & C \left[ (k h^2 + h)\|u\|_{H^2(\Omega)}+ 
(k (kh)^2 + k^2 h) \|u\|_{L^2(\Omega)}\right]\\
& \leq & C ((kh)^2 + kh))\|f\|_{L^2(\Omega)}
 \leq  C kh (1 + kh) \|f\|_{L^2(\Omega)}. 
\end{eqnarray*}
Hence, $k\eta(V_N)$ can be made sufficiently small if $k^2 h$ is sufficiently small. 
We point out that convexity of $\Omega$ is assumed for convenience---under
more stringent conditions on the mesh size $h$, quasioptimality holds for general polygons.  
\qed
\end{proof}
%
%
\section{Approximation with plane, cylindrical, and spherical waves}
\label{sec:approximation-properties-Trefftz-fcts}
Systems of functions that are solutions of a (homogeneous) differential
equation are often called ``Trefftz systems''. Prominent examples
in the context of the Helmholtz equation are, in the two-dimensional
setting, ``generalized harmonic polynomials'' and systems of 
plane waves given by 
\begin{eqnarray}
\label{eq:generalized-harmonic-polynomials}
V^p_{GHP} &:=& \operatorname*{span} \{ J_n(kr) e^{\bi n \varphi}\colon  - p\leq  n \leq n\},\\
\label{eq:plane-waves}
W^p_{PW} &:=& \operatorname*{span} \{ e^{\bi k \omega_n \cdot (x,y)}\colon n =0,\ldots,p-1\}, 
\qquad \omega_n = (\cos \frac{2\pi n}{p}, \sin \frac{2\pi n}{p}); 
\end{eqnarray}
here, $J_n$ is a first kind Bessel function, the functions in $V^p_{GHP}$ are described in 
polar coordinates and the functions of $W^p_{PW}$ in Cartesian coordinates. 
We point out that analogous systems can 
be developed in 3D. These functions are solutions of the homogeneous Helmholtz
equation. For the approximation of a function $u$ that satisfies the 
homogeneous Helmholtz equation on a domain $\Omega\subset \BbbR^2$, one may 
study the ``$p$-version'', i.e., study how well $u$ can be 
approximated from the spaces $V^p_{GHP}$ or $W^p_{PW}$ as 
$p \rightarrow \infty$; 
alternatively, one may study the ``$h$-version'', 
in which, for fixed  $p$, the approximation properties of the spaces 
$V^p_{GHP}$ or $W^p_{PW}$ are expressed in terms of the diameter 
$h = \operatorname*{diam}\Omega$ of a domain under consideration. 
In the way of illustration, we present two types of results: 
\begin{Lemma}[\protect{\cite{melenk95}}]
\label{lemma:exponential-convergence-plane-waves}
Let $\Omega\subset\BbbR^2$ be a simply connected domain and 
$\Omega^\prime\subset\subset\Omega$ be a compact subset. Let $u$ solve 
$-\Delta u - k^2 u =0$ on $\Omega$. Then there exist constants $C$, $b > 0$
(possibly depending on $k$) such that for all $p \ge 2$: 
$$
\inf_{v \in V^p_{GHP}} \|u - v\|_{H^1(\Omega^\prime)} \leq C e^{-b p}, 
\qquad 
\inf_{v \in W^p_{PW}} \|u - v\|_{H^1(\Omega^\prime)} \leq C e^{-b p/\ln p}. 
$$
\end{Lemma}
\begin{proof}
See, e.g., \cite{melenk95} or \cite[Thm.~{5.3}]{melenk05a}. 
\qed
\end{proof} 
\begin{Remark}
{\rm 
Analogs of Lemma~\ref{lemma:exponential-convergence-plane-waves} hold 
if $u$ has only some finite Sobolev regularity. Then, the convergence rates 
are algebraic, 
\cite{melenk95}, \cite[Thm.~{5.4}]{melenk05a}, \cite{hiptmair-moiola-perugia09a}. 
\eremk
}
\end{Remark}
The approximation properties of the spaces $V^p_{GHP}$ and $W^p_{PW}$ can be 
also be studied in an $h$-version setting: 
\begin{Proposition}[\protect{\cite[Thm.~{3.2.2}]{hiptmair-moiola-perugia09a}}]
\label{prop:h-convergence-plane-waves}
Let $\Omega\subset\BbbR^2$ be a domain with diameter $h$ and inscribed 
circle of radius $\rho h$. Let $p = 2 \mu +1$. Assume $k h \leq C_1$. 
Then there exist $C_p > 0$ (depending only on $C_1$, $\rho>0$, $m$, and $p$)
and $v \in W^{2\mu +1}_{PW}$ such that 
$$
\|u - v\|_{j,k,\Omega,\Sigma} 
\leq C_p h^{\mu - j+1} \|u\|_{\mu+1,k,\Omega,\Sigma}, 
\qquad 0 \leq j \leq \mu +1 ,
$$
where $\|v\|_{j,k,\Omega,\Sigma}^2 = \sum_{m=0}^j k^{2(j-m)} |v|^2_{H^m(\Omega)}$. 
\end{Proposition}
Several comments concerning 
Proposition~\ref{prop:h-convergence-plane-waves} are in order: 
\begin{enumerate}
\item 
The constant $C_p$ in Proposition~\ref{prop:h-convergence-plane-waves} 
depends favorably on $p$, and its dependence on $p$ 
can be found in \cite[Thm.~{3.2.3}]{hiptmair-moiola-perugia09a}.  
\item
Proposition~\ref{prop:h-convergence-plane-waves} is formulated for 
the space $W^p_{PW}$ of plane waves---analogous results are valid 
for generalized harmonic polynomials, see 
\cite[Thm.~{2.2.1}]{hiptmair-moiola-perugia09a} for both the 
$h$ and $hp$-version. 
\item
Proposition~\ref{prop:h-convergence-plane-waves} is formulated for the 
two-dimensional case. Similar results are available in 3D, 
\cite{hiptmair-moiola-perugia09a}. 
\item 
The approximation properties of plane waves in terms of the element
size have previously been studied in slightly different norms in 
\cite{cessenat-despres98}. 
\end{enumerate}
\iftechreport 
\begin{Remark}{\rm 
Plane waves and generalized harmonic polynomials represent by no means the
only operator adapted system used in practice. Especially for polygonal
domains, the functions $J_{n \alpha} (kr) \sin (\alpha n \varphi)$, $n \in \BbbN$,
or $J_{n \alpha} (kr) \sin (\alpha n \varphi)$, $n \in \BbbN_0$, for suitable 
$\alpha$ can combine good approximation properties with the option to realize 
homogeneous boundary conditions, \cite{betcke-barnett10}. Further possibilities
include linear combinations of fundamental solutions or, more generally, discretizations 
of layer potentials. We refer to \cite{betcke-barnett10} for a concrete example. 
}\eremk
\end{Remark}
\fi 
\section{Stability of Least Squares and DG methods}
Discrete stability in Section~\ref{sec:discrete-stability-H1-conforming} 
is obtained from stability of the continuous problem by a perturbation argument. 
This approach does not seem to work very well if one aims at using approximation spaces 
that have special features linked to the differential equation under consideration. 
The reason can be seen from the proof of Lemma~\ref{lemma:abstract-quasi-optimality}: 
The adjoint approximation property $\eta(V_N)$ (which needs to be small) measures 
how well certain solutions to the {\em in}homogeneous equation can be 
approximated from the test space. If the ansatz space is based on solutions 
of the homogeneous
equation, then its capabilities to approximate solutions of the inhomogeneous equation
are  clearly limited.  In an $h$-version, the situation is not as severe
as we have just seen in Section~\ref{sec:pum} for the PUM/gFEM. In a pure 
$p$-version setting, however, the techniques of Section~\ref{sec:pum} do 
not seem applicable. 

An option is to leave the setting of Galerkin methods and to work with
formulations with built-in stability properties.  
Such approaches can often be understood as minimizing some residual norm, which 
then provides automatically stability and quasi-optimality (in this residual norm). 
We will illustrate this procedure here by two examples, namely, 
Least Squares methods and DG-methods. Our presentation will highlight an issue 
stemming from this approach, namely, the fact that error estimates in this 
residual norm do not easily lead to error estimates in more classical norms
such as the $L^2(\Omega)$-norm. 
\subsection{Some notation for spaces of piecewise smooth functions}
Let ${\mathcal T}$ be a regular, shape-regular triangulation
of the polygon $\Omega \subset\BbbR^2$. We decompose the set of 
edges ${\mathcal E}$ 
as ${\mathcal E} = \EI \dot \cup \EB$, where 
$\EI$ is the set of edges in $\Omega$ and $\EB$
consists of the edges on $\partial\Omega$. For functions $u:\Omega\rightarrow\BbbR$ 
and $\bsigma:\Omega\rightarrow\BbbR^2$ that are smooth on the elements $K \in {\mathcal T}$, 
we define the jumps and averages as it is customary in DG-settings: 
\begin{itemize}
\item
For $e \in \EI$, let $K_e^+$ and $K_e^-$ be the two elements sharing $e$
and denote by ${\mathbf n}^+$ and ${\mathbf n}^-$ the normal vectors on $e$ 
pointing out of $K_e^+$ and $K_e^-$. 
Correspondingly, we let $u^+$, $u^-$ and $\bsigma^+$ and $\bsigma^-$ 
be traces on $e$ of $u$ and $\bsigma$ from $K^+_e$ and $K^-_e$. We set: 
\begin{eqnarray*}
&& \average{u}|_e:= \frac{1}{2}\left(u^+ + u^-\right),
\qquad
\average{\bsigma}|_e:= \frac{1}{2}\left(\bsigma^+ + \bsigma^-\right), \\
&& \jump{u}|_e := u^+ {\mathbf n}^+ + u^- {\mathbf n}^-,
\qquad
\jump{\bsigma}|_e := \bsigma^+ \cdot {\mathbf n}^+ + \bsigma^- \cdot {\mathbf n}^-.
\end{eqnarray*}
\item
For  {boundary edges} $e \in \EB$ we define
$$
\average{\bsigma}|_e:= \bsigma |_e
\qquad
\jump{u}|_e := u|_e {\mathbf n}
$$
\end{itemize}
With this notation, one can conveniently rearrange certain sums over edges: 
\begin{Lemma}[``DG magic formula'']
\label{lemma:DG-magic-formula}
Let $v:\Omega\rightarrow \BbbR$ and $\bsigma:\Omega\rightarrow \BbbR^2$ be piecewise smooth
on the triangulation ${\mathcal T}$. Then: 
$$
\sum_{K \in {\mathcal T}} \int_{\partial K} v \bsigma \cdot {\mathbf n} =
\int_{\EI} \jump{v} \cdot \average{\bsigma} +
\int_{\EI} \average{v} \cdot \jump{\bsigma} +
\int_{\EB} \jump{v} \cdot \average{\bsigma},
$$
where $\int_{\EI}$ and $\int_{\EB}$ are shorthand notations for the sums of 
integrals over all edges in $\EI$ and $\EB$.
\end{Lemma}
Finally, 
for piecewise smooth functions, $\nabla_h$ denotes the piecewise defined 
gradient. 

\subsection{Stability of least squares methods}
\label{sec:least-squares}
Although Least Squares methods could be based on any space of approximation
spaces, we will concentrate here on the approximation by piecewise solutions 
of the homogeneous Helmholtz equation. With varying focus, this is the
setting of \cite{stojek98, monk-wang99, li07, betcke-barnett10, desmet07} 
and references therein.
We illustrate the procedure for the model problem (\ref{eq:problem})
with $f = 0$. 
The approximation space has the form 
\begin{equation}
\label{eq:VN-least-squares}
V_N = \{u \in L^2(\Omega)\colon u|_K \in V_{N,K} \quad \forall K \in {\mathcal T}\},
\end{equation}
where the spaces $V_{N,K}$ are spaces of solutions of the homogeneous 
Helmholtz equation, e.g., systems of plane waves. For each edge 
$e \in {\mathcal E}$, we select weights $w_{1,e}$, $w_{2,e}>0$ and define
the functional $J:V_N \rightarrow \BbbR$ by 
\begin{equation*} 
J(v):= 
\sum_{e \in \EI} w_{1,e}^2 \|[v]\|^2_{L^2(e)} + w_{2,e}^2 \|[\partial_n u]\|^2_{L^2(e)}
+ \sum_{e \in \EB} w_{2,e}^2 \|g - (\partial_n v + \bi k v)\|^2_{L^2(e)}; 
\end{equation*} 
here $[v]|_e := \jump{v}|_e$ and $[\partial_n v]|_e := \jump{\nabla_h v}|_e$ represent the jumps of 
$v$ and $\partial_n v$ across 
the edge $e$. If the exact solution $u$ of (\ref{eq:problem}) is sufficiently regular, then 
it is a minimizer of $J$ with $J(u) = 0$. In a Least Squares method, $J$ is minimizer over 
a finite dimensional space $V_N$ of the form (\ref{eq:VN-least-squares}). Its variational form reads: 
\begin{eqnarray}
\label{eq:minimization}
\mbox{ find $u_N \in V_N$ s.t.} 
\langle u_N,v\rangle_{J,N} &=& \sum_{e \in \EB} (g , \partial_n v + \bi k v)_{L^2(e)} 
\qquad \forall v \in V_N,
\end{eqnarray}
where 
\begin{eqnarray*}
\lefteqn{
\langle u,v\rangle_{J,N} :=} \\
&&
\sum_{e \in \EI} w_{1,e}^2 ([u],[v])_{L^2(e)} \!+\! w_{2,e}^2 ([\partial_n u],[\partial_n v])_{L^2(e)}
+ 
\sum_{e \in \EB} \!\!w_{2,e}^2 (\partial_n u + \bi k u,\partial_n v + \bi k v)_{L^2(e)}. 
\end{eqnarray*}
The positive semidefinite sesquilinear form $\langle\cdot,\cdot\rangle_{J,N}$ induces in fact a 
norm on $V_N$: To see the definiteness of $\langle\cdot,\cdot\rangle_{J,N}$, we note that 
$v \in V_N$ and  $J(v) = 0$ implies that $v$ is in $C^1(\Omega)$ and elementwise a 
solution of the homogeneous Helmholtz equation. Thus, it is a classical solution of the 
Helmholtz equation on $\Omega$ and satisfies $\partial_n v + \bi k v = 0$ on $\partial\Omega$. 
The uniqueness assertion for (\ref{eq:problem}) with $f = 0$ and $g = 0$ worked out in 
Example~\ref{example:uniqueness} then implies $v = 0$. 
Therefore, the minimization problem (\ref{eq:minimization}) is well-defined. If the 
solution $u$ of (\ref{eq:problem}) satisfies $u \in H^{3/2+\varepsilon}(\Omega)$ for some 
$\varepsilon > 0$, then $J(u) = 0$, and we get quasi-optimality 
of the Least Squares method in the norm $\|\cdot\|_{J,N} = J(\cdot)^{1/2}$: 
\begin{equation}
\label{eq:least-squares-quasioptimality}
\|u - u_N\|_{J,N}^2 = J(u - u_N) = J(u_N)  = \min_{v \in V_N} J(v) = \|u - v\|^2_{J,N}. 
\end{equation}
We mention here that estimates for this minimum can be obtained from (local) estimates in 
classical Sobolev norms as given in Section~\ref{sec:approximation-properties-Trefftz-fcts}
using appropriate trace estimates. 
Turning estimates for $\|u - u_N\|_{J,N} = J(u_N)^{1/2}$ into estimates in terms of more
familiar norms such as $\|u - u_N\|_{L^2(\Omega)}$ is not straight forward. It may be expected
that the norm of the solution operator of the continuous problem appears again; 
the next result, which is closely related to 
\cite{buffa-monk07,hiptmair-moiola-perugia09a,hiptmair-moiola-perugia11,moiola09}, 
illustrates the kind of result one can obtain, in particular in a $p$-version setting: 
\begin{Lemma}[\protect{\cite[Thm.~{3.1}]{monk-wang99}}]
\label{lemma:monk-wang}
Let $\Omega\subset\BbbR^2$ be a polygon. Let $w_{1,e} = k$ and $w_{2,e} = 1$ 
for all edges and $g \in L^2(\partial\Omega)$. Let $u_N \in V_N$ be the minimizer
of $J$, where $V_N$, given by (\ref{eq:VN-least-squares}), consists of elementwise 
solutions of the homogeneous Helmholtz equation. 
\begin{enumerate}[(i)]
\item 
\label{item:lemma:monk-wang-i}
If $\Omega$ is convex, then 
$ \|u - u_N\|^2_{L^2(\Omega)} \leq C k^{-1} \left( (kh)^{-1} + (kh)^{1}\right) 
J(u_N)$.  
\item 
\label{item:lemma:monk-wang-ii}
If $\Omega$ is not convex, then 
\begin{eqnarray*}
\lefteqn{
\|u - u_N\|^2_{L^2(\Omega)} \leq }\\
&&C k^{-1} \left[ (kh)^{-1} + (kh)^{1} \left\{ 1 + \min\{1,kh\}^{-2\beta_{max}} \right\}
\right] (C_{sol}(k)+1)^2 J(u_N),
\end{eqnarray*}
where $C_{sol}(k)$ is defined in (\ref{eq:def-Csol}) and satisfies $C_{sol}(k) = O(k^{5/2})$
by Theorem~\ref{thm:a-priori}. The parameter $\beta_{max}\ge 0$ can be selected arbitrarily
to satisfy the condition $\beta_{max} > 1 - \min_i \frac{\pi}{\omega_i}$, where the $\omega_i$ are the 
interior angles of the polygon.  
\end{enumerate}
\end{Lemma}
\begin{proof}
The result (\ref{item:lemma:monk-wang-i}) is essentially the statement 
of \cite[Thm.~{3.1}]{monk-wang99} in a refined form as given in 
\cite[Lemma~{3.7}]{hiptmair-moiola-perugia11}. 
\iftechreport
The statement 
(\ref{item:lemma:monk-wang-ii}) is a slightly modification of 
(\ref{item:lemma:monk-wang-i}), and we restrict our presentation to that case. 
The key idea is to obtain $L^2(\Omega)$-bounds by a duality argument and use the fact that 
$u - u_N$ solves the homogeneous Helmholtz equation elementwise. More precisely, given
$\varphi \in L^2(\Omega)$ we let $v \in H^1(\Omega)$ solve the adjoint problem 
$$
-\Delta v - k^2 v = \varphi \quad \mbox{ in $\Omega$}, 
\qquad \partial_n v - \bi k v = 0 \quad \mbox{ on $\partial\Omega$}.
$$
By Corollary~\ref{cor:H2-regularity}, the function $v$ is in a weighted $H^2$-space with 
\begin{equation}
\label{eq:monk-wang-a-priori-for-v}
\|v\|_{1,k,\Omega} + k^{-1} \|\Phi_{0,\overrightarrow\beta,k}\nabla^2 v\|_{L^2(\Omega)} 
\leq 
C (C_{sol}(k) +1)\|\varphi\|_{L^2(\Omega)}. 
\end{equation}
Inspection of the arguments underlying the proof of Corollary~\ref{cor:H2-regularity}
shows that the exponents $\beta_j \in [0,1)$ stem from the regularity
theory for the Laplacian with Neumann boundary conditions. Hence, in fact 
$\beta_j \in [0,1/2)$  so that $\nabla_h v$ has an $L^2$-trace on all edges 
of the triangulation (cf.\ Lemma~\ref{lemma:trace-inequality-H11beta}).
For each $K \in {\mathcal T}$ we then have 
\begin{eqnarray}
\label{eq:trace-1}
\|w\|^2_{L^2(\partial K)} &\leq& C \left[ h^{-1} \|w\|^2_{L^2(K)} + h |w|^2_{H^1(K)} \right]
\qquad \forall w \in H^1(K), \\
\label{eq:trace-2}
\|\nabla w\|^2_{L^2(\partial K)} &\leq& C \left[ h^{-1} |w|^2_{H^1(K)} + h |w|^2_{H^2(K)} \right]
\qquad \forall w \in H^2(K),\\
\label{eq:trace-3}
\|\nabla w\|^2_{L^2(\partial K)} &\leq& C \left[ h^{-1} |w|^2_{H^1(K)} +  
h^{1-2\beta} \|r^\beta \nabla^2 w\|^2_{L^2(K)} \right]
\quad \forall w \in H^{2,2}_\beta(K), 
\end{eqnarray}
where, in the last estimate we assume that the origin is at one corner of $K$
and $r$ denotes the distance from that corner. These estimates are obtained with the 
aid of scaling arguments, the multiplicative trace inequality
(see \cite[Prop.~{1.6.3}]{brenner-scott94}), and, in the case of (\ref{eq:trace-3}) 
additionally Lemma~\ref{lemma:trace-inequality-H11beta}. From 
$a b = \max\{a,b\}\min\{a,b\}$
(for $a$, $b \ge 0$), we get $r^\beta = k^{-\beta} (rk)^\beta = 
k^{-\beta} \min\{1,rk\}^\beta \max\{1,rk\}^\beta$. Hence, if ${\mathcal T}_{corner}$ 
denotes the set of elements that abut on the corners of $\Omega$, we get 
$$
\sum_{K \in {\mathcal T}_{corner}} \|\nabla v\|^2_{L^2(\partial K)} \leq C \left[ 
h^{-1} |v|^2_{H^1(B_h)} + 
h \min\{1,hk\}^{-2\beta_{max}} \|\Phi_{0,\overrightarrow\beta,k} \nabla^2 v\|^2_{L^2(B_h)}
\right], 
$$
where $B_h = \cup_{K \in {\mathcal T}_{corner}} K$. Noting that the elements 
$K \in {\mathcal T}\setminus {\mathcal T}_{corner}$ are at least $O(h)$ 
away from the corners allows us to estimate with (\ref{eq:trace-2})
$$
\sum_{K \in {\mathcal T}\setminus{\mathcal T}_{corner}} \!\!\!\!\!\!\!\!\! \|\nabla v\|^2_{L^2(\partial K)} 
\!\leq \! C \!\!\left[ h^{-1} |v|^2_{H^1(\Omega\setminus B_h)} + 
h^{1} \min\{1,kh\}^{-2\beta_{max}} \|\Phi_{0,\overrightarrow\beta,k} \nabla^2 v\|^2_{L^2(\Omega\setminus B_h)}
\right]. 
$$
Hence, we have the two bounds 
\begin{eqnarray*}
\sum_{K \in {\mathcal T}} \|\nabla v\|^2_{L^2(\partial K)} &\leq & C h^{-1} |v|^2_{H^1(\Omega)} + 
C h^1 \min\{1,kh\}^{-2\beta_{max}} \|\Phi_{0,\overrightarrow\beta,k} \nabla^2 v\|^2_{L^2(\Omega)}, \\
\sum_{K \in {\mathcal T}} \|v\|^2_{L^2(\partial K)} &\leq& C h^{-1} \|v\|^2_{L^2(\Omega)} + h|v|^2_{H^1(\Omega)}. 
\end{eqnarray*}
Therefore, recalling $w_{1,e} = k$ and  $w_{2,e} =1$, we obtain from these estimates
and the {\sl a priori} estimate (\ref{eq:monk-wang-a-priori-for-v}) the bound
\begin{eqnarray}
\nonumber 
\lefteqn{
\sum_{e \in {\mathcal E}} w_{2,e}^{-2} \|v\|^2_{L^2(e)} + w_{1,e}^{-2} \|\nabla v\|^2_{L^2(e)} }\\
\label{eq:lemma:monk-wang-200}
&\leq& C k^{-1} \left[ \frac{1}{kh} + kh \left\{ 1 + \min\{1,kh\}^{-2\beta_{max}}\right\}\right] 
(C_{sol}(k)+1)^2 \|\varphi\|^2_{L^2(\Omega)}. 
\end{eqnarray}
The estimate (\ref{eq:lemma:monk-wang-200}) can be used
to bound $|(u - u_N,\varphi)_{L^2(\Omega)}|$: Writing the integral as a sum
over elements, integrating by parts twice and using that $u - u_N$ solves the homogeneous Helmholtz
equation elementwise yields 
\begin{eqnarray*}
(u - u_N,\varphi)_{L^2(\Omega)} &=& \sum_{K \in {\mathcal T}} (u - u_N,-\Delta v - k^2 v)_{L^2(K)}\\
&=& 
\sum_{K \in {\mathcal T}} (\partial_n (u - u_N),v)_{L^2(\partial K)} - 
\sum_{K \in {\mathcal T}} (u - u_N,\partial_n v)_{L^2(\partial K)}
=: \Sigma_1 - \Sigma_2.  
\end{eqnarray*}
The ``DG magic formulas'' of Lemma~\ref{lemma:DG-magic-formula}  produce 
\begin{align*}
\!\!\Sigma_1 \!\!\!&= \!\!\!\!
\sum_{e \in \EI}\!\!\! (\average{\nabla (u - u_N)},\jump{v})_{L^2(e)} \! + \!
(\jump{\nabla (u - u_N)},\average{v})_{L^2(e)} 
\!\!+\!\!\!\! 
\sum_{e \in \EB}\!\! (\average{\nabla (u - u_N)},\jump{v})_{L^2(e)} \\
 \Sigma_2\!\!  &= 
\!\!\! \sum_{e \in \EI}\!\! (\jump{u - u_N},\average{\nabla v})_{L^2(e)} + 
(\average{u - u_N},\jump{\nabla v})_{L^2(e)} + 
\!\!\! \sum_{e \in \EB}\!\! (\jump{u - u_N},\average{\nabla v})_{L^2(e)}.  
\end{align*}
For interior edges $e \in \EI$ we have 
$\jump{v} = 0$ and $\jump{\nabla v} = 0$ as well as $\jump{u} = 0$ and $\jump{\nabla u} = 0$; 
on boundary edges $e \in \EB$ we have with the boundary conditions
satisfied by $u$ and $v$ (i.e., $\partial_n u + \bi k u = g$ and $\partial_n v - \bi k v = 0$)
$$
(\average{\nabla( u - u_N)},\jump{v})_{L^2(e)} 
- (\jump{u - u_N},\average{\nabla v})_{L^2(e)}  
= -  
((\partial_n + \bi k) u_N - g,v)_{L^2(e)}. 
$$
These observations lead to 
\begin{eqnarray*}
\lefteqn{
\left|-(u - u_N,\varphi)_{L^2(\Omega)} \right| = }\\
&&
\left|
\sum_{e \in \EI} (\jump{u_N},\average{\nabla v})_{L^2(e)} 
+ \sum_{e \in \EI} (\jump{\nabla u_N},\average{v})_{L^2(e)} + 
\sum_{e \in \EB} 
((\partial_n + \bi k) u_N - g,v)_{L^2(e)}
\right|\\
&\leq& C \sqrt{ J(u_N) }
\sqrt{\sum_{e \in {\mathcal E}} w_{2,e}^{-2}\|\average{v}\|_{L^2(e)}^2 
+ w_{1,e}^{-2} \|\average{\nabla v}\|^2_{L^2(e)}}, 
\end{eqnarray*}
where, in the last step, we employed the Cauchy-Schwarz inequality for 
sums. From (\ref{eq:lemma:monk-wang-200}) we therefore get 
\begin{eqnarray*}
\lefteqn{
\frac{|(u - u_N,\varphi)_{L^2(\Omega)}|}{\|\varphi\|_{L^2(\Omega)}} \leq } 
\\ &&
C \sqrt{J(u_N)} k^{-1/2} 
\left[ (kh)^{-1/2} + (kh)^{1/2} \min\{1,kh\}^{-\beta_{max}} \right] (C_{sol}(k)+1). 
\end{eqnarray*}
Since $\varphi \in L^2(\Omega)$ is arbitrary, we get the result. 
\qed
\fi 
\end{proof}
\begin{Remark}
\label{rem:monk-wang}
{\rm 
Lemma~\ref{lemma:monk-wang} assumes quasi-uniform meshes and the weights 
$w_{1,e}$, $w_{2,e}$ do not take the edge length into account. This limits 
somewhat it applicability in an $h$-version context. However, the 
result is very suitable for a $p$-version setting. We point out that 
in a case where the $p$-version features only algebraic rates 
of convergence, 
one would have to give the parameters $w_{1,e}$, $w_{2,e}$ a $p$-dependent 
relative weight as opposed to the situation studied in Lemma~\ref{lemma:monk-wang}.
}\eremk
\end{Remark}
\subsection{Stability of plane wave DG and UWVF}
\label{sec:DG}
The framework of Discontinuous Galerkin (DG) methods permits another 
way of deriving numerical schemes that are inherently stable. In a classical, 
piecewise polynomial setting, this is pursued in \cite{feng-wu09,feng-wu11,feng-xing10}; 
related work is in \cite{monk-schoeberl-sinwel10}. Here, we concentrate on 
a setting where the ansatz functions satisfy the homogeneous Helmholtz equation.  
In particular, we study the plane wave DG method, 
\cite{gittelson-hiptmair-perugia09,hiptmair-moiola-perugia11,moiola09}, 
and the closely related Ultra Weak Variational Formulation (UWVF), 
\cite{cessenat-despres98,cessenat-despres03, huttunen-monk07,luostari-huttunen-monk09,buffa-monk07}. 
We point out that the UWVF can be derived in different way. Here, we follow
\cite{buffa-monk07,gittelson-hiptmair-perugia09} in viewing it as a special
DG method. 

Our model problem (\ref{eq:problem}) can be 
reformulated as a first order system by 
introducing the flux $\bsigma := (1/\bi k)\nabla u$: 
\iftechreport %
\begin{subequations}
\begin{eqnarray}
\label{eq:flux-formulation-1}
\bi k \bsigma &=& \nabla u \quad  \mbox{ in $\Omega$},\\
\label{eq:flux-formulation-2}
\bi k u - \nabla \cdot \bsigma &=& 0 \quad \mbox{ in $\Omega$},\\
\label{eq:flux-formulation-3}
\bi k \bsigma \cdot {\mathbf n} + \bi k u &=& g \quad \mbox{ on $\partial\Omega$}. 
\end{eqnarray}
\end{subequations}
\fi 
\iftechreport 
For a mesh ${\mathcal T}$, the weak elementwise formulation of 
(\ref{eq:flux-formulation-1}), (\ref{eq:flux-formulation-2}) is  for every 
$K \in {\mathcal T}$:  
\fi 
\begin{eqnarray*}
\int_K \bi k \bsigma \cdot \overline{\btau} + \int_K u \nabla \cdot \overline{\btau} - 
\int_{\partial K} u \overline{\btau}\cdot {\mathbf n} &= & 0 
\qquad \quad \forall \btau \in H(\operatorname*{div},K),\\
\int_{K} \bi k u \overline{v} + \int_K \bsigma \cdot \nabla \overline{v} - 
\int_{\partial K} \bsigma \cdot {\mathbf n} \overline{v} &=& 
0\qquad \forall v \in H^1(K), 
\end{eqnarray*}
where $H(\operatorname*{div},K) = \{u \in L^2(K)\colon \operatorname*{div} u \in L^2(K)\}$
and ${\mathbf n}$ is the outward pointing normal vector.  
Replacing the spaces 
$H^1(K)$ and $H(\operatorname*{div},K)$ by finite-dimensional 
subsets $V_{N,K} \subset H^1(K)$ and 
$\bSigma_{N,K} \subset H(\operatorname*{div},K)$ and, additionally, 
imposing a coupling between neighboring elements by replacing the 
multivalued traces $u$ and $\bsigma$ on the element edges by single-valued
numerical fluxes $\widehat u_N$, $\widehat \bsigma_N$ to be specified below, 
leads to the problem: Find $(u_N,\bsigma_N) \in V_{N,K} \times \bSigma_{N,K}$ such that 
\begin{eqnarray*}
\label{eq:dg-flux-formulation-1}
\int_K \bi k \bsigma_N \cdot \overline{\btau} + \int_K u_N \nabla \cdot \overline{\btau} - 
\int_{\partial K} \widehat u_N \overline{\btau}\cdot {\mathbf n} &= & 0 
\quad \forall \btau \in \bSigma_{N,K},\\
\label{eq:dg-flux-formulation-2}
\int_{K} \bi k u_N \overline{v} + \int_K \bsigma_N \cdot \nabla \overline{v} - 
\int_{\partial K} \widehat \bsigma_N \cdot {\mathbf n} \overline{v} &=& 
0 \quad \forall v \in V_{N,K}. 
\qquad 
\end{eqnarray*}
The variable $\bsigma_N$ can be eliminated by making the assumption
that $\nabla V_{N,K}\subset\bSigma_{N,K}$ for all $K\in {\mathcal T}$
and then selecting the test function $\btau = \nabla v$ on each 
element. This yields after an integration by parts:  
\begin{equation}
\label{eq:DG-formulation-elementwise}
\int_K \nabla u_N \nabla \overline{v} - k^2 u_N \overline{v} - 
\int_{\partial K} (u_N - \widehat u_N) \partial_n \overline{v} - 
\bi k \widehat\bsigma_N\cdot{\bf n} \overline{v} = 0 \qquad \forall K \in {\mathcal T}.
\end{equation}
Since $V_N = \{u \in L^2(\Omega)\colon u|_K \in V_{N,K} \forall K \in {\mathcal T}\}$ 
consists of discontinuous functions without any interelement continuity
imposed across the element edges, 
(\ref{eq:DG-formulation-elementwise}) is equivalent to the sum over the 
elements: 
Find $u_N \in V_N$ such that  for all $v \in V_N$
\begin{equation}
\label{eq:plane-wave-DG}
\sum_{K \in {\mathcal T}} 
\int_K \nabla u_N \cdot \nabla \overline{v} - k^2 u_N \overline{v} + 
\int_{\partial K} ( \widehat u_N - u_N) \nabla \overline{v} \cdot {\mathbf n} - 
\int_{\partial K} \bi k\widehat \bsigma_N \cdot {\mathbf n} \overline{v} = 0.
\end{equation}
This formulation is now completed by specifying the fluxes $\widehat u_N$ 
and $\widehat \bsigma_N$, which at the same time takes care of the 
boundary condition 
\iftechreport 
(\ref{eq:flux-formulation-3}): 
\fi 
\iftechreport 
\begin{subequations}
\label{eq:flux-choices}
\begin{itemize}
\item For interior edges $e \in \EI$
\begin{eqnarray}
\widehat \bsigma_N = \frac{1}{\bi k} \average{\nabla_h u_N}-\alpha \jump{u_N}, 
\qquad  
\widehat u_N = \average{u_N}-\beta \frac{1}{\bi k}\jump{\nabla_h u_N} . 
\end{eqnarray}
\item For boundary edges $e \in \EB$
\begin{eqnarray}
\widehat \bsigma_N &=& \frac{1}{\bi k} \nabla_h u_N-\frac{1-\delta}{\bi k} 
\left(\nabla_h u_N + \bi k u_N {\mathbf n} - g {\mathbf n} \right).\\
\widehat u_N &=& u_N - \frac{\delta}{\bi k }
\left(\nabla_h u \cdot{\mathbf n} + \bi k u_N  -g
\right).
\end{eqnarray}
\end{itemize}
\end{subequations}
\fi 
Different choices of the parameters $\alpha$, $\beta$, $\delta$
lead to different methods analyzed in the literature. For example: 
\begin{enumerate}
\item
$\alpha=  \beta = \delta= 1/2$: this is the UWVF as analyzed in 
\cite{cessenat-despres98,cessenat-despres03, huttunen-monk07,luostari-huttunen-monk09,buffa-monk07}
if the spaces $V_{N,K}$ consist of a space $W^p_{PW}$ of plane waves. 
\item
$\alpha = O(p/(kh \log p))$, \quad $\beta = O((kh \log p)/p)$, \quad
$\delta = O((kh \log p)/p)$: this choice is introduced and advocated in 
\cite{hiptmair-moiola-perugia11,moiola09} in conjunction with $V_{N,K} = W^p_{PW}$. 
\end{enumerate}
With these choices of fluxes, the formulation (\ref{eq:plane-wave-DG}) takes the form
\begin{equation}
\label{eq:DG-formulation}
\mbox{ Find $u_N \in V_N$ s.t.} \quad A_N(u_N,v) = l(v) \qquad \forall v \in V_N, 
\end{equation}
where the sesquilinear form $A_N$ and the linear form $l$ are given by 
\begin{align}
\nonumber 
& A_N(u,v) \!=\!\! \int_\Omega \nabla_h u \cdot \nabla_h \overline{v} \!-\! k^2 u \overline{v} 
-\!\! \int_{\EI}\!\! \jump{u}\average{\nabla_h \overline{v}}
-\!\! \int_{\EI}\!\! \average{\nabla_h u} \jump{\overline{v}}
-\!\! \int_{\EB}\!\!\! \delta u \partial_n \overline{v}  
-\!\! \int_{\EB}\!\!\!\! \delta \partial_n u \overline{v} \\
\label{eq:AN-first-version}
&-\frac{1}{\bi k }\int_{\EI} \beta \jump{\nabla_h u}\jump{\nabla_h \overline{v}}
-\frac{1}{\bi k }\int_{\EB} \delta \partial_n u \partial_n \overline{v}
+\bi k \int_{\EI} \alpha \jump{u}\jump{\overline{v}}
+\bi k \int_{\EB} (1-\delta) u \overline{v}\\
\nonumber 
&l(v) = - \frac{1}{\bi k } \int_{\EB} \delta g \partial_n \overline{v} 
+ \int_{\EB} (1-\delta)g \overline{v}. 
\end{align}
So far, the choice of the spaces $V_{N,K}$ is arbitrary. If the approximation 
spaces $V_{N,K}$ (more precisely: the test spaces) consist of piecewise solutions 
of the homogeneous Helmholtz equation, then a further integration by parts 
is possible to eliminate all volume contributions in $A_N$. Indeed, 
Lemma~\ref{lemma:DG-magic-formula} produces 
\begin{eqnarray*}
\sum_{K \in {\mathcal T}} \int_K \nabla u \cdot \nabla \overline{v} - k^2  u \overline{v} 
&=& \sum_{K \in {\mathcal T}} \int_{\partial K} u \nabla \overline{v} {\mathbf n} 
 = \int_{\EI} \jump{u} \average{\nabla \overline{v}}
+ \average{u} \jump{\nabla \overline{v}} 
+ \int_{\EB} \jump{u} \average{\nabla \overline{v}}
\end{eqnarray*}
so that $A_N$ simplifies to 
\begin{eqnarray*}
A_N(u,v) &= & \int_{\EI} \average{u}\jump{\nabla_h \overline{v}} 
   +
\bi\frac{1}{k}\int_{\EI}\beta \jump{\nabla_h u}\jump{\nabla_h \overline{v}}
- \int_{\EI} \average{\nabla_h u}\jump{\overline{v}} 
 + 
\bi k \int_{\EI} \alpha \jump{u}\jump{\overline{v}} 
\\
&& \mbox{} + \int_{\EB} (1-\delta) u \partial_n \overline{v} 
  + 
\bi \frac{1}{k}\int_{\EB} \delta \partial_n u \partial_n \overline{v} 
- \int_{\EB} \delta \partial_n u \overline{v} 
 + 
\bi k  \int_{\EB} (1-\delta)u  \overline{v}. 
\end{eqnarray*}
Next, we make the important observation that 
$\operatorname*{Im} A_N$ induces a norm on the space $V_N$
if $\alpha$, $\beta > 0$ and $\delta \in (0,1)$. Indeed: 
\begin{enumerate}
\item
$\alpha$, $\beta > 0$ and $\delta \in (0,1)$ implies 
$\operatorname*{Im} A_N(v,v) \ge 0 \quad \forall v \in V_N$
by inspection of (\ref{eq:AN-first-version}).
\item 
$\operatorname*{Im} A_N(v,v) = 0$ and the fact that $V_N$ consists of elementwise
solutions of the homogeneous Helmholtz equation 
implies as in the case of $\langle \cdot,\cdot\rangle_{J,N}$
in Section~\ref{sec:least-squares} that $v \in C^1(\Omega)$ solves 
the homogeneous
Helmholtz equation and $\partial_n v  =  v = 0$ on $\partial\Omega$; the uniqueness assertion
of Example~\ref{example:uniqueness} then proves $v \equiv  0$. 
\end{enumerate}
This is at the basis of the convergence analysis. Introducing 
\begin{eqnarray*}
{\|u\|^2_{DG}} &:=& \sqrt{\operatorname*{Im} A_N(u,u)} = 
\frac{1}{k}\|\beta^{1/2}\jump{\nabla_h u}\|^2_{L^2(\EI)}
+ \|\alpha^{1/2} \jump{u}\|^2_{L^2(\EI)} 
\\
&&\phantom{\sqrt{\operatorname*{Im} A_N(u,u)}} \quad \mbox{}
+ \frac{1}{k} \|\delta^{1/2} \partial_n u\|^2_{L^2(\EB)} + k 
\|(1-\delta)^{1/2} u\|^2_{L^2(\EB)},\\
{\|u\|^2_{DG,+} } &:=& {\|u\|^2_{DG} } + k \|\beta^{-1/2} \average{u}\|^2_{L^2(\EI)}
+  k^{-1} \|\alpha^{-1/2} \average{u}\|^2_{L^2(\EI)}+ 
k \|\delta^{-1/2}  u\|^2_{L^2(\EB)},
\end{eqnarray*}
we can formulate coercivity and continuity results: 
\begin{Proposition}[\protect{\cite{buffa-monk07,hiptmair-moiola-perugia11}}]
\label{prop:DG-quasi-optimality}
Let $V_N$ consist of piecewise solutions of the homogeneous Helmholtz equation. 
Then $\|\cdot\|_{DG}$ is a norm on $V_N$ and for some $C>0$ depending solely on
the choice of $\alpha$, $\beta > 0$, and $\delta \in (0,1)$: 
\begin{eqnarray*}
\operatorname*{Im} A_N(u,u) &=& \|u\|^2_{DG} \qquad \forall u \in V_N, \\
|A_N(u,v)| &\leq & C \|u\|_{DG,+} \|v\|_{DG}\qquad \forall u,v \in V_N 
\end{eqnarray*}
Let the solution of $u$ of (\ref{eq:problem}) (with $f = 0$) satisfy 
$u \in H^{3/2+\varepsilon}(\Omega)$ for some $\varepsilon > 0$. Then, 
by consistency of $A_N$, 
the solution $u_N \in V_N$ of (\ref{eq:DG-formulation}) satisfies 
the following quasioptimality estimate for some $C > 0$ independent of $k$: 
\begin{equation}
\label{eq:DG-quasi-optimality}
\|u - u_N\|_{DG} \leq C \inf_{v \in V_N} \|u - v\|_{DG,+}. 
\end{equation}
\end{Proposition}
Several comments are in order: 
\begin{enumerate}
\item
The UWVF of \cite{cessenat-despres98} featured quasi-optimality in a residual
type norm. We recall that the UWVF is a DG method for the particular choice 
$\alpha = \beta = \delta = 1/2$. 
\item 
When $V_N$ consists (elementwise) of systems of plane waves or generalized
harmonic polynomials, then the infimum in (\ref{eq:DG-quasi-optimality})
can be estimated using approximation results on the elements by taking appropriate
traces. This is worked out in detail in
\cite{hiptmair-moiola-perugia09a,hiptmair-moiola-perugia11,moiola09}
and earlier in an $h$-version setting in \cite{cessenat-despres98} (see also 
\cite{buffa-monk07}).
\item 
The $\|\cdot\|_{DG}$-norm controls the error on the skeleton ${\mathcal E}$ only. 
The proof of Lemma~\ref{lemma:monk-wang} shows how error estimates in such norms 
can be used to obtain estimates for $\|u  -u_N\|_{L^2(\Omega)}$; 
we refer again to \cite{buffa-monk07} where this worked out for the UWVF and to 
\cite{hiptmair-moiola-perugia09a,hiptmair-moiola-perugia11,moiola09} where 
the case of the plane wave DG is studied. As pointed out in Remark~\ref{rem:monk-wang}, 
quasi-uniformity of the underlying mesh ${\mathcal T}$ is an important ingredient
for the arguments of Lemma~\ref{lemma:monk-wang}. 
\end{enumerate}
It is noteworthy that Proposition~\ref{prop:DG-quasi-optimality} does
not make any assumptions on the mesh size $h$ and the space $V_N$ 
except that it consist of 
piecewise solutions of the homogeneous Helmholtz equation. 
Optimal error estimates are possible in an 
$h$-version setting, where the number of plane waves per element is kept fixed:  
\begin{Proposition}[\cite{gittelson-hiptmair-perugia09}]
Let $\Omega$ be convex. Assume that 
$V_{N,K} = W^{2\mu+1}_{Pw}$ ($\mu \ge 1$ fixed) 
for all $K \in {\mathcal T}$. Assume that $\alpha$ is of the form 
$\alpha = {\tt a}/(kh)$ and that $\beta > 0$, $\delta \in (0,1/2)$. 
Then there exist ${\tt a}_0$, $c_0$, $C > 0$ (all independent of $h$ and $k$) 
such that if ${\tt a} \ge {\tt a}_0$ and $k^2 h \leq c_0$, then following 
error bound is true: 
$$
\|u - u_N\|_{1,DG} \leq C \inf_{v \in V_N} \|u - v\|_{1,DG,+}; 
$$
here, $\|\cdot\|_{1,DG}$ and $\|\cdot\|_{1,DG,+}$ are given by 
$\|v\|^2_{1,DG}:= \sum_{K \in {\mathcal T}} |v|^2_{H^1(K)} + k^2 \|v\|^2_{L^2(K)} + \|v\|^2_{DG}$
and 
$\|v\|^2_{1,DG,+}:= \sum_{K \in {\mathcal T}} |v|^2_{H^1(K)} + k^2 \|v\|^2_{L^2(K)} + \|v\|^2_{DG,+}$. 
\end{Proposition}
\begin{proof}
The proof follows by inspection of the procedure in \cite[Sec.~5]{gittelson-hiptmair-perugia09}
and is stated in \cite[Props.~{4.2}, {4.3}]{moiola09}. The essential ingredients
of the proof are: (a) inverse estimates for systems of plane waves that have been made in 
available in \cite{gittelson-hiptmair-perugia09} so that techniques of standard DG methods 
can be used to treat $A_N$; (b) use of duality arguments as in 
Lemma~\ref{lemma:abstract-quasi-optimality} to treat the $L^2$-norm of the error; 
(c) the fact that in an $h$-version setting, 
plane waves have some approximation power for arbitrary functions in $H^2$
(this is analogous to Lemma~\ref{lemma:PUM-approximation}). 
\qed
\end{proof}
\iftechreport
\section{Remarks on 1D}
\label{sec:1D}
The 1D situation is rather special in that pollution can be completely 
eliminated; the underlying reason is that the space of solutions of 
the homogeneous Helmholtz equation is finite-dimensional (two dimensional, 
in fact). We illustrate this for the following model problem: 
\begin{equation}
\label{eq:model-problem-1d}
-u^{\prime\prime} - k^2 u  =f \quad \mbox{ in $\Omega = (0,1)$}, 
\qquad u(0) = 0, \qquad 
u^\prime(1) - \bi  k u(1) = 0. 
\end{equation}
Let ${\mathcal T}$ be a mesh on $\Omega$ with nodes $0 = x_0 < x_1 <\cdots < x_N = 1$. 
We assume that the mesh size $h:= \max_i (x_{i+1} - x_i)$ satisfies $kh < \pi$. 
For each node $x_i$, let $\psi_i \in H^1(\Omega)$ be defined by the conditions
$$
\psi_i(x_j) = \delta_{ij}, 
\qquad (-\psi_i^{\prime\prime} - k^2 \psi_i)|_K = 0 
\qquad \forall K \in {\mathcal T} 
$$
and let $V_N^{opt} = \operatorname*{span} \{\psi_i\colon i=1,\ldots,N\}$.
Thus, the functions $\psi_i$ are piecewise solutions of the 
homogeneous Helmholtz equation. 
The Galerkin method based on $V_N^{opt}$ is:  
\begin{equation}
\label{eq:pollution-free-1D-FEM}
\mbox{ Find $u_N \in V_N^{opt}$ s.t.} 
\quad \int_\Omega u_N^\prime \overline{v}^\prime - k^2 u_N \overline{v} 
- \bi k u_N(1) \overline{v}(1)= \int_\Omega f \overline{v} 
\qquad \forall v \in V_N^{opt}. 
\end{equation}
The Galerkin method based on $V_N^{opt}$ is nodally exact: 
\begin{Lemma}
\label{lemma:1d-optimal}
There exist constants $C_1$, $C_2 > 0$ independent of $k$ such that 
the following is true for $kh \leq C_1$: 
\begin{enumerate}[(i)]
\item 
\label{item:lemma:1d-optimal-i}
The functions $\psi_i$ are well-defined. 
\item
\label{item:lemma:1d-optimal-ii}
The method (\ref{eq:pollution-free-1D-FEM}) is nodally exact.
\item
\label{item:lemma:1d-optimal-iii}
For $f \in L^2(\Omega)$ there holds 
$
\|u - u_N\|_{1,k,\Omega} \leq C_2 (hk) \|f\|_{L^2(\Omega)}. 
$
\end{enumerate}
\end{Lemma}
\begin{proof}
Elementary considerations show that for $kh < \pi$, the functions 
$\psi_i$ are well-defined. 

The most interesting feature of Lemma~\ref{lemma:1d-optimal} is
the nodal exactness. To that end, we note that the 
Green's function for (\ref{eq:model-problem-1d}) is 
$$
G(x,y) = \frac{1}{k} 
\begin{cases}
\sin kx \, e^{\bi k y} & 0 <  x < y \\
\sin ky \, e^{\bi k x} & y <  x < 1. 
\end{cases}
$$
Let $e \in H^1(\Omega)$
satisfy $e(0) = 0$ and the Galerkin orthogonality condition
\begin{equation}
\label{eq:1D-galerkin-orthogonality}
C(e,v):= \int_\Omega e^\prime \overline{v}^\prime - k^2 e \overline{v} 
- \bi k e(1) \overline{v}(1) = 0
\qquad \forall v \in V^{opt}_N. 
\end{equation}
The key observation is that for each $x_i$, $i=1,\ldots,N$, 
the function $v_i:= G(\cdot,x_i) \in V^{opt}_N$ since is a solution
of the homogeneous Helmholtz equation on $(0,x_i) \cup (x_i,1)$, 
it satisfies $G(0,x_i) = 0$. Furthermore, we have 
$v_i^\prime(x_N) - \bi k v_i(x_N) = 0$. 
Hence, we get from the Galerkin orthogonality (\ref{eq:1D-galerkin-orthogonality})
by an integration by parts: 
\begin{eqnarray*}
0 &=& 
\int_\Omega e^\prime \overline{v}_i^\prime - k^2 e \overline{v}_i 
+ \bi k e(1) \overline{v}_i(1) \\
&=& 
\int_{0}^{x_i} e^\prime \overline{v}_i^\prime - k^2 e \overline{v}_i 
+ \int_{x_i}^{x_N} e^\prime \overline{v}_i^\prime - k^2 e \overline{v}_i 
- \bi k e(1) \overline{v}_i(1) \\
&=&
\int_{0}^{x_i} e (-\overline{v}_i^{\prime\prime} - k^2 \overline{v}_i) 
+ e(x_i) \overline{v}_i^\prime(x_i) 
- e(x_i) \overline{v}_i^\prime(x_i)  \\
&& \quad \mbox{}
+ \int_{x_i}^{1} e (-\overline{v}_i^{\prime\prime} - k^2  \overline{v}_i )
+ e(1) \overline{v}_i^\prime(1)
- \bi k e(1) \overline{v}_i(1) \\
&=& e(x_i) [\overline{v}^\prime](x_i) 
+ e(1) \overline{v}_i^\prime(1)
- \bi k e(1) \overline{v}_i(1) \\
&=& e(x_i) [\overline{v}_i^\prime](x_i); 
\end{eqnarray*}
here, we have employed the standard notation for the jump of a piecewise 
smooth function $w$: $[w^\prime](x_i):= 
\lim_{x \rightarrow x_i-} w^\prime(x) -
\lim_{x \rightarrow x_i+} w^\prime(x)$.  Since $[v_i^\prime](x_i) \ne 0$, 
we conclude 
$$
e(x_i) = 0 \qquad \forall i \in \{1,\ldots,N\}. 
$$
Hence, the FEM (\ref{eq:pollution-free-1D-FEM}) is nodally exact. 

The above  argument also shows that any $e \in V^{opt}_N$ satisfying 
(\ref{eq:1D-galerkin-orthogonality}) must satisfy $e(x_i) = 0$ for all 
$i \in \{1,\ldots,N\}$. Hence, by the definition of $V^{opt}_N$ as the 
span of the functions $\psi_i$, we conclude $e \equiv 0$. Thus, the 
kernel of the linear systems described by (\ref{eq:pollution-free-1D-FEM})
is trivial. By the usual dimension argument, we have unique solvability. 

We have the {\sl a priori} bound
\begin{equation}
\label{eq:lemma:1d-optimal-10}
\|u\|_{1,k,\Omega} \leq C \|f\|_{L^2(\Omega)}
\end{equation}
for the solution $u$ of (\ref{eq:model-problem-1d})
(as for the model problem (\ref{eq:problem}), this can be shown
using the test function $v = x u^\prime$; an alternative proof
based on the Green's function and the representation 
$$
u(x) = 
\int_\Omega G(x,y) f(y)\,dy
$$ 
is given in
\cite[Thm.~{4.4}]{ihlenburg98}). If one denotes  by $\varphi_i$ the
classical piecewise linear hat function associated with node $x_i$,
then one has by Taylor expansion
\begin{equation}
\label{eq:lemma:1d-optimal-20}
\|\varphi_i - \psi_i\|_{L^\infty(K)} 
\leq C ( k h_K)^2, 
\qquad 
\|(\varphi_i - \psi_i)^\prime\|_{L^\infty(K)} 
\leq C  k^2 h_K,
\qquad \forall K \in {\mathcal T}, 
\end{equation}
where $h_K =\operatorname*{diam} K$. 
The approximation properties follow easily from the nodal exactness.
Specifically, denoting by $I u$ the classical piecewise linear interpolant
of $u$ and by $\widetilde I u \in V_N^{opt}$ the nodal interpolant
determined by $\widetilde I u(x_i) = u(x_i)$ for all $i \in \{0,\ldots,N\}$,
we have the well-known estimate
$$
\|u - I u\|_{1,k,\Omega} \leq C (kh^2 + h)\|u^{\prime\prime}\|_{L^2(\Omega)}
\leq C kh ( 1 + kh)\|f\|_{L^2(\Omega)}, 
$$
where we used the differential equation and
the bound (\ref{eq:lemma:1d-optimal-10})
to estimate $\|u^{\prime\prime}\|_{L^2(\Omega)} 
\leq \|f\|_{L^2(\Omega)} + k^2 \|u\|_{L^2(\Omega)} 
\leq C k \|f\|_{L^2(\Omega)}$. Next, we estimate the difference
$Iu - \widetilde Iu$. The multiplicative trace inequality takes the form
\begin{equation}
\label{eq:lemma:1d-optimal-30}
h_K \|w\|^2_{L^\infty(K)} \leq C \left[ \|w\|^2_{L^2(K)} + h_K \|w\|_{L^2(K)} \|w^\prime\|_{L^2(K)}
\right] \qquad \forall w \in H^1(K). 
\end{equation}
Hence, the estimates (\ref{eq:lemma:1d-optimal-20}), (\ref{eq:lemma:1d-optimal-30})
imply
\begin{eqnarray*}
\|Iu - \widetilde I u\|^2_{1,k,\Omega} &\leq& 
C \sum_{K \in {\mathcal T}} k^2 (kh_K)^2 (1 + (kh_K)^2) h_K\|u\|^2_{L^\infty(K)} \\
& \leq & C (kh)^2 (1 + (kh)^2) 
\left[ 
k^2 \|u\|^2_{L^2(\Omega)} + k^2 h^2 \|u^\prime\|^2_{L^2(\Omega)}
\right]
\\
& \leq &
C (kh)^2 (1 + (kh)^2)^2 \|u\|^2_{1,k,\Omega}. 
\end{eqnarray*}
An appeal to (\ref{eq:lemma:1d-optimal-10}) concludes the argument.
\qed
\end{proof}
Several comments are in order concerning the stability of the method: 
\begin{enumerate}
\item 
In the 1D situation, the good stability properties of 
high order Galerkin FEM can alternatively be understood in light 
of Lemma~\ref{lemma:1d-optimal}:
Applying the Galerkin method to a classical high order method and then 
condensing out the degrees of freedom corresponding to internal
shape functions (``bubbles''), leads to a linear system that is 
identical to the one obtained by using shape functions 
$\psi_i^p$, $i=0,\ldots,N$, that satisfy 
$
\psi_i^p(x_j) = \delta_{ij}
$
and additionally 
$$
\int_K (\psi_i^p)^\prime  \overline{v}^\prime - k^2 \psi_i^p \overline{v} = 0 
\qquad \forall v \in H^1_0(K) \cap {\mathcal P}_p
$$
Since on a fixed mesh ${\mathcal T}$, we have 
$\lim_{p\rightarrow \infty} \psi_i^p = \psi_i$, better stability properties
of higher order methods may reasonably be hoped for. 
\item 
The system matrix of the Galerkin FEM based on the space $V_N^{opt}$ is a tridiagonal
matrix. The same matrix can also be obtained in different ways. Consider, for example, 
the sesquilinear form 
\begin{equation}
\label{eq:stabilized-bilinear-form}
B(u,v):= \int_\Omega u^\prime \overline{v}^\prime - k u \overline{v} - \bi k u(1) \overline{v}(1) 
+ \sum_{K \in {\mathcal T}} \tau_K \int_\Omega L_k u L_k \overline{v}, 
\end{equation}
where $L_k = - \frac{d^2}{dx^2} - k^2$. For a suitable choice of the parameters $\tau_K$ 
in dependence on $k$ and $h_K$, the system matrix resulting from this $B$ using the classical piecewise
linear hat functions leads to the same matrix as the Galerkin method based on the shape
functions $\psi_i$, $i=1,\ldots,N$. In 1D, it is therefore possible to design nodally exact
methods based on the stabilization techniques in the from (\ref{eq:stabilized-bilinear-form}). 
In \cite{babuska-sauter00}, a nodally exact method is derived using other techniques. 

\end{enumerate}

\fi 

{\bf Acknowledgement:} Financial support by the 
{\sl Vienna Science and Technology Fund} (WWTF) is gratefully acknowledged.
\iftechreport
\section*{Appendix}
\addcontentsline{toc}{section}{Appendix}
\setcounter{equation}{0}
\setcounter{section}{0}
\setcounter{Theorem}{0}
\def\theequation{A.\arabic{equation}}
\def\theTheorem{A.\arabic{Theorem}}

\input{appendix}
\fi 
%
%
\bibliographystyle{plain}
\bibliography{nummech}
\end{document}

%% file: appendix.tex
For the reference triangle $\widehat K:= \{(x,y)\colon 0 < x < 1, 0 < y < 1-x\}$
and $\beta \in [0,1)$ the following two lemmas require the spaces 
$H^{1,1}_{\beta}(\widehat K)$, $H^{2,2}_\beta(\widehat K)$ as well as 
the Besov spaces $B^s_{2,\infty}(\widehat K)$. The spaces $B^s_{2,\infty}(\widehat K)$
are defined by interpolation using the $K$-functional (see, e.g., 
\cite[Chap.~{12}]{brenner-scott94}). For $m \in \{1,2\}$, the spaces 
$H^{m,m}_\beta(\widehat K)$ are determined by the norm 
$
\|u\|^2_{H^{m,m}_\beta(\widehat K)} := \|u\|^2_{H^{m-1}(\widehat K)} + \|r^\beta \nabla^m u\|^2_{L^2(\widehat K)},
$
where $r$ denotes the distance from the origin.
\begin{Lemma}
\label{lemma:besov-space-embedding}
Let $\widehat K$ be the  
reference triangle. Let $\beta \in [0,1)$. Then the embeddings 
$H^{2,2}_{\beta}(\widehat K) \subset B^{2-\beta}_{2,\infty}(\widehat K)$ and 
$H^{1,1}_{\beta}(\widehat K) \subset B^{1-\beta}_{2,\infty}(\widehat K)$  
are continuous. 
The embeddings $H^{2,2}_{\beta}(\widehat K) \subset 
H^{2-\beta-\varepsilon}(\widehat K)$ and 
$H^{1,1}_{\beta}(\widehat K) \subset H^{1-\beta-\varepsilon}(\widehat K)$ are compact for 
all $\varepsilon > 0$. 
\end{Lemma}
\begin{proof}
Since the case $\beta = 0$ corresponds to classical Sobolev spaces, we 
restrict our attention here to the situation $\beta \in (0,1)$. 
The argument follows ideas presented in \cite[Thm.~{2.1}]{babuska-osborn91}
and \cite{babuska-kellogg79a}. 
We start with the following two Hardy inequalities for sufficiently
smooth functions
\begin{eqnarray}
\label{eq:lemma:besov-space-embedding-1}
\|r^{\beta - 1}\nabla u\|_{L^2(\widehat K)} &\leq& C \|u\|_{H^{2,2}_{\beta}(\widehat K)},\\ 
\label{eq:lemma:besov-space-embedding-2}
\|r^{\beta - 2} (u - u(0))\|_{L^2(\widehat K)} &\leq & C \|u\|_{H^{2,2}_{\beta}(\widehat K)}; 
\end{eqnarray}
here, 
(\ref{eq:lemma:besov-space-embedding-1}) is shown, for example, in \cite[Lemma~{A.1.7}]{melenk02}
and (\ref{eq:lemma:besov-space-embedding-2}) follows from combining \cite[Lemma~{4.2}]{babuska-kellogg79a}
with (\ref{eq:lemma:besov-space-embedding-1}). 
Noting that \cite[(2.2)]{babuska-kellogg79a} states the continuous
embedding $H^{2,2}_{\beta}(\widehat K) \subset C(\overline{\widehat K})$, we have that $u(0)$ 
in (\ref{eq:lemma:besov-space-embedding-2}) is indeed well-defined. 

We employ the real method of interpolation and write $B^{2-\beta}_{2,\infty} = (L^2,H^2)_{1-\beta/2,\infty}$. 
Our method of proof consists in showing that for $\theta = 1-  \beta/2$ we have 
$$
\sup_{t \in (0,1)} t^{-\theta} K(t,\widetilde u) \leq C \|u\|_{H^{2,2}_{\beta}(\widehat K)}, 
\qquad \widetilde u:= u - u(0), 
$$
for some $C > 0$ independent of $u$. To that end, 
we proceed as in the proof of \cite[Lemma~{2.1}]{babuska-kellogg79a}. 
For every $\delta > 0$, let $\chi_\delta \in C^\infty_0(\BbbR^2)$ with 
$\chi \equiv 1$ on $B_{\delta/2}(0)$ and 
$\operatorname*{supp} \subset \chi_\delta B_\delta(0)$ as well as 
$\|\nabla^j \chi_\delta\|_{L^\infty(\BbbR^2)} \leq C \delta^{-j}$, $j \in \{0,1,2\}$. We define the splitting 
$$
\widetilde u = \chi_\delta \widetilde u + (1 - \chi_\delta) \widetilde u 
=: u_1 + u_2
$$
Then from (\ref{eq:lemma:besov-space-embedding-1}) and (\ref{eq:lemma:besov-space-embedding-2})
\begin{align*}
&\|\chi_\delta \widetilde u\|_{L^2(\widehat K)}  \leq 
C \|\widetilde u\|_{L^2(\widehat K\cap B_{\delta}(0))} 
\leq \delta^{2-\beta} \|r^{\beta - 2} \widetilde u\|_{L^2(\widehat K)} 
\leq C \delta^{2-\beta} \|u\|_{H^{2,2}_{\beta}(\widehat K)},  \\
& |(1 - \chi_\delta) \widetilde u|_{H^2(\widehat K)} \leq \\
&
C \delta^{-2} \|\widetilde u\|_{L^2((\widehat K \cap B_{\delta}(0)) \setminus B_{\delta/2}(0))} + 
C \delta^{-1} \|\nabla \widetilde u\|_{L^2((\widehat K \cap B_{\delta}(0) )\setminus B_{\delta/2}(0))} + 
C \|\nabla^2 \widetilde u\|_{L^2(\widehat K \setminus B_{\delta/2}(0))} \\
&\leq C \delta^{-2+2-\beta} \|r^{\beta-2}\widetilde u\|_{L^2(\widehat K)}+
C \delta^{-1+1-\beta} \|r^{\beta-1} \nabla \widetilde u\|_{L^2(\widehat K)} + 
C \delta^{-\beta}\|r^\beta \nabla^2 \widetilde u\|_{L^2(\widehat K)} \\
&\leq C \delta^{-\beta} \|u\|_{H^{2,2}_{\beta}(\widehat K)}. 
\end{align*}
From this, we can infer for any $\delta \in (0,1)$
$$
K(t,\widetilde u) \leq \|u_1\|_{L^2(\widehat K)} + t \|u_2\|_{H^2(\widehat K)}
\leq C \|u\|_{H^{2,2}_{\beta}(\widehat K)}\left[ \delta^{2-\beta} + t \delta^{-\beta}\right].
$$
Selecting $\delta = t^{1/2}$ gives 
$K(t,\widetilde u) \leq C t^{1-\beta/2} \|u\|_{H^{2,2}_{\beta}(\widehat K)}$. 
Finally, the compactness assertions of the embeddings follows from the compactness
of the embeddings $B^{s}_{2,\infty} \subset B^{s^\prime}_{2,2} = H^{s^\prime}$ for 
$s^\prime < s$. 
\qed
\end{proof}
\begin{Lemma}
\label{lemma:trace-inequality-H11beta}
Let $\beta \in [0,1/2)$ and $\widehat K$ be the reference triangle. Then 
there exists $C > 0$ such that for all $u \in H^{1,1}_\beta(\widehat K)$
there holds 
$
\|u\|_{L^2(\partial \widehat K)} \leq C \left[\|u\|_{L^2(\widehat K)} + 
\|r^\beta \nabla u\|_{L^2(\widehat K)}\right] . 
$
\end{Lemma}
\begin{proof}
For each $s > 1/2$, we have the 
inequality $\|u\|_{L^2(\partial\widehat K)} \leq C_s \|u\|_{H^s(\widehat K)}$. 
From the embedding $H^{1,1}_\beta(\widehat K) \subset H^{1-\beta}(\widehat K)$
of Lemma~\ref{lemma:besov-space-embedding}, we then get 
$
\|u\|_{L^2(\partial\widehat K)} \leq C\|u\|_{H^s(\widehat K)} 
\leq C \left[ \|u\|_{L^2(\widehat K)} + \|r^\beta \nabla u\|_{L^2(\widehat K)}\right]. 
$
\qed
\end{proof}
\begin{Lemma}
\label{lemma:singular-fcts}
Let $\beta \in [0,1)$ and $\Omega\subset\BbbR^2$ be a finite sector 
with apex at the origin. Let $u \in C^\infty(\Omega)$ satisfy 
\begin{eqnarray*}
\|\Phi_{n,\beta,1} \nabla^{n+2} u\|_{L^2(\Omega)} &\leq& C_u \gamma_u n! 
\qquad \forall n \in \BbbN_0.
\end{eqnarray*}
Then, for $k \ge k_0 > 0$, there exist constants $C$, $\gamma > 0$
(depending only on $\beta$, $\Omega$, $\gamma_u$, and $k_0$) such that 
\begin{eqnarray*}
\|\Phi_{n,\beta,k} \nabla^{n+2} u\|_{L^2(\Omega)} &\leq& 
C C_u k^{-(2-\beta)} \gamma^n \max\{n,k\}^{n+2}
\qquad \forall n \in \BbbN_0.
\end{eqnarray*}
\end{Lemma}
\begin{proof}
Lemma~\ref{lemma:properties-of-weights} yields 
$$
\frac{1}{\max\{n,k\}^{n+2}} \Phi_{n,\beta,k}(x) 
\leq C k^{-(2-\beta)} \gamma^n \frac{1}{n!} \Phi_{n,\beta,1}(x) 
\qquad\forall x \in \Omega, 
$$
where $C$, $\gamma > 0$ are independent of $n$ and $k$. The result
now follows.
\qed
\end{proof}
\begin{Lemma}
\label{lemma:properties-of-weights}
Let $\beta \in [0,1)$. Then for $0 < r < R$ and all $n \in \BbbN_0$
$$
\min\left( 1, \frac{r}{\min\left\{1, \frac{n+1}{k+1}\right\}}\right)^{n+\beta}
\frac{1}{\max\{n,k\}^{n+2}}
\leq C k^{-(2-\beta)} \gamma^n r^{n+\beta} \frac{1}{n^{n+2}}
$$
\end{Lemma}
\begin{proof}
We denote the left-hand side by $lhs$ and consider several cases. 

{\em case 1: $n\leq k$ and $r(k+1) \leq n+1$}: 
\begin{eqnarray*}
lhs &=& 
\left(\frac{(k+1)r}{n+1}\right)^{n+\beta}\frac{1}{k^{n+2}}  \\
& = & 
r^{n+\beta} \frac{1}{n^{n+2}} \left(\frac{n}{n+1}\right)^{n+2}(n+1)^{2-\beta}
\left(\frac{k+1}{k}\right)^{n+2} (k+1)^{-(2-\beta)}\\
&\leq& C \gamma^n k^{-(2-\beta)} r^{n+\beta} \frac{1}{n^{n+2}}
\end{eqnarray*}
for suitable $C$, $\gamma > 0$ if we assume that $k \ge k_0 > 0$. 
\newline 
{\em case 2: $n\leq k$ and $r(k+1) >  n+1$}:
\begin{eqnarray*}
lhs &=& \frac{1}{k^{n+2}}  = 
\frac{1}{k^{2-\beta}} \frac{1}{k^{n+\beta}}
= 
\frac{1}{k^{2-\beta}} \left(\frac{k+1}{k}\right)^{n+\beta} \frac{1}{(k+1)^{n+\beta}}\\
&\leq& 
\frac{1}{k^{2-\beta}} \left(\frac{k+1}{k}\right)^{n+\beta} \left(\frac{r}{n+1}\right)^{n+\beta}\\
&\leq & 
C \gamma^n k^{-(2-\beta)} r^{n+\beta} \frac{1}{n^{n+2}}
\end{eqnarray*}
for suitable $C$, $\gamma > 0$. 
\newline 
{\em case 3: $n>k$}: Then, for $0 < r < R$
\begin{eqnarray*} 
lhs &=& \left(\min\{1,r\}\right)^{n+\beta} \frac{1}{n^{n+2}}
\leq   r^{n+\beta} \frac{1}{n^{n+2}} 
\leq k^{-(2-\beta)} r^{n+\beta} \frac{1}{n^{n+2}} n^{2-\beta}\\
& \leq & 
C k^{-(2-\beta)} r^{n+\beta} \frac{1}{n^{n+2}} \gamma^n 
\end{eqnarray*} 
for suitable $C$, $\gamma > 0$. 
\qed
\end{proof}
\section*{Further results and proofs} 
\begin{Example}
{\rm 
\label{ex:1d-pollution-L2}
In Example~\ref{ex:1d-pollution}, we studied the convergence behavior 
of the $h$-FEM in the $H^1(\Omega)$-seminorm. 
In Fig.~\ref{fig:1D-example-L2} we present the corresponding
results for the convergence in the $L^2(\Omega)$-norm 
by plotting $\|u - u_N\|_{L^2(\Omega)}/\|u\|_{L^2(\Omega)}$ vs. the 
number of degrees of freedom per wavelength $N_\lambda$. 
For $p = 1$, we observe 
$$
\frac{\|u - u_N\|_{L^2}}{\|u\|_{L^2}} \approx C k N_\lambda^{-2}, 
\qquad N_\lambda \rightarrow \infty,
$$ 
which is
in agreement with the analysis given in \cite[Sec.~{4.6.4}]{ihlenburg98}. 
The cases $p > 1$ seem to behave differently as we observe 
$$
\frac{\|u  -u_N\|_{L^2}}{\|u\|_{L^2}} \approx C N_\lambda^{-(p+1)}, 
\qquad N_\lambda \rightarrow \infty
$$
\begin{figure}
\includegraphics[width=0.5\textwidth]{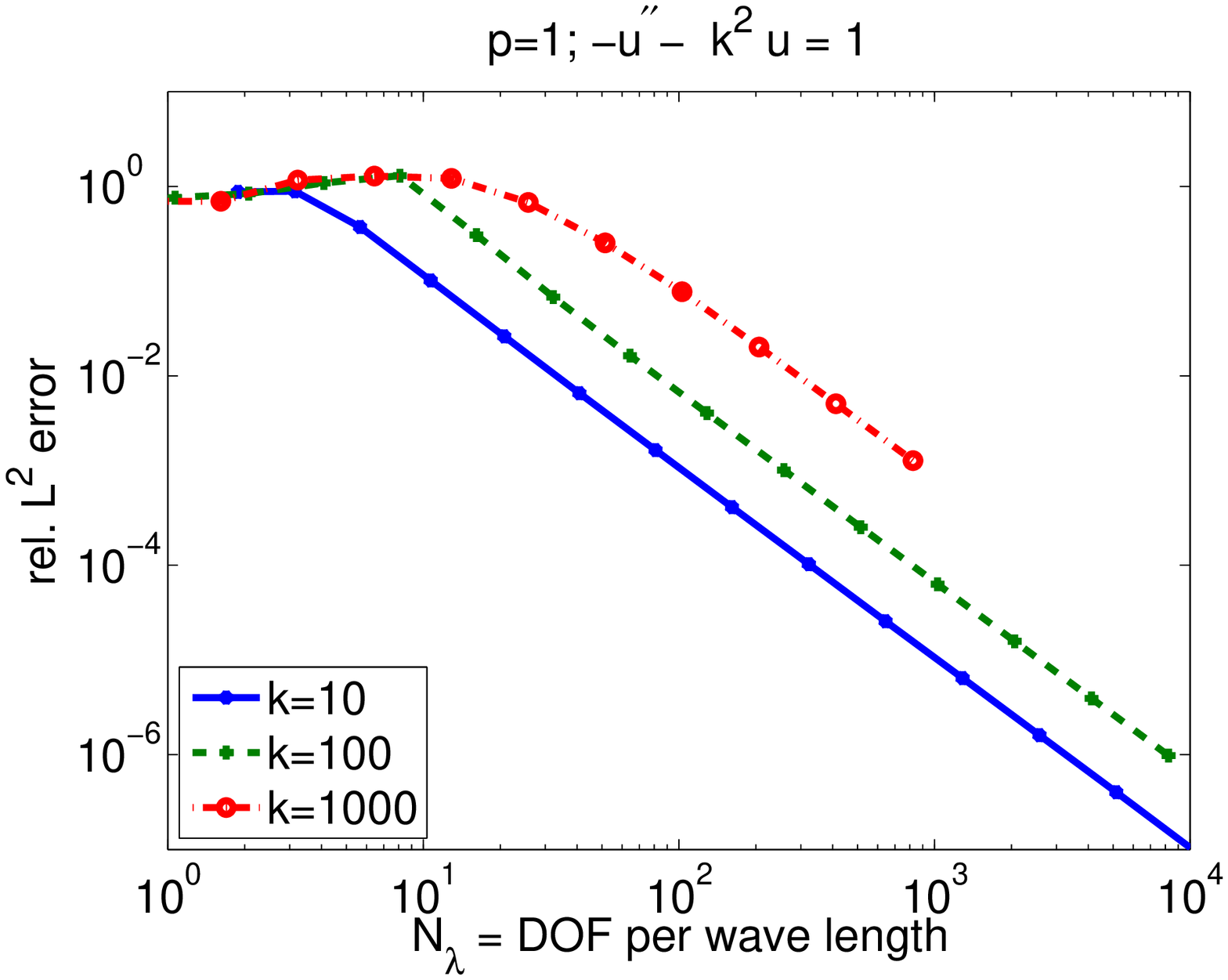}
\includegraphics[width=0.5\textwidth]{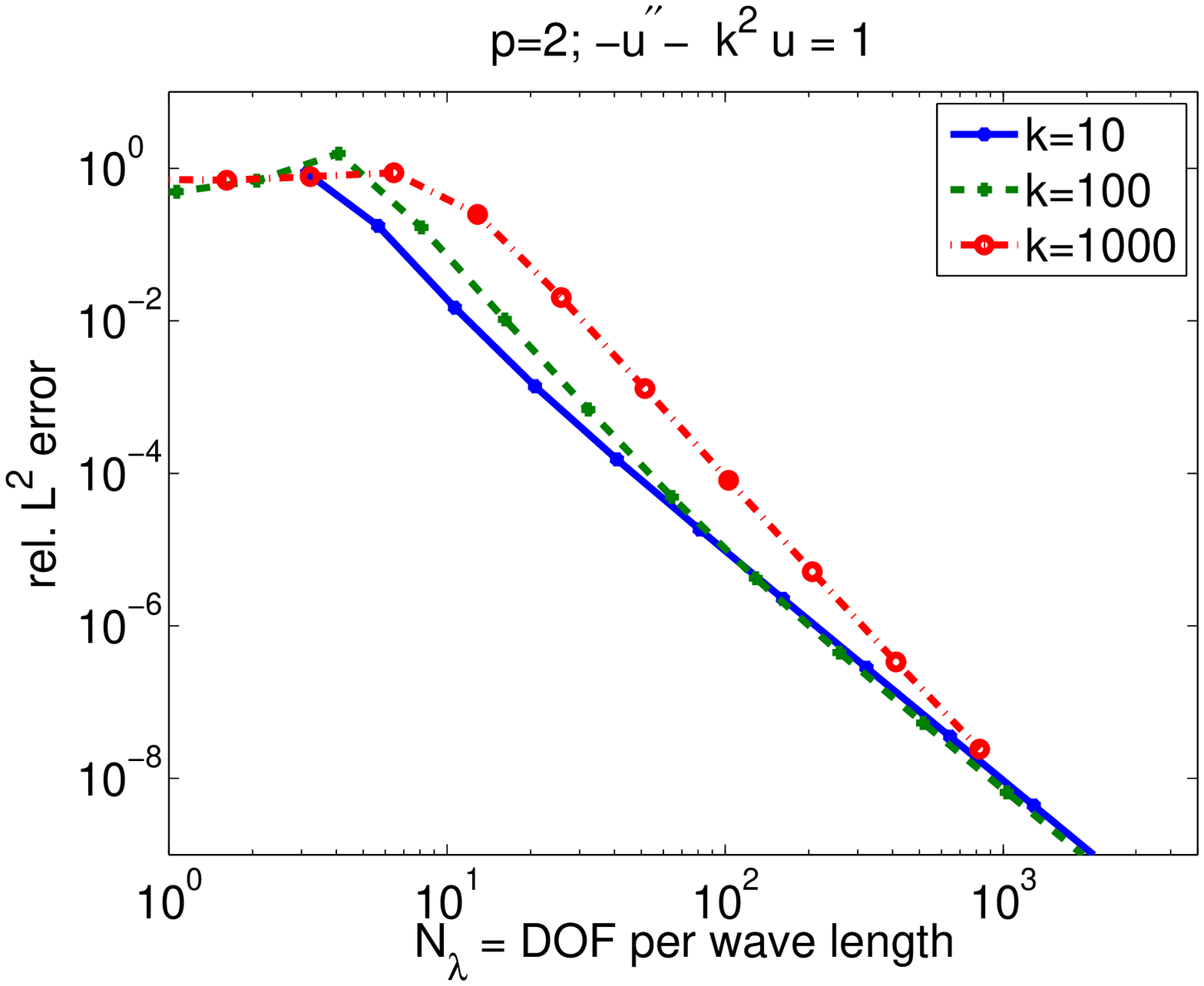}
\includegraphics[width=0.5\textwidth]{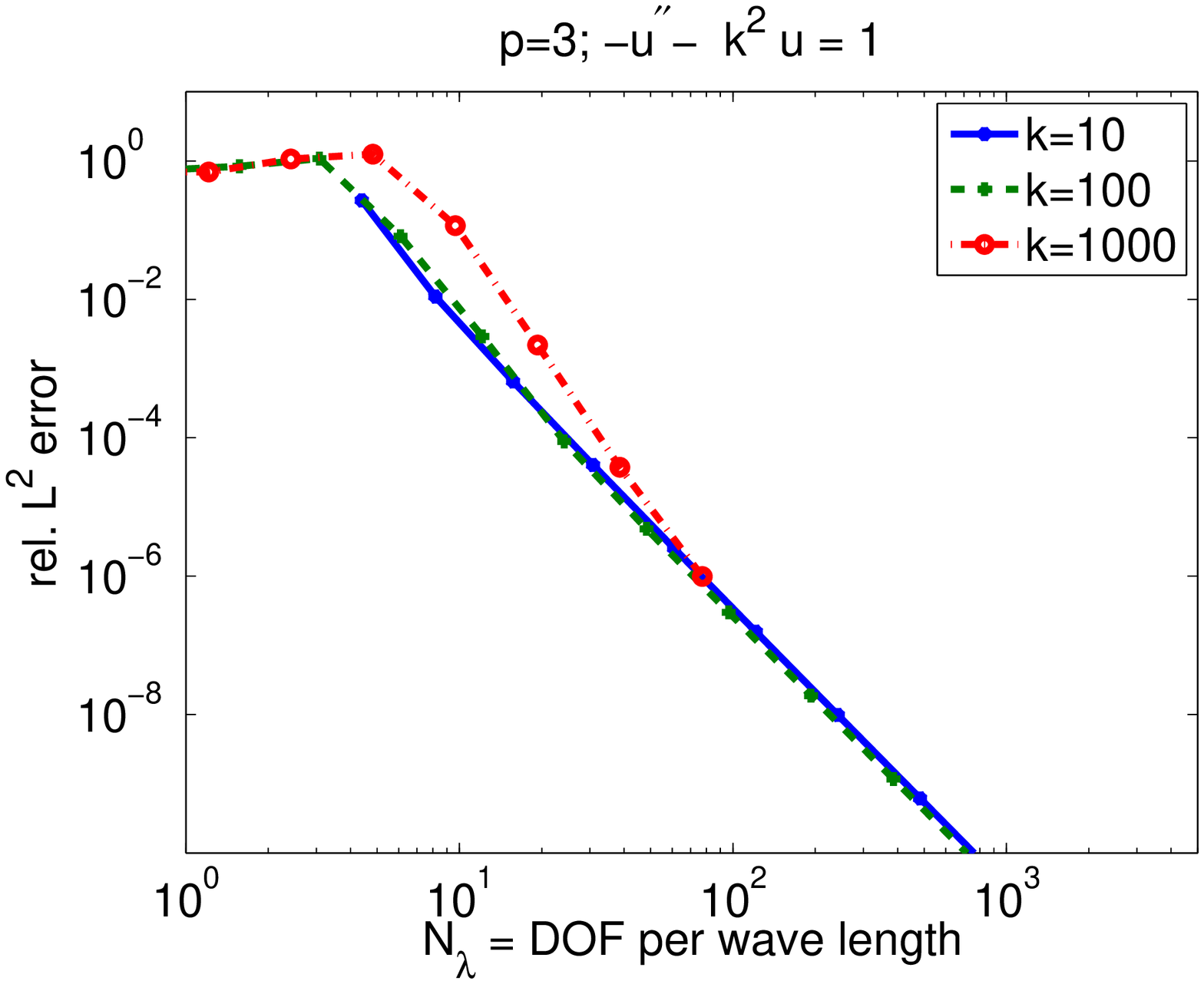}
\includegraphics[width=0.5\textwidth]{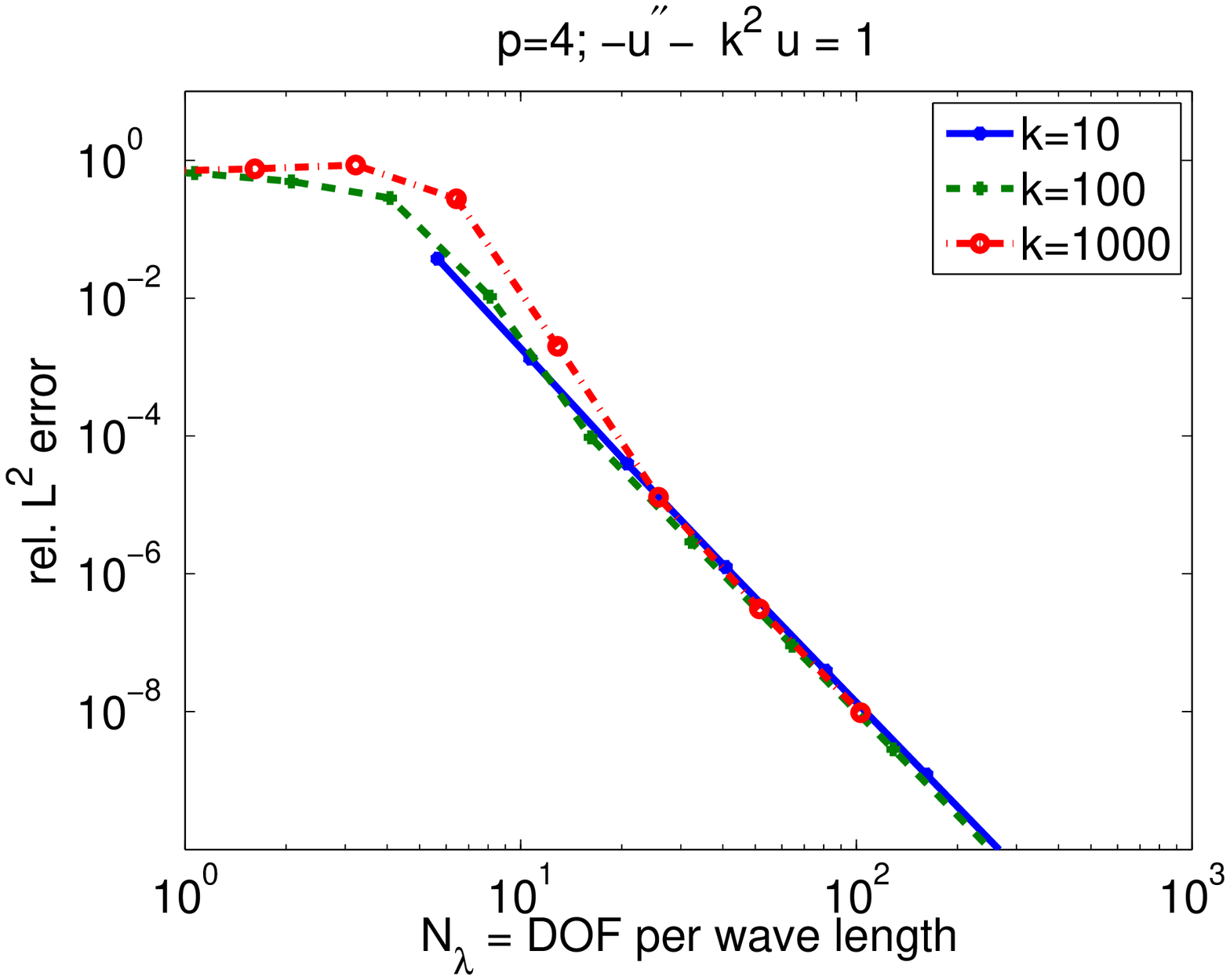}
\caption{\label{fig:1D-example-L2} Performance of $h$-FEM for (\ref{eq:1D-model}).
Top: $p=1$, $p=2$. Bottom: $p=3$, $p=4$
(cf.~Examples~\ref{ex:1d-pollution-L2}, \ref{ex:1d-pollution}).}
\end{figure}
\eremk
}
\end{Example}

\begin{numberedproof}{Remark~\ref{rem:pum-approximation}}
By interpolation 
using the $K$-functional 
we can write $H^{1+\theta} = (H^1,H^2)_{\theta,2}$ for $\theta \in (0,1)$. 
Hence, every $u \in H^{1+\theta}(\Omega)$ can be decomposed as 
$
u = u_1 + u_2
$
with 
\begin{equation}
\label{eq:rem-pum-approximation-10}
\|u_1\|_{H^1(\Omega)} \leq t^{\theta} \|u\|_{H^{1+\theta}(\Omega)}, 
\qquad 
\|u_2\|_{H^2(\Omega)} \leq t^{\theta-1} \|u\|_{H^{1+\theta}(\Omega)}, 
\end{equation}
where $t > 0$ is arbitrary. The proof of Lemma~\ref{lemma:PUM-approximation}
shows that $u_1$ and $u_2$ can be approximated from $V_N$ as follows: 
\begin{eqnarray*}
\inf_{v \in V_N} \|u_2 - v\|_{L^2(\Omega)} + h \|\nabla (u_2 - v)\|_{L^2(\Omega)} &\leq& 
h^2 \|u_2\|_{H^2(\Omega)} + (kh)^2 \|u_2 \|_{L^2(\Omega)} \\
\inf_{v \in V_N} \|u_1 - v\|_{L^2(\Omega)} + h \|\nabla (u_1 - v)\|_{L^2(\Omega)} &\leq& 
h \|u_1\|_{H^1(\Omega)} + (kh)^2 \|u_1 \|_{L^2(\Omega)}.  
\end{eqnarray*}
Using $t = h$ in (\ref{eq:rem-pum-approximation-10}) we therefore get 
\begin{eqnarray*}
\inf_{v \in V_N} \|u - v\|_{L^2(\Omega)} + h \|\nabla (u - v)\|_{L^2(\Omega)} &\leq& 
h^{1+\theta} \|u\|_{H^{1+\theta}(\Omega)} + (kh)^2 \left[ \|u_1\|_{L^2(\Omega)} + \|u_2 \|_{L^2(\Omega)}\right].
\end{eqnarray*} 
The decomposition $u = u_1 + u_2$ and the 
triangle inequality yield $\|u_1\|_{L^2(\Omega)} + \|u_2\|_{L^2(\Omega)} \leq 
\|u\|_{L^2(\Omega)} + 2 \|u_1\|_{L^2(\Omega)} \leq 
\|u\|_{L^2(\Omega)} + 2 \|u_1\|_{H^1(\Omega)} \leq 
\|u\|_{L^2(\Omega)} + 2 h^{\theta} \|u\|_{H^{1+\theta}(\Omega)}$. Combining these estimates, we obtain  
\begin{eqnarray*}
\inf_{v \in V_N} \|u - v\|_{L^2(\Omega)} + h \|\nabla (u_2 - v)\|_{L^2(\Omega)} &\leq& 
\left(h^{1+\theta} + (kh)^2 h^\theta\right)\|u\|_{H^{1+\theta}(\Omega)} + (kh)^2 \|u\|_{L^2(\Omega)}, 
\end{eqnarray*} 
which concludes the proof. 
\end{numberedproof}

\begin{Lemma}
\label{lemma:approximation-of-H^22beta-fcts}
Let $\beta \in [0,1)$. Then, 
for every $p \in \BbbN$ there exists a linear operator 
$\pi_p: H^{2,2}_{\beta}(\widehat K) \rightarrow {\mathcal P}_p$ that admits 
an ``element-by-element construction''
in the sense of \cite[Def.~{5.3}]{melenk-sauter10} with the following approximation property: 
$$
p\|u - \pi_p u\|_{L^2(\widehat K)} + \|u - \pi_p u\|_{H^1(\widehat K)} \leq 
C p^{-(1-\beta)} \|r^\beta \nabla^2 u\|_{L^2(\widehat K)}, 
$$
where $C > 0$ is independent of $p$ and $u$. 
\end{Lemma}
\begin{proof}
Inspection of the proof of \cite[Thm.~B.4]{melenk-sauter10} shows that the operator 
$\pi_p$ constructed there does in fact not depend on the regularity parameter $s > 1$. It has (as stated
in \cite[Thm.~{B.4}]{melenk-sauter10}), the approximation property  
\begin{equation}
\label{eq:lemma:approximation-of-H^22beta-fcts-1}
p\|u - \pi_p u\|_{L^2(\widehat K)} + \|u - \pi_p u\|_{H^1(\widehat K)} \leq 
C p^{-(s-1)} |u|_{H^s(\widehat K)} \qquad \forall u \in H^s(\widehat K), 
\end{equation}
if $p \ge s-1$. 
Upon writing the Besov space $B^{s}_{2,\infty}$ as an interpolation space  
$B^{s}_{2,\infty} = (H^2(\widehat K),H^{1}(\widehat K))_{s-1,\infty}$ for $s \in (1,2)$, 
we can infer for $s \in (1,2)$ from (\ref{eq:lemma:approximation-of-H^22beta-fcts-1}) 
the slightly stronger statement 
\begin{equation}
\label{eq:lemma:approximation-of-H^22beta-fcts-2}
p\|u - \pi_p u\|_{L^2(\widehat K)} + \|u - \pi_p u\|_{H^1(\widehat K)} \leq 
C p^{-(s-1)} \|u\|_{B^s_{2,\infty}(\widehat K)} \qquad \forall u \in B^s_{2,\infty}(\widehat K). 
\end{equation}
Appealing to Lemma~\ref{lemma:besov-space-embedding} then yields 
\begin{equation}
\label{eq:lemma:approximation-of-H^22beta-fcts-3}
p\|u - \pi_p u\|_{L^2(\widehat K)} + \|u - \pi_p u\|_{H^1(\widehat K)} \leq 
C p^{-(s-1)} \|u\|_{H^{2,2}_{\beta}(\widehat K)}.  
\end{equation}
We replace the full $H^{2,2}_{\beta}(\widehat K)$ norm by the seminorm in the standard
way by a compactness argument. Since $H^{2,2}_\beta(\widehat K)$ is compactly
embedded in $H^1(\widehat K)$ (see, e.g., \cite[Lemma~{4.19}]{schwab98}) one obtains 
$\inf_{v \in {\mathcal P}_1} \|u - v\|_{H^{2,2}_{\beta}(\widehat K)} \leq 
C \|r^\beta \nabla^2 u\|_{L^2(\widehat K)}$. The proof is completed by noting that 
(\ref{eq:lemma:approximation-of-H^22beta-fcts-1}) implies that 
$\pi_p$ reproduces linear polynomials.  
\qed
\end{proof}